\documentclass[11pt,letterpaper]{article}
\RequirePackage[backend=biber,bibstyle=authoryear,citestyle=authoryear,
natbib=true,sorting=nyt,sortcites=true,maxnames=10,minnames=10,
hyperref=true,abbreviate=false,date=long,isbn=true,eprint=true,
uniquename=init,giveninits=true]{biblatex}
\usepackage{hyphenat,graphicx}
\usepackage[colorlinks,linktocpage,breaklinks]{hyperref}
\makeatletter
\def\@filecolor{blue}
\def\@linkcolor{blue}
\def\@citecolor{blue}
\def\@urlcolor{blue}
\let\@old@citep\citep
\let\@old@citet\citet
\let\@old@citeauthor\citeauthor
\def\citep{\@old@citep*}
\def\citet{\@old@citet*}
\def\citeauthor{\@old@citeauthor*}
\let\cite\citep
\makeatother
\usepackage[top=1.5in,bottom=1.2in,left=1.3in,right=1.3in,letterpaper]%
  {geometry}
\makeatletter
\def\ps@headings{%
    \let\@mkboth\@gobbletwo
    \def\@oddhead{\hss\scshape\shorttitle\hss\reset@font\rmfamily\thepage}
    \def\@evenhead{\reset@font\rmfamily\thepage\hss\scshape\shortauthors\hss}
    \let\@oddfoot\@empty\let\@evenfoot\@empty}
\@twosidetrue
\makeatother
\usepackage[neveradjust]{paralist}
\makeatletter
\@definecounter{@ADL@savenum}
\newcommand\savenum{\setcounter{@ADL@savenum}%
  {\the\@nameuse{c@\@listctr}}}
\newcommand\resumenum{\setcounter{\@listctr}{\arabic{@ADL@savenum}}}
\makeatother
\makeatletter
\usepackage{xspace}
\DeclareRobustCommand\onedot{\futurelet\@let@token\@onedot}
\def\@onedot{\ifx\@let@token.\else.\null\fi\xspace}
\makeatother
\def\eg{e.g\onedot,}
\def\ie{i.e\onedot,}
\def\cf{cf\onedot}

\def\resp{resp\onedot}
\usepackage{theorem}
\newtheorem{theorem}{Theorem}[section]
\newtheorem{proposition}[theorem]{Proposition}
\newtheorem{lemma}[theorem]{Lemma}
\newtheorem{corollary}[theorem]{Corollary}
\newtheorem{prooflemma}{Lemma}[theorem]
\newtheorem{proofsublemma}{Sublemma}[theorem]

{\theorembodyfont{\normalfont\rmfamily}
\newtheorem{definition}[theorem]{Definition}
\newtheorem{remark}[theorem]{Remark}
\newtheorem{example}[theorem]{Example}
\newtheorem{examples}[theorem]{Examples}}
\makeatletter
\newcommand\pushright{\protect\@ADL@pushright}
\newcommand\@ADL@pushright[1]{{\ifvmode\null\hfill{#1}\par\else\ifmmode%
  \@ADLmaths@pushright{\hbox{#1}}\else\ifinner\@ADLhbox@pushright{#1}%
  \else\@ADLparag@pushright{#1}\fi\fi\fi}}
\newcommand\@ADLmaths@pushright[1]{{\ifinner\@ADLhbox@pushright{#1}\else%
  \tag*{$#1$}\fi}}
\newcommand\@ADLparag@pushright[1]{{\parfillskip=0pt\widowpenalty=10000%
  \displaywidowpenalty=10000\finalhyphendemerits=0\@ADLhbox@pushright#1\par}}
\newcommand\@ADLhbox@pushright{\unskip\nobreak\hfil\penalty50\hskip.2em%
  \null\hfill\hfill}
\newenvironment{proof}{\trivlist\item[\hskip\labelsep\textit{Proof:}\/]%
  \@ADLsave@set@qed\xspace\normalfont\rmfamily}
  {\qed\@ADLrestore@qed\endtrivlist}
\newenvironment{subproof}{\trivlist\item[\hskip\labelsep\textit{Proof:}\/]%
  \@ADLsave@set@subqed\normalfont\rmfamily}
  {\subqed\@ADLrestore@subqed\endtrivlist}
\newif\if@ADL@qed\@ADL@qedfalse
\newcommand\qed{\protect\@ADL@qed{$\blacksquare$}}
\newcommand\@ADL@qed[1]{\if@ADL@qed\global\@ADL@qedfalse%
  \pushright{#1}\else\ifhmode\ifinner\else\par\fi\fi\fi}
\newcommand\@ADLrestore@qed{\global\let\if@ADL@qed\@ADLsaved@ifqed}
\newcommand\@ADLsave@set@qed{\let\@ADLsaved@ifqed
  \if@ADL@qed\global\@ADL@qedtrue}
\newcommand\eqqed{\tag*{\qed}}
\newif\if@ADL@subqed\@ADL@subqedfalse
\newcommand\subqed{\protect\@ADL@subqed{$\blacktriangledown$}}
\newcommand\@ADL@subqed[1]{\if@ADL@subqed\global\@ADL@subqedfalse%
  \pushright{#1}\else\ifhmode\ifinner\else\par\fi\fi\fi}
\newcommand\@ADLrestore@subqed{\global\let\if@ADL@subqed\@ADLsaved@ifsubqed}
\newcommand\@ADLsave@set@subqed{\let\@ADLsaved@ifsubqed
  \if@ADL@subqed\global\@ADL@subqedtrue}
  
\makeatother
\makeatletter
\newif\if@ADL@oprocend\@ADL@oprocendfalse
\newcommand\oprocend{\@ADLsave@set@oprocend
  \protect\@ADL@oprocend{$\bullet$}\@ADLrestore@oprocend}
\newcommand\@ADL@oprocend[1]{\if@ADL@oprocend\global\@ADL@oprocendfalse%
  \pushright{#1}\else\ifhmode\ifinner\else\par\fi\fi\fi}
\newcommand\@ADLrestore@oprocend{\global
  \let\if@ADL@oprocend\@ADLsaved@ifoprocend}
\newcommand\@ADLsave@set@oprocend{\let\@ADLsaved@ifoprocend\if@ADL@oprocend%
  \global\@ADL@oprocendtrue}
\newcommand\eqoprocend{\tag*{\oprocend}}
\makeatother
\usepackage[reqno,centertags]{amsmath}
\usepackage{amssymb,xspace,mathrsfs}
\usepackage[mathcal]{euscript}
\let\subsetneq\subset
\let\subset\subseteq

\let\supset\supseteq
\newcommand\real{\mathbb{R}}
\newcommand\realp{\real_{>0}}
\newcommand\realnn{\real_{\ge0}}
\newcommand\realnp{\real_{\le0}}
\newcommand\rational{\mathbb{Q}}
\newcommand\integer{\mathbb{Z}}
\newcommand\integerp{\integer_{>0}}
\newcommand\integernn{\integer_{\ge0}}
\newcommand\eul{\textup{e}}
\newcommand\setdef[2]{\{#1\;|\enspace#2\}}
\newcommand\asetdef[2]{\left\{#1\immediate\vphantom{#2}\;\right|
  \left.\immediate\vphantom{#1}\enspace#2\right\}}
\newcommand\slnorm{\lvert}
\newcommand\srnorm{\rvert}
\newcommand\snorm[1]{\slnorm #1\srnorm}
\newcommand\asnorm[1]{\left\slnorm #1\right\srnorm}
\newcommand\dlnorm{\lVert}
\newcommand\drnorm{\rVert}
\newcommand\dnorm[1]{\dlnorm #1\drnorm}
\newcommand\ifam[1]{(#1)}
\newcommand\map[3]{#1\colon#2\rightarrow#3}
\newcommand\mapdef[5]{\begin{aligned}
  #1\colon&\begin{aligned}[t]#2\end{aligned}\rightarrow
  \begin{aligned}[t]#3\end{aligned}\\&\begin{aligned}[t]#4\end{aligned}
  \mapsto\begin{aligned}[t]#5\end{aligned}\end{aligned}}
\newcommand{\scirc}{\mathbin{\mathchoice
  {\xcirc\scriptstyle}
  {\xcirc\scriptstyle}
  {\xcirc\scriptscriptstyle}
  {\xcirc\scriptscriptstyle}
}}
\newcommand{\xcirc}[1]{\vcenter{\hbox{$#1\circ$}}}  
\def\mappings[#1]#2#3{\mapschar^{#1}(#2;#3)}
\makeatletter
\def\interval{\@ifnextchar({\@ADL@openleftint}{\@ADL@closedleftint}}
\def\@ADL@openleftint(#1,#2{(#1,#2%
  \@ifnextchar){\@ADL@openrightint}{\@ADL@closedrightint}}
\def\@ADL@closedleftint[#1,#2{[#1,#2%
  \@ifnextchar){\@ADL@openrightint}{\@ADL@closedrightint}}
\def\@ADL@openrightint){)}
\def\@ADL@closedrightint]{]}
\makeatother
\newenvironment{keywords}{\quote\small\textbf{Keywords.}}{\endquote}
\newenvironment{AMS}{\quote\small\textbf{AMS Subject Classifications (2010).}}
   {\endquote}
\newcommand\defn[1]{{\normalfont\bfseries\emph{\mathversion{bold}#1}}}
\newcommand\eqdef{\triangleq}
\newcommand\mathupper[1]{\textup{#1}}
\newcommand\subscr[2]{#1_{\textup{#2}}}
\newcommand\supscr[2]{#1^{\textup{#2}}}
\newcommand\interior{\operatorname{int}}
\newcommand\closure{\operatorname{cl}}

\newcommand\ts[1]{\mathcal{#1}}
\newcommand\nbhd[1]{\mathcal{#1}}

\renewcommand\d[1]{{\normalfont\textrm{d}}#1}
\newcommand\mapschar{\mathupper{C}}
\newcommand\ChRec{\mathupper{ChRec}}
\newcommand\dist{\mathupper{dist}}
\newcommand\NWnd{\mathupper{NWnd}}
\newcommand\Orb{\mathupper{Orb}}
\newcommand\id{\operatorname{id}}
\newcommand\supp{\operatorname{supp}}

\newcommand\tdomain{\mathbb{T}}
\newcommand\tdomainp{\mathbb{T}_{>0}}
\newcommand\tdomainnn{\mathbb{T}_{\ge0}}
\newcommand\tdomainn{\mathbb{T}_{<0}}
\newcommand\tdomainnp{\mathbb{T}_{\le0}}

\newcommand\sC{\mathscr{C}}
\newcommand\sNC{\mathscr{NC}}
\newcommand\sO{\mathscr{O}}
\newcommand\sP{\mathscr{P}}
\newcommand\sT{\mathscr{T}}

\newcommand\ol[1]{\overline{#1}}
\newcommand\ul[1]{\underline{#1}}
\newcommand\Phidisc[1]{\supscr{\Phi}{d,$#1$}}

\makeatletter
\def\oball{\@ifnextchar[{\@ADL@oballd}{\@ADL@oball}}
\def\@ADL@oballd[#1]#2#3{\mathsf{B}_{#1}(#2,#3)}
\newcommand\@ADL@oball[2]{\mathsf{B}(#1,#2)}
\makeatother

\title{The Fundamental Theorem of Dynamical Systems:\\all at once and all in
the same place}
\author{Andrew D.\ Lewis\thanks{Professor, Department of Mathematics and
Statistics, Queen's University, Kingston, ON K7L 3N6, Canada,
email:~\texttt{andrew.lewis@queensu.ca}}}
\newcommand\shorttitle{The Fundamental Theorem of Dynamical Systems}
\newcommand\shortauthors{A.\ D.\ Lewis}
\date{2025/08/13}
\begin{document}
\maketitle

\begin{abstract}
The so-called Fundamental Theorem of Dynamical Systems\textemdash{}which
(1)~relates attractors and repellers to the chain recurrent set and (2)~gives
the existence of a complete Lyapunov function\textemdash{}can be seen as a
means of separating out ``recurrent'' and ``transient'' dynamics.  An
overview of this theorem is given in its various guises,
continuous-time/discrete-time and flows/semiflows.  As part of this overview,
a unified approach is developed for working simultaneously with both the
continuous-time and discrete-time frameworks for topological dynamics.
Additionally, a complete Lyapunov function is provided for the first time for
continuous-time flows and semiflows.
\end{abstract}

\tableofcontents

\section{Introduction}

The behaviour of dynamical systems is too complicated to admit any sort of
very useful classification.  However, one of the basic sorts of
classifications one might seek is to separate dynamics that is
``steady-state'' from dynamics that is ``transient.''  The exact meaning of
``steady state'' is rather unclear at the outset,~\eg~should it be
equilibrium dynamics? period dynamics?  After some thought, it becomes clear
that ``recurrence'' is a good general concept for capturing ``steady-state.''
However, there is a hierarchy of notions of recurrence, with the most
elementary notion being that of a fixed point.  As one ascends (or descends,
depending on one's point of view) through the hierarchy, one naturally
wonders whether the hierarchy has finitely many steps, or whether there are
just more and more subtle notions of recurrence that become increasingly
incomprehensible.  One way to view the Fundamental Theorem of Dynamical
systems\textemdash{}the term coined by \citet{DEN:95} for the initial work of
\citet{CC:78}\textemdash{}is that it says that this hierarchy can be
terminated with the notion of recurrence known as ``chain recurrence.''  Of
course, this resolution of the difference between transient and recurrent
dynamics hardly resolves the problem of completely understanding dynamics:
chain recurrent dynamics can be extremely complicated.  Nonetheless, this
Fundamental Theorem of Dynamical Systems provides genuine insights into the
general structure of topological dynamical systems.

\subsection{Contribution}

In this paper, we shall provide a comprehensive overview of chain recurrence
and the Fundamental Theorem of Dynamical Systems.  We feel as if this
overview is warranted because (1)~the essential results are scattered across
many papers and (2)~there are a few places in the literature where some
fundamental misconceptions lead to misleading or erroneous statements.  The
reason for the theory being scattered across many papers is simple: the
notion of chain recurrence and the statement of the Fundamental Theorem of
Dynamical Systems applies to $2\times2=4$ classes of dynamical systems:
continuous-time/discrete-time and flows/semiflows.  Sometimes the differences
in the way various classes are handled can be passed off to the Latin
abbreviation ``\cf'' without incurring too much of a loss; typically the
differences between ``flow'' and ``semiflow'' fall into this category.  Other
times, a proof in one case is simply inapplicable to another case; the
construction of a complete Lyapunov function in the discrete-time case being
one such example since the construction is not immediately helpful in the
continuous-time case.  What we wish to do, therefore, is to present the
theory in a complete and unified way (as much as this is possible).  Most of
the techniques we use will be familiar to those familiar with the theory.
However, we hope that presenting this theory in a unified way will have a
\emph{per se} benefit\@.

\subsection{Historical developments}

Let us give an historical overview of the development of this theory.  The
terminology we use is given precise definitions in
Section~\ref{sec:background}\@.  A reader completely new to the ideas we
present here may benefit from first reading this section.  As we have
indicated, the theory is initially given by \citet{CC:78}\@, working with
continuous-time flows on compact metric spaces.  From this initial work,
there arises a few natural extensions: (1)~the inclusion of noncompact
spaces; (2)~the adaptation to discrete-time; (3)~the adaptation to semiflows.
Note that extensions (2)~and (3)~are essential for the theory to be applied
to the important setting of the dynamics of a continuous mapping.  These
adaptations were carried out in a series of papers during the 1980's and
90's.  When carrying out extensions of Conley's initial work, there are a few
different aspects of the work that can be extended, and normally these
extensions are not carried out together, but piecemeal.  These aspects are
the following.
\begin{compactenum}
\item \emph{Chain recurrence:} The Fundamental Theorem of Dynamical Systems
is connected with the particular recurrence notion of chain recurrence.
There are different  definitions of chain recurrence, depending on the class
of dynamical system with which one is working.  As well, chain recurrence,
and particularly the connected notion of chain transitivity, has its own
properties that can be important to understand, independently of the
Fundamental Theorem of Dynamical Systems.

\item \emph{The Conley decomposition:} Part of the Fundamental Theorem of
Dynamical Systems is a decomposition of the state space into (a)~a subset on
which the dynamics is chain recurrent and (b)~the complement on which the
dynamics is gradient-like.

\item \emph{Complete Lyapunov functions:} Another part of the Fundamental
Theorem of Dynamical Systems concerns the existence of a complete Lyapunov
function that decreases along trajectories, and is constant on the chain
recurrent set.
\end{compactenum}
We shall consider how all three of these aspects have developed, either
together or separately.

As we indicated above, chain recurrence is part of a theory of recurrence for
topological dynamical systems.  Other notions of recurrence include limit
sets, Poincar\'e recurrent sets, and nonwandering sets.  A notion of ``weak
nonwandering point'' for ordinary differential equations is introduced by
\citet{ANS/VAD:73}\@, and this notion can apparently be shown to be
equivalent to chain recurrence as introduced by~\cite{CC:78} (we have seen
this ``on the internet'' but know of no precise reference).  The (fairly
straightforward) adaptation to discrete-time semiflows seems to have first
been given by \citet{LB/JEF:85}\@, where connections to other forms of
recurrence is proved in some special cases.  The discrete-time setting is
also considered in~\cite{RE:89} and reflections are made on numerical aspects
of dynamics.  An expository presentation is given by~\cite{JEF:88} that
considers chain recurrence in relation to other recurrence notions.

As introduced by \citeauthor{CC:78}\@, chain recurrence is defined for
compact metric spaces.  In its original form, for noncompact spaces the
notion of chain recurrence is a metric notion, not a topological notion.  We
shall see an instance of this in Example~\ref{eg:metric-dependent}\@.  An
approach to making the concept less connected to the metric was introduced
in~\cite{MH:91}\@.  In this work, the setting is discrete-time semiflows on
locally compact metric spaces.  As well as giving a definition of chain
recurrence in this setting, \citeauthor{MH:91} also shows that the Conley
decomposition holds.  The local compactness is relaxed to general metric
spaces in~\cite{MH:92}\@, still in the discrete-time semiflow framework.
Also in this work, the existence of a complete Lyapunov function is proved
for locally compact, second countable state spaces.  The extension to general
metric spaces is developed further by~\citet{MH:95} to include the
continuous-time framework as well as the discrete-time.  Additionally,
\citeauthor{MH:95} gives two useful alternative characterisations of
continuous-time chain recurrence; for instance, chain recurrence for a
continuous-time semiflow is related to chain recurrence for its time-one map.
The work~\cite{MH:95} contained a few errors that were corrected
by~\cite{SKC/CKC/JSP:02}\@, while introducing a few other errors themselves.
We hope we have corrected these errors in our presentation.  The existence of
a complete Lyapunov function is proved by~\cite{MH:98} for discrete-time
semiflows on separable metric spaces.  This is extended in \cite{MP:11} to
the continuous-time case, although the proof has errors that we correct in
our presentation.  For flows defined by smooth ordinary differential
equations, \citet{SH/S:21} show the existence of smooth complete Lyapunov
functions.

Extensions of chain recurrence and the Fundamental Theorem of Dynamical
Systems beyond the metric space setting are also possible, although we do not
consider these in our overview.  A thorough development of these notions for
uniform spaces is given in the manuscript~\cite{EA/JW:17}\@.  Chain
recurrence and related notions are considered for discrete-time semiflows on
arbitrary topological spaces in~\cite{LSB/WAC:92}\@.  This treatment is
extended to continuous-time semiflows in~\cite{PO:05}\@.  An ``open cover''
framework for chain recurrence is presented by \citet{MP/LABSM:06}\@.

Recent developments in the concepts of chain recurrence and the Fundamental
Theorem of Dynamical Systems include the applications to linear dynamical
systems, including those in
infinite-dimensions~\cite{MBA/GEM/RV:22,NCB/AP:24,}\@, control
theory~\cite{FC/AJS/ECV:24}\@, and hybrid dynamical
systems~\cite{MDK/PG/DEK:21}\@.

\section{Dynamical systems and related concepts}\label{sec:background}

In this section, we give our definitions for the classes of dynamical systems
we use throughout the paper, as well as a few of the concepts from
topological dynamics to which we refer.  As the reader will see, we unite the
continuous- and discrete-time in our presentation as much as this is
possible.  This has a variety of conceptual benefits.  We do not, however, go
as far as \citet{EA:93} who uses general relations in place of continuous-
and discrete-time.

\subsection{Flows and semiflows}

We use the symbol $\tdomain$ to stand for either $\real$ or $\integer$\@.  If
we do not explicitly specify, we intend that $\tdomain$ can be either of the
two possibilities.  For $t_0\in\tdomain$ and for an interval
$I\subset\real$\@, We denote
\begin{align*}
&\tdomain_{>t_0}=\setdef{t\in\tdomain}{t>t_0},\\
&\tdomain_{\ge{}t_0}=\setdef{t\in\tdomain}{t\ge t_0},\\
&\tdomain_{<t_0}=\setdef{t\in\tdomain}{t<t_0},\\
&\tdomain_{\le{}t_0}=\setdef{t\in\tdomain}{t\le t_0},\\
&\tdomain_I=\setdef{t\in\tdomain}{t\in I}.
\end{align*}

Let us give the definition of the dynamical systems we work with.
\begin{definition}
Let $(\ts{X},\sO)$ be a topological space.
\begin{compactenum}[(i)]
\item A \defn{topological semiflow} on $\ts{X}$ is a continuous mapping
$\map{\Phi}{\tdomainnn\times\ts{X}}{\ts{X}}$ satisfying
\begin{compactenum}[(a)]
\item $\Phi(0,x)=x$\@, $x\in\ts{X}$\@, and
\item $\Phi(t_1,\Phi(t_2,x))=\Phi(t_1+t_2,x)$\@, $t_1,t_2\in\tdomainnn$\@,
$x\in\ts{X}$\@.
\end{compactenum}
\item A \defn{topological flow} on $\ts{X}$ is a continuous mapping
$\map{\Phi}{\tdomain\times\ts{X}}{\ts{X}}$ satisfying
\begin{compactenum}[(a)]
\item $\Phi(0,x)=x$\@, $x\in\ts{X}$\@, and
\item $\Phi(t_1,\Phi(t_2,x))=\Phi(t_1+t_2,x)$\@, $t_1,t_2\in\tdomain$\@,
$x\in\ts{X}$\@.
\end{compactenum}
\end{compactenum}
If $\tdomain=\real$\@, the topological semiflow or flow is
\defn{continuous-time}\@, and otherwise it is
\defn{discrete-time}\@.\oprocend
\end{definition}

\begin{remark}\label{rem:dt-flow-map}
Discrete-time flows and semiflows correspond to the dynamics of continuous
mappings and homeomorphisms, respectively.
\begin{compactenum}
\item If $\Phi$ is a discrete-time topological semiflow, then we can define
$\phi(x)=\Phi(1,x)$ and easily very that
\begin{equation*}
\Phi(t,x)=\phi^t(x)\eqdef\underbrace{\phi\scirc\cdots\scirc\phi}_{t\ \textrm{times}}(x),\qquad(t,x)\in\integernn\times\ts{X},
\end{equation*}
where it is understood that $\phi^0=\id_{\ts{X}}$\@.

\item Similarly, if $\Phi$ is a discrete-time topological flow and if
$\phi(x)=\Phi(1,x)$\@, then $\phi$ is an homeomorphism with
$\phi^{-1}(x)=\Phi(-1,x)$\@.  Here we have $\Phi(t,x)$ as above for
$t\in\integernn$\@, and
\begin{equation*}
\Phi(-t,x)=\phi^{-t}(x)\eqdef
\underbrace{\phi^{-1}\scirc\cdots\scirc\phi^{-1}}_{t\ \textrm{times}}(x),\qquad(t,x)\in\integernn\times\ts{X},
\end{equation*}
\end{compactenum}
Note that our use of ``flow'' or ``semiflow'' in the discrete-time case is
nonstandard.  However, there is a unifying benefit to using the same
terminology for both the continuous- and discrete-time settings.\oprocend
\end{remark}

There are somewhat surprising relationships between chain recurrence for a
continuous-time flow or semiflow, and chain recurrence for the discrete-time
flow or semiflow defined by its time-one map.  For this reason, we make the
following definition.
\begin{definition}
Let $(\ts{X},\sO)$ be a topological space and let $\Phi$ be a topological
continuous-time topological flow (\resp~semiflow) on $\ts{S}$\@.  The
\defn{$T$-discretisation} of $\Phi$ is the discrete-time topological flow
(\resp~semiflow) $\Phidisc{T}$ on $\ts{X}$ defined by requiring that
\begin{equation*}
\Phidisc{T}(1,x)=\Phi(T,x),\qquad x\in\ts{X},
\end{equation*}
\cf~Remark~\ref{rem:dt-flow-map}\@.\oprocend
\end{definition}

It is evident that, if $\Phi$ is a topological flow, then its restriction to
$\tdomainnn\times\ts{X}$ is topological semiflow.  If $\Phi$ is a topological
flow (\resp~semiflow) and if $t\in\tdomain$ (\resp~$t\in\tdomainnn$), then we
have the homeomorphism (\resp~continuous mapping)
\begin{equation*}
\mapdef{\Phi_t}{\ts{X}}{\ts{X}}{x}{\Phi(t,x).}
\end{equation*}
In like manner, if $x\in\ts{X}$\@, then
\begin{equation*}
\mapdef{\Phi^{+,x}}{\tdomainnn}{\ts{X}}{t}{\Phi(t,x)}
\end{equation*}
is the \defn{forward trajectory} of $x$ and, if $\Phi$ is a flow, then
\begin{equation*}
\mapdef{\Phi^x}{\tdomain}{\ts{X}}{t}{\Phi(t,x)}
\end{equation*}
is the \defn{trajectory} of $x$\@.

Let us also define the notion of orbits.
\begin{definition}
Let $(\ts{X},\sO)$ be a topological space and let $\Phi$ be a topological
flow or semiflow on $\ts{X}$\@.  Let $A\subset\ts{X}$\@.
\begin{compactenum}[(i)]
\item The \defn{forward orbit} of $A$ is
\begin{equation*}
\Orb^+(A)=\setdef{\Phi(t,x)}{t\in\tdomainnn,\ x\in A}.
\end{equation*}
\item The \defn{backward orbit} of $A$ is
\begin{equation*}
\Orb^-(A)=\setdef{x\in\ts{X}}
{\Phi(t,x)\in A\ \textrm{for some}\ t\in\tdomainnn}.
\end{equation*}\savenum
\end{compactenum}
If $\Phi$ is a topological flow, then
\begin{compactenum}[(i)]\resumenum
\item the \defn{orbit} of $A$ is
\begin{equation*}\eqoprocend
\Orb(A)=\setdef{\Phi(t,x)}{t\in\tdomain,\ x\in A}.
\end{equation*}
\end{compactenum}
\end{definition}

If $A=\{x\}$ for some $x\in\ts{X}$\@, we abbreviate
\begin{equation*}
\Orb^+(x)=\Orb^+(\{x\}),\enspace\Orb^-(x)=\Orb^-(\{x\}),\enspace
\Orb(x)=\Orb(\{x\}).
\end{equation*}
Note that, for a topological flow $\Phi$\@, the backward orbit of $A$ is
\begin{equation*}
\Orb^-(A)=\setdef{\Phi(t,x)}{t\in\tdomainnp,\ x\in A}.
\end{equation*}
The difference in how one should view the backward orbit for a flow and a
semiflow will come up in our treatment of the Conley decomposition.

We shall make use of the following invariance notions.
\begin{definition}
Let $(\ts{X},\sO)$ be a topological space and let $\Phi$ be a topological
flow or semiflow.  Let $A\subset\ts{X}$\@.
\begin{compactenum}[(i)]
\item The set $A$ is \defn{forward-invariant} for $\Phi$ if
$\Orb^+(A)\subset A$\@.
\item If $\Phi$ is a flow, the set $A$ is \defn{invariant} for $\Phi$ if
$\Orb(A)\subset A$\@.\oprocend
\end{compactenum}
\end{definition}

The following elementary lemma will account for some differences in how
certain proofs work for flows versus semiflows.
\begin{lemma}\label{lem:comp-inv}
Let\/ $(\ts{X},\sO)$ be a topological space and let\/ $\Phi$ be a topological
flow on\/ $\ts{X}$\@.  If\/ $A\subset X$ is invariant for\/ $\Phi$\@, then\/
$\ts{X}\setminus A$ is invariant for\/ $\Phi$\@.
\begin{proof}
Let $x\in\ts{X}\setminus A$ and let $t\in\tdomain$\@.  If $\Phi(t,x)\in A$\@,
then $x=\Phi(-t,x)\in A$ by invariance of $A$\@.  Therefore, we must have
$\Phi(t,x)\in\ts{X}\setminus A$\@.
\end{proof}
\end{lemma}

Note that the lemma is false for forward-invariance and semiflows; think of a
semiflow for which every trajectory ends up at a fixed point in finite time.

\subsection{Attracting sets, repelling sets, and basins}

The Conley decomposition relates the chain recurrent set (which we define in
Section~\ref{sec:chain-stuff}) to attracting sets and their basins (for
semiflows) or attracting sets and repelling sets (for flows), In this section
we introduce the necessary terminology.
\begin{definition}
Let $(\ts{X},\sO)$ be a topological space and let $\Phi$ be a topological
flow or semiflow on $\ts{X}$\@.  A nonempty subset $\nbhd{T}\subset\ts{X}$ is
a \defn{trapping region} for $\Phi$ if there exists\/ $T\in\tdomainp$ such
that
\begin{equation*}\eqoprocend
\closure(\Phi(\tdomain_{\ge{}T}\times\nbhd{T}))\subset\interior(\nbhd{T}).
\end{equation*}
\end{definition}

\begin{remark}\label{rem:special-trap}
Different definitions and terminology can be found for what we call a
trapping region.  For the continuous-time case, sometimes ``preattractor'' is
used.  In the discrete-time case, the terminology ``attractor block'' is
sometimes used.  Again, we prefer unified terminology.  Let us also examine
more serious differences in the definitions that can be found in the
literature.
\begin{compactenum}
\item \label{enum:special-trap1} Note that, in the discrete-time case, our
definition of a trapping region is \emph{not} the usual definition.  The
usual definition in the discrete-time case is that
$\closure(\Phi_1(\nbhd{T}))\subset\interior(\nbhd{T})$\@.  We shall call a
subset $\nbhd{T}$ a \defn{strong trapping region} for the discrete-time flow
or semiflow $\Phi$\@.  Note that a strong trapping region is a trapping
region.  Note, also, that a strong trapping region is forward-invariant,
which leads to some simplifications for strong trapping regions compared to
trapping regions.  It is easy to build examples for which the two definitions
of trapping region are not the same.

\item \label{enum:special-trap2} Somewhat in keeping with the notion of a
strong trapping region in the discrete-time case, one can consider a notion
of trapping region in the continuous-time case where the requirement is that
there exists $T\in\realp$ such that
$\closure(\Phi_T(\nbhd{T}))\subset\interior(\nbhd{T})$\@.  In the case when
$\ts{X}$ is compact, this condition agrees with our notion of a trapping
region.  This is proved for flows by \citet[page~33, C]{CC:78}\@; we prove
this here for semiflows as well.  Let $T\in\realp$ be such that
$\closure(\Phi_T(\nbhd{T}))\subset\interior(\nbhd{T})$\@.  Because metric
spaces are normal~\cite[Example~15.3(c)]{SW:70}\@, let $\nbhd{N}$ be an open
set for which
\begin{equation*}
\Phi_T(\closure(\nbhd{T}))\subset\closure(\Phi_T(\nbhd{T}))\subset
\nbhd{N}\subset\closure(\nbhd{N})\subset\interior(\nbhd{T}).
\end{equation*}
Note that $\Phi_T(\nbhd{T})$ is compact, being a closed subset of a
compact space.  By this compactness and by the continuity of the dynamics,
there exists $\delta\in\realp$ such that
\begin{equation*}
\Phi(\interval({T-\delta},{T+\delta})\times\closure(\nbhd{T}))
\subset\nbhd{N}.
\end{equation*}
Let $T'=\frac{T}{2\delta}$\@.  If $t\ge T'$\@, then $t$ is a finite sum
$t=t_1+\dots+t_k$ where $t_j\in\interval({T-\delta},{T+\delta})$\@,
$j\in\{1,\dots,k\}$\@.  Note that
\begin{equation*}
\Phi_{t_k}(\closure(\nbhd{T}))\subset\nbhd{N}\subset\interior(\nbhd{T})
\subset\closure(\nbhd{T}).
\end{equation*}
Thus
\begin{equation*}
\Phi_{t_{k-1}}\scirc\Phi_{t_k}(\closure(\nbhd{T}))
\subset\Phi_{t_{k-1}}(\closure(\nbhd{T}))\subset\nbhd{N}
\subset\interior(\nbhd{T})\subset\closure(\nbhd{T}).
\end{equation*}
Inductively,
\begin{equation*}
\Phi_t(\closure(\nbhd{T}))=
\Phi_{t_1}\scirc\Phi_{t_2}\scirc\Phi_{t_k}(\closure(\nbhd{T}))
\subset\Phi_{t_1}(\closure(\nbhd{T}))\subset\nbhd{N}.
\end{equation*}
Therefore,
\begin{align*}
\closure(\Phi(\interval[T',\infty)\times\nbhd{T}))
\subset&\;\closure\left(\bigcup_{t\in\interval[T',\infty)}
\Phi_t(\nbhd{T})\right)\\
\subset&\;\closure\left(\bigcup_{t\in\interval[T',\infty)}
\Phi_t(\closure(\nbhd{T}))\right)\\
\subset&\;\closure(\nbhd{N})\subset
\interior(\nbhd{T}),
\end{align*}
and so $\nbhd{T}$ is a trapping region.

\citet[Example~1]{MH:95} gives a simple continuous-time flow which shows that
compactness is required for this assertion.\oprocend
\end{compactenum}
\end{remark}

Trapping regions give rise to attracting sets and (in the case of flows)
repelling sets.  Note that we refrain from calling these ``attractors'' and
``repellers'' since one normally wants to reserve this terminology for
attracting sets and repelling sets with some minimality or transitivity
property.
\begin{definition}
Let\/ $(\ts{X},\sO)$ be a topological space and let\/ $\Phi$ be a topological
flow or semiflow.
\begin{compactenum}[(i)]
\item An \defn{attracting set} for $\Phi$ is a subset $A$ such that there
exists a trapping region $\nbhd{T}$ for which
\begin{equation*}
A=\bigcap_{t\in\tdomainnn}\closure(\Phi(\tdomain_{\ge{}t}\times\nbhd{T})).
\end{equation*}\savenum
\end{compactenum}
Suppose now that $\Phi$ is a flow.
\begin{compactenum}[(i)]\resumenum
\item A \defn{repelling set} for $\Phi$ is a subset $R$ such that there
exists a trapping region $\nbhd{T}$ for which
\begin{equation*}
R=\bigcap_{t\in\tdomainnp}\closure(\Phi(\tdomain_{\le{}t}\times
(\ts{X}\setminus\nbhd{T}))).
\end{equation*}
\item An \defn{attracting-repelling pair} for $\Phi$ is a pair $(A,R)$ of
subsets such that there exists a trapping region $\nbhd{T}$ for which
\begin{equation*}\eqoprocend
A=\bigcap_{t\in\tdomainnn}\closure(\Phi(\tdomain_{\ge{}t}\times\nbhd{T})),\quad
R=\bigcap_{t\in\tdomainnp}\closure(\Phi(\tdomain_{\le{}t}\times(\ts{X}
\setminus\nbhd{T}))).
\end{equation*}
\end{compactenum}
\end{definition}

If we wish to indicate that an attracting set or a repelling set comes from a
particular trapping region $\nbhd{T}$\@, we may write $A_{\nbhd{T}}$ or
$R_{\nbhd{T}}$\@, respectively.  It is entirely possible that attracting and
repelling sets may be empty.
\begin{examples}
\begin{compactenum}
\item Let $\ts{X}=\interval[0,1]$ and let $\Phi$ be the continuous-time
flow on $\ts{X}$ obtained by restricting to $\ts{X}$ the flow associated with
the ordinary differential equation
\begin{equation*}
\dot{x}(t)=x(t)(1-x(t))
\end{equation*}
for on $\real$\@.  Then the trapping region
$\nbhd{T}=\interval({\frac{1}{2}},1]$ has the attracting set
$A_{\nbhd{T}}=\{1\}$ and the repelling set $R_{\nbhd{T}}=\{0\}$\@.  

\item On $\ts{X}=\interval(-\infty,0]$\@, consider the continuous-time flow
$\Phi(t,x)=x\eul^{-t}$\@.  Then the trapping region $\nbhd{T}=\interval(-1,0]$
has the attracting set $A_{\nbhd{T}}=\{0\}$ and the repelling set
$R_{\nbhd{T}}=\emptyset$\@.

\item On $\ts{X}=\interval[0,\infty)$\@, consider the continuous-time flow
$\Phi(t,x)=x\eul^t$\@.  Then the trapping region
$\nbhd{T}=\interval(1,\infty)$ has the attracting set
$A_{\nbhd{T}}=\emptyset$ and the repelling set $R_{\nbhd{T}}=\{0\}$\@.

\item On $\ts{X}=\real$\@, consider the continuous-time flow
$\Phi(t,x)=x+t$\@.  Then the trapping region $\nbhd{T}=\interval(0,\infty)$
has the attracting set $A_{\nbhd{T}}=\emptyset$ and the repelling set
$R_{\nbhd{T}}=\emptyset$\@.\oprocend
\end{compactenum}
\end{examples}

Let us show that trapping regions can be taken to be either open or closed.
\begin{proposition}\label{prop:clop-trapping}
Let\/ $(\ts{X},\sO)$ be a topological space and let\/ $\Phi$ be an autonomous
topological flow or semiflow.  If\/ $\nbhd{T}$ is a trapping region, then
\begin{equation*}
A_{\nbhd{T}}=A_{\closure(\nbhd{T})}=A_{\interior(\nbhd{T})}.
\end{equation*}
If\/ $\Phi$ is a flow, then
\begin{equation*}
R_{\nbhd{T}}=R_{\closure(\nbhd{T})}=R_{\interior(\nbhd{T})}.
\end{equation*}
\begin{proof}
First we show that $\interior(\nbhd{T})$ and $\closure(\nbhd{T})$ are
trapping regions.  Let $T\in\tdomainp$ be such that
\begin{equation*}
\closure(\Phi(\tdomain_{\ge{}T}\times\nbhd{T}))\subset\interior(\nbhd{T}).
\end{equation*}
This, in particular, implies that 
\begin{equation*}
\closure(\Phi(\tdomain_{\ge{}T}\times\interior(\nbhd{T})))\subset
\interior(\nbhd{T})
\end{equation*}
and
\begin{equation*}
\closure(\Phi(\tdomain_{\ge{}T}\times\closure(\nbhd{T})))\subset
\closure(\closure(\Phi(\tdomain_{\ge{}T}\times\nbhd{T})))=
\closure(\Phi(\tdomain_{\ge{}T}\times\nbhd{T}))\subset
\interior(\nbhd{T})\subset\interior(\closure(\nbhd{T})),
\end{equation*}
which shows that $\interior(\nbhd{T})$ and $\closure(\nbhd{T})$ are indeed
trapping regions.

Now note that
\begin{equation*}
\Phi_T(\closure(\nbhd{T}))\subset
\closure(\Phi(\tdomain_{\ge{}T}\times\nbhd{T}))\subset\interior(\nbhd{T}),
\end{equation*}
and so
\begin{equation*}
\Phi_{T+t}(\closure(\nbhd{T}))=\Phi_t(\Phi_T(\closure(\nbhd{T})))
\subset\Phi_t(\interior(\nbhd{T})),\qquad t\in\tdomainnn.
\end{equation*}
Now calculate
\begin{equation}\label{eq:Aclint1}
\begin{aligned}
A_{\nbhd{T}}=&\;
\bigcap_{t\in\tdomainnn}\closure(\Phi(\tdomain_{\ge{}t}\times\nbhd{T}))
\subset\bigcap_{t\in\tdomainnn}\closure(\Phi(\tdomain_{\ge{}t}\times
\closure(\nbhd{T})))\\
\subset&\;\bigcap_{t\in\tdomainnn}\closure(\Phi(\tdomain_{\ge{}T+t}\times
\closure(\nbhd{T})))
\subset\bigcap_{t\in\tdomainnn}
\closure(\Phi(\tdomain_{\ge{}t}\times\interior(\nbhd{T})))
\subset A_{\nbhd{T}},
\end{aligned}
\end{equation}
and so, in particular,
\begin{equation*}
A_{\nbhd{T}}=\bigcap_{t\in\tdomainnn}
\closure(\Phi(\tdomain_{\ge{}t}\times\closure(\nbhd{T})))
=\bigcap_{t\in\tdomainnn}
\closure(\Phi(\tdomain_{\ge{}t}\times\interior(\nbhd{T}))),
\end{equation*}
which is the first part of the result.

For the second, note that $\closure(\Phi_T(A))=\Phi_T(\closure(A))$ and
$\interior(\Phi_T(A))=\Phi_T(\interior(A))$ for any $A\subset\ts{X}$ since
$\Phi_T$ is an homeomorphism when $\Phi$ is a flow.  Also note that
$\Phi_T(\ts{X}\setminus A)=\ts{X}\setminus\Phi_T(A)$ for any $A\subset\ts{X}$
since $\Phi_T$ is a bijection.  Therefore,
\begin{equation*}
\begin{aligned}[t]
\closure(\Phi_T(\nbhd{T}))\subset\interior(\nbhd{T})
\end{aligned}\implies
\begin{aligned}[t]
\ts{X}\setminus\interior(\nbhd{T})\subset&\;
\ts{X}\setminus\closure(\Phi_T(\nbhd{T}))
=\interior(\ts{X}\setminus\Phi_T(\nbhd{T}))\\
=&\;\interior(\Phi_T(\ts{X}\setminus\nbhd{T}))=
\Phi_T(\interior(\ts{X}\setminus\nbhd{T})).
\end{aligned}
\end{equation*}
Therefore,
\begin{equation*}
\Phi_{-T}(\closure(\ts{X}\setminus\nbhd{T}))=
\Phi_{-T}(\ts{X}\setminus\interior(\nbhd{T}))\subset
\interior(\ts{X}\setminus\nbhd{T}).
\end{equation*}
Thus, by a modification of the arguments in the previous part of the proof,
we have
\begin{align*}
R_{\nbhd{T}}=&\;
\bigcap_{t\in\tdomainnp}\closure(\Phi(\tdomain_{\le{}t}\times
(\ts{X}\setminus\nbhd{T})))
\subset\bigcap_{t\in\tdomainnp}\closure(\Phi(\tdomain_{\le{}t}\times(\ts{X}
\setminus\interior(\nbhd{T}))))\\
\subset&\;\bigcap_{t\in\tdomainnp}
\closure(\Phi(\tdomain_{\le{}t-T}\times(\closure(\ts{X}\setminus\nbhd{T}))))
\subset\bigcap_{t\in\tdomainnp}
\closure(\Phi(\tdomain_{\le{}t}\times(\interior(\ts{X}\setminus\nbhd{T}))))
\subset R_{\nbhd{T}},
\end{align*}
giving the result for repelling sets.
\end{proof}
\end{proposition}

The definitions of attracting sets and repelling sets suggest a ``duality''
between these notions, depending on whether time goes forwards or backwards.
The following result makes this suggestion precise.  We denote by
$\Phi^\sigma$ the time-reversed flow of a flow $\Phi$ defined by $\Phi^\sigma(t,x)=\Phi(-t,x)$\@.
\begin{proposition}\label{prop:att<->rep}
Let\/ $(\ts{X},\d)$ be a metric space and let\/ $\Phi$ be a topological flow
on\/ $\ts{X}$\@.  If\/ $\nbhd{T}$ is a trapping region with
attracting-repelling pair\/ $(A,R)$\@, then\/ $\ts{X}\setminus\nbhd{T}$ is a
trapping region for\/ $\Phi^\sigma$ and\/ $(R,A)$ is the attracting-repelling
pair associated to\/ $\ts{X}\setminus\nbhd{T}$\@.
\begin{proof}
Let $T\in\tdomainp$ be such that
$\closure(\Phi(\tdomain_{\ge{}T}\times\nbhd{T}))\subset\interior(\nbhd{T})$\@.
Denote
\begin{equation*}
\nbhd{N}=\ts{X}\setminus
\closure\left(\bigcup_{t\in\tdomain_{\interval[0,T]}}\Phi_{-t}(\nbhd{T})\right).
\end{equation*}
The following facts about $\nbhd{N}$ are useful.
\begin{prooflemma}
The following statements hold:
\begin{compactenum}[(i)]
\item \label{pl:att<->rep1}
$\closure(\nbhd{N})\subset\interior(\ts{X}\setminus\nbhd{T})$\@;
\item \label{pl:att<->rep2}
$\Phi_{-t}(\ts{X}\setminus\interior(\nbhd{T}))\subset\nbhd{N}$\@,\/
$t\in\tdomain_{\ge{}2T}$\@.
\end{compactenum}
\begin{subproof}
\eqref{pl:att<->rep1} Since metric spaces are
normal~\cite[Example~15.3(c)]{SW:70}\@, let $\nbhd{N}'$ be an open set such
that
\begin{equation*}
\Phi_T(\closure(\nbhd{T}))\subset\nbhd{N}'\subset\closure(\nbhd{N}')
\subset\nbhd{T}.
\end{equation*}
Thus
\begin{equation*}
\closure(\nbhd{T})\subset\Phi_{-T}(\nbhd{N}')\subset
\closure(\Phi_{-T}(\nbhd{N}'))\subset\Phi_{-T}(\nbhd{T}),
\end{equation*}
and so
\begin{equation*}
\ts{X}\setminus\closure(\nbhd{T})\supset
\ts{X}\setminus\Phi_{-T}(\nbhd{N}')\supset
\ts{X}\setminus\Phi_{-T}(\nbhd{T})\supset
\ts{X}\setminus
\closure\left(\bigcup_{t\in\tdomain_{\interval[0,T]}}\Phi_{-t}(\nbhd{T})\right)
=\nbhd{N}.
\end{equation*}
Since the set $\ts{X}\setminus\Phi_{-T}(\nbhd{N}')$ is closed, we obtain
\begin{equation*}
\closure(\nbhd{N})\subset\ts{X}\setminus\closure(\nbhd{T})=
\interior(\ts{X}\setminus\nbhd{T}),
\end{equation*}
as desired.

\eqref{pl:att<->rep2} Let $t\in\tdomain_{\ge{}2T}$\@.  Then
\begin{align*}
\Phi_t(\ts{X}\setminus\nbhd{N})=&\;\closure\left(\Phi_t\left(
\bigcup_{t'\in\tdomain_{\interval[0,T]}}\Phi_{-t'}(\nbhd{T})\right)\right)
=\closure\left(\bigcup_{t'\in\tdomain_{\interval[0,T]}}
\Phi_{t-t'}(\nbhd{T})\right)\\
\subset&\;\closure\left(\bigcup_{s\in\tdomain_{\ge{}T}}
\Phi_s(\nbhd{T})\right)\subset\interior(\nbhd{T});
\end{align*}
here we have used the fact that $\Phi_t$ commutes with closure since it is an
homeomorphism and commutes with union since it is a bijection.  Therefore,
\begin{equation*}
\ts{X}\setminus\nbhd{N}\subset\Phi_{-t}(\interior(\nbhd{T}))=
\ts{X}\setminus(\ts{X}\setminus\Phi_{-t}(\interior(\nbhd{T}))),
\end{equation*}
whereupon $\Phi_{-t}(\ts{X}\setminus\interior(\nbhd{T}))\subset\nbhd{N}$\@,
as claimed.
\end{subproof}
\end{prooflemma}

The proof of the proposition is now straightforward.  First of all
\begin{align*}
\closure(\Phi^\sigma(\tdomain_{\ge{}2T}\times(\ts{X}\setminus\nbhd{T})))
\subset&\;\closure\left(\bigcup_{t\in\tdomain_{\le{}-2T}}
\Phi_t(\ts{X}\setminus\nbhd{T})\right)\\
\subset&\;\closure\left(\bigcup_{t\in\tdomain_{\le{}-2T}}
\Phi_t(\closure(\ts{X}\setminus\nbhd{T}))\right)\\
=&\;\closure(\nbhd{N})\subset\interior(\ts{X}\setminus\nbhd{T}),
\end{align*}
which shows that $\ts{X}\setminus\nbhd{T}$ is a trapping region for
$\Phi^\sigma$\@.  Also,
\begin{equation*}
A_{\nbhd{T}}=\bigcap_{t\in\realnn}\Phi_t(\nbhd{T})=
\bigcap_{t\in\realnp}
\Phi_{-t}(\ts{X}\setminus(\ts{X}\setminus\nbhd{T}))=
R_{\ts{X}\setminus\nbhd{T}}
\end{equation*}
and
\begin{equation*}
R_{\nbhd{T}}=\bigcap_{t\in\realnp}\Phi_t(\ts{X}\setminus\nbhd{T})=
\bigcap_{t\in\realnn}\Phi_{-t}(\ts{X}\setminus\nbhd{T})=
A_{\ts{X}\setminus\nbhd{T}},
\end{equation*}
where $A_{\ts{X}\setminus\nbhd{T}}$ and $R_{\ts{X}\setminus\nbhd{T}}$ denote
the attracting and repelling sets, respectively, for $\Phi^\sigma$
associated to the trapping region $\ts{X}\setminus\nbhd{T}$\@.
\end{proof}
\end{proposition}

The following properties of attracting and repelling sets will be used.
\begin{proposition}\label{prop:attract-repell-props}
Let\/ $(\ts{X},\sO)$ be a topological space and let\/ $\Phi$ be a topological
flow or semiflow.  Let\/ $A\subset\ts{X}$ be an attracting set and, when\/
$\Phi$ is a flow, let\/ $R$ be a repelling set.  Then the following
statements hold:
\begin{compactenum}[(i)]
\item \label{pl:attract-repell-props1} $A$ is closed;
\item \label{pl:attract-repell-props2} if\/ $\Phi$ is a flow, then\/ $A$
and\/ $R$ are closed;
\item \label{pl:attract-repell-props3} $A$ and\/ $R$ (\resp~$A$) are
invariant (\resp~forward-invariant).
\end{compactenum}
\begin{proof}
\eqref{pl:attract-repell-props1} and~\eqref{pl:attract-repell-props2} follow
since $A$ is (or $A$ and $R$ are) an intersection of closed sets.

\eqref{pl:attract-repell-props3} We let $\nbhd{T}$ be a trapping region with
attracting set $A$\@.  Let $x\in A$ and let $s\in\tdomainnn$\@.  Then, for
$t\in\tdomainnn$\@, we have
$x\in\closure(\Phi(\tdomain_{\ge{}t}\times\nbhd{T}))$\@.  Since
\begin{equation*}
\Phi_s(\closure(\Phi(\tdomain_{\ge{}t}\times\nbhd{T})))\subset
\closure(\Phi_s(\Phi(\tdomain_{\ge{}t}\times\nbhd{T})))=
\closure(\Phi(\tdomain_{\ge{}s+t}\times\nbhd{T})),
\end{equation*}
we have $\Phi_s(x)\in\closure(\Phi(\tdomain_{\ge{}s+t}\times\nbhd{T}))$\@.
As this holds for every $t\in\tdomainnn$\@, we have
\begin{equation*}
\Phi_s(x)\in\bigcap_{t\in\tdomainnn}
\closure(\Phi(\tdomain_{\ge{}s+t}\times\nbhd{T}))=A.
\end{equation*}
Thus $A$ is always forward-invariant.

Now suppose that $\Phi$ is a flow and let $s\in\tdomainn$\@.  In this case,
for $t\in\tdomainnn$\@, we have
$x\in\closure(\Phi(\tdomain_{\ge{}t-s}\times\nbhd{T}))$\@.  Since
\begin{equation*}
\Phi_s(\closure(\Phi(\tdomain_{\ge{}t-s}\times\nbhd{T})))\subset
\closure(\Phi_s(\Phi(\tdomain_{\ge{}t-s}\times\nbhd{T}))=
\closure(\Phi(\tdomain_{\ge{}t}\times\nbhd{T})),
\end{equation*}
we have $\Phi_s(x)\in\closure(\Phi(\tdomain_{\ge{}t}\times\nbhd{T}))$\@.  As
this holds for every $t\in\tdomainnn$\@, we have
\begin{equation*}
\Phi_s(x)\in\bigcap_{t\in\tdomainnn}
\closure(\Phi(\tdomain_{\ge{}t}\times\nbhd{T}))=A.
\end{equation*}

One similarly shows that $R$ is invariant in case $\Phi$ is a flow.
\end{proof}
\end{proposition}

Of course, the notion of a repelling set is not applicable to semiflows, by
our definition.  However, sometimes it is useful to think of
$\ts{X}\setminus\Orb^-(\nbhd{T})$ as being the ``repelling set'' for a
trapping region $\nbhd{T}$\@.  The following simple lemma indicates that this
is a reasonable way to think of repelling sets for semiflows.
\begin{lemma}\label{lem:repelling-negorb}
Let\/ $(\ts{X},\d)$ be a metric space and let\/ $\Phi$ be a topological flow
on\/ $\ts{X}$\@.  If\/ $\nbhd{T}$ is an open trapping region with
attracting-repelling pair\/ $(A,R)$\@, then\/
$R=\ts{X}\setminus\Orb^-(\nbhd{T})$\@.
\begin{proof}
If $x\in R$ and $t\in\tdomainnn$\@, then $\Phi(t,x)\in R$ by
Proposition~\ref{prop:attract-repell-props}\eqref{pl:attract-repell-props3}\@.
Since $R\subset\ts{X}\setminus\nbhd{T}$\@, $\Phi(t,x)\not\in\nbhd{T}$ and so
$x\not\in\Orb^-(\nbhd{T})$\@.  Next let $x\not\in R$ and let
$t\in\tdomainnn$\@.  Since $\Phi$ is a flow (and so $\Phi_t$ is an
homeomorphism for all $t\in\tdomain$) and since $\ts{X}\setminus\nbhd{T}$ is
closed, we have
\begin{equation*}
R=\bigcap_{s\in\tdomainnp}\Phi(\tdomain_{\le{}s}\times(\ts{X}\setminus\nbhd{T})).
\end{equation*}
Thus $x\not\in\Phi_{-t}(\ts{X}\setminus\nbhd{T})$\@.  However, $\ts{X}$ is
the disjoint union of $\Phi_{-t}(\ts{X}\setminus\nbhd{T})$ and
$\Phi_{-t}(\nbhd{T})$ (again since $\Phi$ is a flow).  Therefore, $x\in\Phi_{-t}(\nbhd{T})\subset\Orb^-(\nbhd{T})$\@.
\end{proof}
\end{lemma}

\section{Chains, chain recurrence, and chain
equivalence}\label{sec:chain-stuff}

In this section we define chains and associated notions such as chain
recurrence and chain equivalence.  We begin by enumerating useful properties
of so-called error functions, and then, after introducing the definitions for
chains and related notions, give the proof of \citet{MH:95} of equivalent
characterisations of chain equivalence.  \citeauthor{MH:95} actually
characterises chain \emph{recurrence}\@, but the proofs can be adapted to
chain equivalence.  Moreover, we shall see that these characterisations of
chain equivalence are essential to our proof of the existence of complete
Lyapunov functions for continuous-time flows and semiflows.

\subsection{Error functions}\label{subsec:error-functions}

The important observation of \cite{MH:92} was that one can replace the
constant $\epsilon$'s in the usual definition of chain recurrence (see the
definitions below) with positive continuous functions.  The use of
nonconstant functions, in combination with the metric, is reminiscent of the
construction of the so-called fine topology for the space of continuous
functions on a metric space~\cite{RAM/SK/VJ:18}\@.  In this section we
collect a few useful technical results for positive continuous functions.

First we give a descriptive name to the set positive continuous functions.
By $\mappings[0]{\ts{A}}{\ts{B}}$\@, we denote the space of continuous
functions from the topological space $\ts{A}$ to the topological space
$\ts{B}$\@.
\begin{definition}
For a topological space $(\ts{X},\sO)$\@, an \defn{error function} is an
element of $\mappings[0]{\ts{X}}{\realp}$\@.\oprocend
\end{definition}

Our first result is a general result concerning approximations of
semicontinuous functions by continuous function.
\begin{lemma}\label{lem:uclcsqueeze}
Let\/ $(\ts{X},\sO)$ be a paracompact topological space and let\/
$\map{\ul{f},\ol{f}}{\ts{X}}{\real}$ be upper (\resp~lower) semicontinuous
functions satisfying\/ $\ul{f}(x)<\ol{f}(x)$\@,\/ $x\in\ts{X}$\@.  Then there
exists a continuous function\/ $\map{f}{\ts{X}}{\real}$ such that
\begin{equation*}
\ul{f}(x)<f(x)<\ol{f}(x),\qquad x\in\ts{X}.
\end{equation*}
\begin{proof}
For $q\in\rational$\@, denote
\begin{equation*}
\nbhd{O}_q=\setdef{x\in\ts{X}}{\ul{f}<q}\cup
\setdef{x\in\ts{X}}{\ol{f}(x)>q},
\end{equation*}
this set being open.  Since, for each $x\in\ts{X}$\@, there exists
$q\in\rational$ such that $\ul{f}(x)<q<\ol{f}(x)$\@, it follows that
$\ifam{\nbhd{O}_q}_{q\in\rational}$ is an open cover of $\ts{X}$\@.  Let
$\ifam{\phi_q}_{q\in\rational}$ be a partition of unity subordinate to
$\ifam{\nbhd{O}_q}_{q\in\rational}$ and define
$f=\sum_{q\in\rational}q\phi_q$\@.  By local finiteness, $f$ is continuous.
Also, for $x\in\ts{X}$\@, let $q_1,\dots,q_k\in\rational$ be such that
$x\in\supp(\phi_q)$ if and only if $q\in\{q_1,\dots,q_k\}$\@.  Therefore,
$x\in\cap_{j=1}^k\nbhd{O}_{q_j}$ and so $\ul{f}(x)<q_j<\ol{f}(x)$\@,
$j\in\{1,\dots,k\}$\@.  Therefore,
\begin{equation*}
\ul{f}(x)=\ul{f}(x)\sum_{j=1}^k\phi_{q_j}(x)<
\sum_{j=1}^kq_j\phi_{q_j}(x)=f(x)<\ol{f}(x)\sum_{j=1}^k\phi_{q_j}(x)=
\ol{f}(x),
\end{equation*}
as desired.
\end{proof}
\end{lemma}

Next we give a sort of approximation lemma for functions in
$\mappings[0]{\ts{X}}{\realp}$\@.  For a metric space $(\ts{X},\d)$\@,
$\oball[\d{}]{r}{x}$ denotes the ball of redius $r\in\realp$ centred at
$x\in\ts{X}$\@.
\begin{lemma}\label{lem:Papprox}
Let\/ $(\ts{X},\d)$ be a metric space.  Then, for\/
$\varepsilon\in\mappings[0]{\ts{X}}{\realp}$\@, there exists\/
$\delta\in\mappings[0]{\ts{X}}{\realp}$ such that
\begin{equation*}
\d(x,y)<\delta(x)\implies
\tfrac{1}{2}\varepsilon(x)<\varepsilon(y)<\tfrac{3}{2}\varepsilon(x).
\end{equation*}
\begin{proof}
First we give a $\delta$ that gives rise to the lower bound.  Let
$\map{\beta}{\ts{X}}{\realp}$ be defined by
\begin{equation*}
\beta(x)=\sup\setdef{\eta\in\realp}{\textrm{there exists}\
\alpha\in\interval({\tfrac{1}{2}},1)\ \textrm{such that}\
\d(x,y)<\eta\implies\varepsilon(y)>\alpha\varepsilon(x)}.
\end{equation*}
Fix $\alpha\in\interval({\frac{1}{2}},1)$ so that
$\alpha\varepsilon(x)<\varepsilon(x)$ for $x\in\ts{X}$\@.  By continuity of
$\varepsilon$ at $x\in\ts{X}$\@, there exists $\eta\in\realp$ be such that
\begin{equation*}
\d(x,y)<\eta\implies\varepsilon(y)<\alpha\varepsilon(x).
\end{equation*}
Therefore, $\beta(x)\ge\eta>0$\@.

We claim that $\beta$ is lower semicontinuous.  Let $x\in\ts{X}$ and let
$a<\beta(x)$\@.  By definition of $\beta$\@, let $\eta\in\realp$ be such that
$a<\eta$ and such that there exists $\alpha\in\interval({\frac{1}{2}},1)$
for which
\begin{equation*}
\d(x,y)<\eta\implies\varepsilon(y)>\alpha\varepsilon(x).
\end{equation*}
Let $\gamma\in\interval({\frac{1}{2}},\alpha)$\@.  Since $\gamma<\alpha$\@,
$\gamma^{-1}(\alpha-\gamma)\varepsilon(x)>0$\@.  Let
$\xi\in\interval(0,{\gamma^{-1}(\gamma-\alpha)\varepsilon(x)})$\@.  Let
$\nbhd{V}$ be a neighbourhood of $x$ such that
\begin{equation*}
y\in\nbhd{V}\implies\varepsilon(y)<\varepsilon+\xi.
\end{equation*}
Then we compute, for $y\in\nbhd{V}$\@,
\begin{equation}\label{eq:Papprox1}
\gamma\varepsilon(y)<\gamma\varepsilon(x)+\gamma\xi<
\gamma\varepsilon(x)+(\gamma-\alpha)\varepsilon(x)=\alpha\varepsilon(x).
\end{equation}
Let $\zeta\in\interval(a,\eta)$ and let $r\in\interval(0,{\eta-\zeta})$ be
sufficiently small that $\nbhd{U}\eqdef\oball[\d{}]{r}{x}\subset\nbhd{V}$\@.
We claim that, if $y\in\nbhd{U}$\@, then
$\oball[\d{}]{\zeta}{y}\subset\oball[\d{}]{\eta}{x}$\@.  Indeed, let
$y\in\nbhd{U}$ and let $z\in\oball[\d{}]{\zeta}{y}$\@.  Then
\begin{equation*}
\d(x,z)\le\d(x,y)+\d(y,z)\le r+\zeta<\eta-\zeta+\zeta=\eta,
\end{equation*}
as claimed.  Thus, if $y\in\nbhd{U}$ and if $z\in\oball[\d{}]{y}{\zeta}$\@,
then $\d(x,z)<\eta$ implying that $\varepsilon(z)>\alpha\varepsilon(x)$\@.
Let $y\in\nbhd{U}\subset\nbhd{V}$\@.  Then, as in~\eqref{eq:Papprox2}\@,
$\gamma\varepsilon(y)<\alpha\varepsilon(x)$\@, which immediately gives
$\varepsilon(z)>\gamma\varepsilon(y)$\@.  Since
$\gamma\in\interval(0,{\frac{1}{2}})$\@, the definition of $\beta$ implies
that we must have $\beta(y)\ge\zeta>a$\@.  This gives the claimed lower
semicontinuity of $\beta$\@.

By Lemma~\ref{lem:uclcsqueeze}\@, let $\delta\in\mappings[0]{\ts{X}}{\realp}$
be such that $\delta(x)\in\interval(0,{\beta(x)})$ for $x\in\ts{X}$\@.  The
definition of $\beta$ then implies that there exists $\eta\in\realp$ such
that $\delta(x)<\eta$ and such that
\begin{equation*}
\d(x,y)<\eta\implies\varepsilon(y)>\alpha(x)
\end{equation*}
for some $\alpha\in\interval({\frac{1}{2}},1)$\@.  Thus we have
\begin{equation*}
\d(x,y)<\delta(x)<\eta\implies\varepsilon(y)>\alpha\varepsilon(x)>
\tfrac{1}{2}\varepsilon(x),
\end{equation*}
as desired.

Now we give $\delta$ that gives rise to the upper bound.  The argument is
similar to that in the first part of the proof, but we give it for
completeness.  Let $\map{\beta}{\ts{X}}{\realp}$ be defined by
\begin{equation*}
\beta(x)=\sup\setdef{\eta\in\realp}{\textrm{there exists}\
\alpha\in\interval(1,{\tfrac{3}{2}})\ \textrm{such that}\
\d(x,y)<\eta\implies\varepsilon(y)<\alpha\varepsilon(x)}.
\end{equation*}
Let $\alpha\in\interval(1,{\frac{3}{2}})$\@.  Since $\varepsilon$ is
continuous at $x\in\ts{X}$\@, there exists $\eta\in\realp$ be such that
\begin{equation*}
\d(x,y)<\eta\implies\varepsilon(y)<\alpha\varepsilon(x).
\end{equation*}
By definition of $\beta$\@, we have $\beta(x)\ge\eta>0$\@.

We claim that $\beta$ is lower semicontinuous.  Let $x\in\ts{X}$ and let
$a<\beta(x)$\@.  By definition of $\beta$\@, let $\eta\in\realp$ be such that
$a<\eta$ and such that there exists $\alpha\in\interval(1,{\frac{3}{2}})$
for which
\begin{equation*}
\d(x,y)<\eta\implies\varepsilon(y)<\alpha\varepsilon(x).
\end{equation*}
Let $\gamma\in\interval(\alpha,{\frac{3}{2}})$\@.  Since $\gamma>\alpha$\@,
$\gamma^{-1}(\gamma-\alpha)\varepsilon(x)>0$\@.  Let
$\xi\in\interval(0,{\gamma^{-1}(\gamma-\alpha)\varepsilon(x)})$\@.  Let
$\nbhd{V}$ be a neighbourhood of $x$ such that
\begin{equation*}
y\in\nbhd{V}\implies\varepsilon(y)>\varepsilon-\xi.
\end{equation*}
Then we compute, for $y\in\nbhd{V}$\@,
\begin{equation}\label{eq:Papprox2}
\gamma\varepsilon(y)>\gamma\varepsilon(x)-\gamma\xi>
\gamma\varepsilon(x)-(\gamma-\alpha)\varepsilon(x)=\alpha\varepsilon(x).
\end{equation}
Let $\zeta\in\interval(a,\eta)$ and let $r\in\interval(0,{\eta-\zeta})$ be
sufficiently small that $\nbhd{U}\eqdef\oball[\d{}]{r}{x}\subset\nbhd{V}$\@.
We claim that, if $y\in\nbhd{U}$\@, then
$\oball[\d{}]{\zeta}{y}\subset\oball[\d{}]{\eta}{x}$\@.  Indeed, let
$y\in\nbhd{U}$ and let $z\in\oball[\d{}]{\zeta}{y}$\@.  Then
\begin{equation*}
\d(x,z)\le\d(x,y)+\d(y,z)\le r+\zeta<\eta-\zeta+\zeta=\eta,
\end{equation*}
as claimed.  Thus, if $y\in\nbhd{U}$ and if $z\in\oball[\d{}]{y}{\zeta}$\@,
then $\d(x,z)<\eta$ implying that $\varepsilon(z)<\alpha\varepsilon(x)$\@.
Let $y\in\nbhd{U}\subset\nbhd{V}$\@.  Then, as in~\eqref{eq:Papprox2}\@,
$\alpha\varepsilon(x)<\gamma\varepsilon(y)$\@, which immediately gives
$\varepsilon(z)<\gamma\varepsilon(y)$\@.  Since
$\gamma\in\interval(1,{\frac{3}{2}})$\@, the definition of $\beta$ implies
that we must have $\beta(y)\ge\zeta>a$\@.  This gives the claimed lower
semicontinuity of $\beta$\@.

By Lemma~\ref{lem:uclcsqueeze}\@, let $\delta\in\mappings[0]{\ts{X}}{\realp}$
be such that $\delta(x)\in\interval(0,{\beta(x)})$ for $x\in\ts{X}$\@.  The
definition of $\beta$ then implies that there exists $\eta\in\realp$ such
that $\delta(x)<\eta$ and such that
\begin{equation*}
\d(x,y)<\eta\implies\varepsilon(y)<\alpha(x)
\end{equation*}
for some $\alpha\in\interval(1,{\frac{3}{2}})$\@.  Thus we have
\begin{equation*}
\d(x,y)<\delta(x)<\eta\implies\varepsilon(y)<\alpha\varepsilon(x)<
\tfrac{3}{2}\varepsilon(x),
\end{equation*}
as desired.

By choosing the min of the $\delta$'s giving rise to the lower and upper
bound, we can ensure a $\delta$ that simultaneously gives rise to both
bounds.
\end{proof}
\end{lemma}

Next we show how a positive function behaves under continuous maps.  The
proof has an entirely similar flavour to the preceding lemma.
\begin{lemma}\label{lem:general-epsilon-delta}
Let\/ $(\ts{X},\d)$ and\/ $(\ts{X}',\d')$ be metric spaces and let\/
$\phi\in\mappings[0]{\ts{X}}{\ts{X}'}$\@.  Then, for\/
$\varepsilon\in\mappings[0]{\ts{X}'}{\realp}$\@, there exists\/
$\delta\in\mappings[0]{\ts{X}}{\realp}$ such that
\begin{equation*}
\d(x,y)<\delta(x)\implies\d'(\phi(x),\phi(y))<\varepsilon(\phi(x)).
\end{equation*}
\begin{proof}
For $x\in\ts{X}$\@, define
\begin{equation*}
\beta(x)=\sup\setdef{\eta\in\realp}{\textrm{there exists}\
\alpha\in\interval(0,{\varepsilon(\phi(x))})\ \textrm{such that}\ \phi(\oball[\d{}]{\eta}{x})\subset\oball[\d']{\alpha}{\phi(x)}},
\end{equation*}
this making sense by continuity of $\phi$\@.  Let $x\in\ts{X}$ and let
$\alpha\in\interval(0,{\varepsilon(\phi(x))})$\@.  Then there exists
$\eta\in\realp$ such that
\begin{equation*}
\phi(\oball[\d{}]{\eta}{x}\subset\oball[\d']{\alpha}{\phi(x)},
\end{equation*}
and so we can conclude that $\beta(x)\ge\eta>0$\@.

Now we shall show that $\beta$ is lower semicontinuous.  Let $x\in\ts{X}$
and let $a<\beta(x)$\@.  Let $\eta\in\realp$ satisfy $a<\eta$ and
\begin{equation*}
\phi(\oball[\d{}]{\eta}{x})\subset\oball[\d']{\alpha}{\phi(x)}
\end{equation*}
for some $\alpha\in\interval(0,{\varepsilon(\phi(x))})$\@.  Let
$\gamma\in\interval(\alpha,{\varepsilon(\phi(x))})$\@.  By continuity of
$\phi$ and $\varepsilon\scirc\phi$\@, let $\nbhd{U}$ be a neighbourhood of
$x$ such that
\begin{compactenum}
\item $\d(\phi(x),\phi(y))<\gamma-\alpha$ and
\item $\varepsilon(\phi(y))>\gamma$
\end{compactenum}
for $y\in\nbhd{U}$\@.  Let $y\in\nbhd{U}$\@.  Note that
\begin{equation*}
\d(z,\phi(x))<\alpha\implies\d(z,\phi(y))\le\d(\phi(x),\phi(y))+\d(\phi(x),z)
<\gamma-\alpha+\alpha=\gamma,
\end{equation*}
whence $\oball[\d']{\alpha}{\phi(x)}\subset\oball[\d']{\gamma}{\phi(y)}$\@.
Next choose $\xi\in\interval(a,\eta)$ and choose
$b\in\interval(0,{\eta-\xi})$ so that $\oball[\d{}]{b}{x}\subset\nbhd{U}$\@.
For $y\in\nbhd{U}$ we have
\begin{equation*}
\d(z,y)<\xi\implies\d(z,x)\le\d(x,y)+\d(y,z)<b+\xi<\eta-\xi+\xi=\eta,
\end{equation*}
whereupon $\oball[\d{}]{\xi}{y}\subset\oball[\d{}]{\eta}{x}$\@.  We also have
\begin{equation*}
\phi(\oball[\d{}]{\xi}{y})\subset\phi(\oball[\d{}]{\eta}{x})
\subset\oball[\d']{\alpha}{\phi(x)}\subset\oball[\d']{\gamma}{\phi(y)}.
\end{equation*}
From this we deduce that $\beta(y)\ge\xi>a$\@, which gives lower
semicontinuity of $\beta$\@.

By Lemma~\ref{lem:uclcsqueeze}\@, let $\delta\in\mappings[0]{\ts{X}}{\realp}$
be such that $\delta(x)\in\interval(0,{\beta(x)})$ for $x\in\ts{X}$\@.  Let
$x\in\ts{X}$\@.  Since $\delta(x)<\beta(x)$\@, there exists $\eta\in\realp$
such that $\delta(x)<\eta$ and
\begin{equation*}
\phi(\oball[\d{}]{\eta}{x})\subset\oball[\d']{\alpha}{\phi(x)}
\end{equation*}
for some $\alpha\in\interval(0,{\varepsilon(\phi(x))})$\@.  Then we have
\begin{equation*}
\d(x,y)<\delta(x)\implies y\in\oball[\d{}]{\eta}{x},
\end{equation*}
which in turn gives
\begin{equation*}
\phi(y)\in\phi(\oball[\d{}]{\eta}{x})\subset\oball[\d']{\alpha}{\phi(x)}
\implies\d'(\phi(y),\phi(x))<\alpha<\varepsilon(\phi(x)),
\end{equation*}
as desired.
\end{proof}
\end{lemma}

Our final technical lemma that is of interest to us is the following.
\begin{lemma}\label{lem:PPhiestimate}
Let\/ $(\ts{X},\d)$ be a metric space and let\/ $\Phi$ be a topological flow
or semiflow on\/ $\ts{X}$\@.  Let\/
$\varepsilon\in\mappings[0]{\ts{X}}{\realp}$ and\/ $T\in\tdomainp$\@.  Then
there exists\/ $\delta\in\mappings[0]{\ts{X}}{\realp}$ such that
\begin{equation*}
\d(x,y)<\delta(x)\implies\d(\Phi(t,x),\Phi(t,y))<
\varepsilon(\Phi(t,x)),\qquad t\in\tdomain_{\interval[0,T]}.
\end{equation*}
\begin{proof}
We record a lemma (sometimes referred to as the integral continuity
condition) for use later in the proof.
\begin{proofsublemma}\label{sublem:metric-flow}
Let\/ $(\ts{X},\d)$ be a metric space and let\/ $\Phi$ be a topological flow
or semiflow on\/ $\ts{X}$\@.  Let\/ $x_1\in\ts{X}$ and\/
$T\in\tdomainp$\@.  Then, for\/ $\epsilon\in\realp$\@, there exists\/
$\delta\in\realp$ such that
\begin{equation*}
\d(x_1,x_2)<\delta\implies\d(\Phi(t,x_1),\Phi(t,x_2))<\epsilon,
\qquad t\in\tdomain_{\interval[0,T]}.
\end{equation*}
\begin{subproof}
For concreteness, we suppose that $\Phi$ is a semiflow.  For flows, the proof
is the same, with the only change being the domain of $\Phi$\@.

Let $\epsilon\in\realp$\@, let $x_1\in\ts{X}$\@, and let\/ $T\in\tdomainp$\@.
Note that
\begin{equation*}
(t,x_2)\mapsto\d(\Phi(t,x_1),\Phi(t,x_2))
\end{equation*}
is continuous since it is a composition of the continuous maps
\begin{align*}
&\tdomainnn\times\ts{X}\ni(t,x_2)\mapsto(t,(t,x_2))
\in\tdomainnn\times(\tdomainnn\times\ts{X}),\\
&(s,(t,x_2))\mapsto(\Phi(s,x_1),\Phi(t,x_2)),\\
&\ts{X}\times\ts{X}\ni(y_1,y_2)\mapsto\d(y_1,y_2)\in\realnn.
\end{align*}
For each $t\in\tdomain_{\interval[0,T]}$\@, there is a neighbourhood
$\nbhd{U}_t$ of $x_1$ and an open (relatively in $\tdomain$) interval $I_t$
around $t$ such that
\begin{equation*}
\d(\Phi(s,x),\Phi(s,x_1))<\epsilon,\qquad(s,x)\in I_t\times\nbhd{U}_t,
\end{equation*}
by continuity.  Note that $\ifam{I_t}_{t\in\tdomain_{\interval[0,T]}}$ is an
open cover of $\tdomain_{\interval[0,T]}$ and so there exist
$t_1,\dots,t_k\in\tdomain_{\interval[0,T]}$ such that
$\tdomain_{\interval[0,T]}\subset\cup_{t_j}I_{t_j}$\@.  Let
$\nbhd{U}=\cap_{j=1}^k\nbhd{U}_{t_j}$\@.  For $t\in\tdomain_{\interval[0,T]}$
and $x_2\in\nbhd{U}$\@, we have $(t,x_2)\in I_{t_j}\times\nbhd{U}_{t_j}$ for
some $j\in\{1,\dots,k\}$\@, whence
\begin{equation*}
\d(\Phi(t,x_1),\Phi(t,x_2))<\epsilon.
\end{equation*}
The result follows by taking $\delta$ sufficiently small that $\oball[\d{}]{\delta}{x_1}\subset\nbhd{U}$\@.
\end{subproof}
\end{proofsublemma}

For $x,y\in\ts{X}$\@, denote
\begin{equation*}
M_x=\inf\setdef{\varepsilon(\Phi(t,x))}{t\in\tdomain_{\interval[0,T]}}
\end{equation*}
and
\begin{equation*}
\rho(x,y)=\sup\setdef{\d(\Phi(t,x),\Phi(t,y))}{t\in\tdomain_{\interval[0,T]}},
\end{equation*}
and define $\map{\beta}{\ts{X}}{\realp}$ by
\begin{equation*}
\beta(x)=\sup\setdef{\eta\in\realp}{\textrm{there exists}\
\alpha\in\interval(0,1)\ \textrm{such that}\
\d(x,y)<\eta\implies\rho(x,y)<\alpha M_x}.
\end{equation*}

Let $\alpha\in\interval(0,1)$ so $\alpha M_x>0$\@.  Let $\eta\in\realp$ be
such that
\begin{equation*}
\d(x,y)<\eta\implies\d(\Phi(t,x),\Phi(t,y))<\alpha M_x,\qquad
t\in\tdomain_{\interval[0,T]};
\end{equation*}
this is possible by Sublemma~\ref{sublem:metric-flow}\@.  Note that, if
$\rho(x,y)<\alpha M_x$\@, then $\beta(x)\ge\eta>0$ by definition of
$\beta$\@.

We claim that $\beta$ is lower semicontinuous.  Let $x\in\ts{X}$ and let
$a<\beta(x)$\@.  By definition of $\beta$\@, let $\eta\in\realp$ satisfy
$a<\eta$ and satisfy
\begin{equation*}
\d(x,y)<\eta\implies\rho(x,y)<\alpha M_x
\end{equation*}
for some $\alpha\in\interval(0,1)$\@.

We claim that, for $\zeta\in\realp$\@, there exists a neighbourhood
$\nbhd{V}$ of $x$ such that
\begin{equation*}
\varepsilon(\Phi(t,x))-\zeta<\varepsilon(\Phi(t,y))<
\varepsilon(\Phi(t,x))+\zeta,\qquad
(t,y)\in\tdomain_{\interval[0,T]}\times\nbhd{V}.
\end{equation*}
To see this, note that continuity of $\varepsilon\scirc\Phi$ implies that,
for $t\in\tdomain_{\interval[0,T]}$\@, there is a neighbourhood $\nbhd{V}_t$
of $x$ and an open (relatively in $\tdomain$) set $I_t$ about $t$ such that
\begin{equation*}
\snorm{\varepsilon(\Phi(t,x))-\varepsilon(\Phi(s,y))}<\frac{\zeta}{2},
\qquad(s,y)\in I_t\times\nbhd{V}_t.
\end{equation*}
By compactness of $\tdomain_{\interval[0,T]}$\@, let
$t_1,\dots,t_k\in\tdomain_{\interval[0,T]}$ be such that
$\tdomain_{\interval[0,T]}\subset\cup_{j=1}^kI_{t_j}$\@.  Denote
$\nbhd{V}=\cap_{j=1}^k\nbhd{V}_{t_j}$\@.  Let
$(t,y)\in\tdomain_{\interval[0,T]}\times\nbhd{V}$\@.  Then $t\in I_{t_j}$ for
some $j\in\{1,\dots,k\}$\@.  We also have $x,y\in\nbhd{V}_{t_j}$\@.  Then
\begin{align*}
&\snorm{\varepsilon(\Phi(t_j,x))-\varepsilon(\Phi(t,x))},
\snorm{\varepsilon(\Phi(t_j,x))-\varepsilon(\Phi(t,y))}<\frac{\zeta}{2}\\
\implies&\snorm{\varepsilon(\Phi(t,y))-\varepsilon(\Phi(t,x))}<\zeta,
\end{align*}
establishing our claim.

Let $\gamma,\kappa\in\realp$ satisfy
\begin{equation*}
0<\alpha<\gamma<\kappa<1,\quad\frac{\kappa-\gamma}{\kappa}M_x>0.
\end{equation*}
As per the preceding paragraph, there exists a neighbourhood $\nbhd{V}$ of
$x$ such that, for $(t,y)\in\tdomain_{\interval[0,T]}\times\nbhd{V}$\@,
\begin{align*}
\varepsilon(\Phi(t,y))>&\;\varepsilon(\Phi(t,x))-
\frac{\kappa-\gamma}{\kappa}M_x\\
\ge&\;M_x-\frac{\kappa-\gamma}{\kappa}M_x=\frac{\gamma}{\kappa}M_x.
\end{align*}
This gives the inequalities
\begin{equation}\label{eq:PPhiestimate1}
M_y>\frac{\gamma}{\kappa}M_x\implies\gamma M_x<\theta M_y,\qquad
(\gamma-\alpha)M_x>0
\end{equation}
for $y\in\nbhd{V}$\@.

Since $\gamma-\alpha\in\interval(0,1)$\@, let $\xi\in\realp$ be such that
\begin{equation*}
\d(x,y)<\xi\implies\rho(x,y)<(\gamma-\alpha)M_x;
\end{equation*}
as above, this is possible by Sublemma~\ref{sublem:metric-flow}\@.

Let $\zeta\in\realp$ be such that, if $0<r<\min\{\eta-\zeta,\xi\}$\@, then
$\nbhd{U}\eqdef\oball[\d{}]{r}{x}\subset\nbhd{V}$\@.  We claim that
$\oball[\d{}]{\zeta}{y}\subset\oball[\d{}]{\eta}{x}$ if $y\in\nbhd{U}$\@.
Indeed, if $\d(y,z)<\zeta$\@, then
\begin{equation*}
\d(x,z)\le\d(x,y)+\d(y,z)<\eta-\zeta+\zeta=\eta,
\end{equation*}
as claimed.  Thus, if $y\in\nbhd{U}$ and $z\in\oball[\d{}]{\zeta}{y}$\@, then
\begin{equation*}
\d(x,z)<\eta\implies\rho(x,z)<\alpha M_x.
\end{equation*}
If $y\in\nbhd{U}$ then
\begin{equation*}
\d(x,y)<b<\xi\implies\rho(x,y)<(\gamma-\alpha)M_x.
\end{equation*}
By~\eqref{eq:PPhiestimate1}\@, $\gamma M_x<\kappa M_y$ if $y\in\nbhd{U}$\@.
Thus, if $y\in\nbhd{U}$ and $z\in\oball[\d{}]{\zeta}{y}$\@, then
\begin{equation*}
\rho(y,z)\le\rho(y,x)+\rho(x,z)<(\gamma-\alpha)M_x+\alpha M_x=\gamma M_x<
\kappa M_y.
\end{equation*}
The definition of $\beta$ then gives $\beta(y)\ge\zeta>a$\@, giving the
desired lower semicontinuity of $\beta$\@.

By Lemma~\ref{lem:uclcsqueeze}\@, let $\delta\in\mappings[0]{\ts{X}}{\realp}$
be such that $0<\delta(x)<\beta(x)$\@, $x\in\ts{X}$\@.  The definition of
$\beta$ implies that, for a given $x\in\ts{X}$\@, there exists
$\eta\in\realp$ such that $\eta>\delta(x)$ such that
\begin{equation*}
\d(x,y)<\eta\implies\rho(x,y)<\alpha M_x
\end{equation*}
for some $\alpha\in\interval(0,1)$\@.  The definitions of $\rho$ and $M_x$
give
\begin{equation*}
\d(\Phi(t,x),\Phi(t,y))\le\rho(x,y)<\alpha M_x<
\alpha\varepsilon(\Phi(t,x)),\qquad t\in\tdomain_{\interval[0,T]},
\end{equation*}
which is the desired conclusion.
\end{proof}
\end{lemma}

\subsection{Chain recurrence}

We next introduce chains and chain recurrence.  The notion of chain we use is
the following.
\begin{definition}
Let $(\ts{X},\d)$ be a metric space and let $\Phi$ be a topological flow or
semiflow on $\ts{X}$\@.  For $x,y\in\ts{X}$\@,
$\varepsilon\in\mappings[0]{\ts{X}}{\realp}$\@, and $T\in\tdomainp$\@, an
\defn{$(\varepsilon,T)$-chain} for $\Phi$ from $x$ to $y$ is two finite
sequences
\begin{equation*}
x_0,x_1,\dots,x_k,\quad t_0,t_1,\dots,t_{k-1}
\end{equation*}
with
\begin{compactenum}[(i)]
\item $x_0,x_1,\dots,x_k\in\ts{X}$\@,
\item $t_0,t_1,\dots,t_{k-1}\in\tdomain_{\ge{}T}$\@,
\item $x_0=x$ and $x_k=y$\@, and
\item $\d(\Phi(t_j,x_j),x_{j+1})<\varepsilon(\Phi(t_j,x_j))$\@, $j\in\{0,1,\dots,k-1\}$\@.
\end{compactenum}
For such an $(\varepsilon,T)$-chain, its \defn{length} is $k$\@.  An
\defn{$\varepsilon$-$T$-chain} is an $(\varepsilon,T)$ chain
\begin{equation*}
x_0,x_1,\dots,x_k,\quad t_0,t_1,\dots,t_{k-1}
\end{equation*}
for which $t_j=T$\@, $j\in\{0,1,\dots,k-1\}$\@.\oprocend
\end{definition}

\begin{figure}[htbp]
\centering
\includegraphics[width=0.8\hsize]{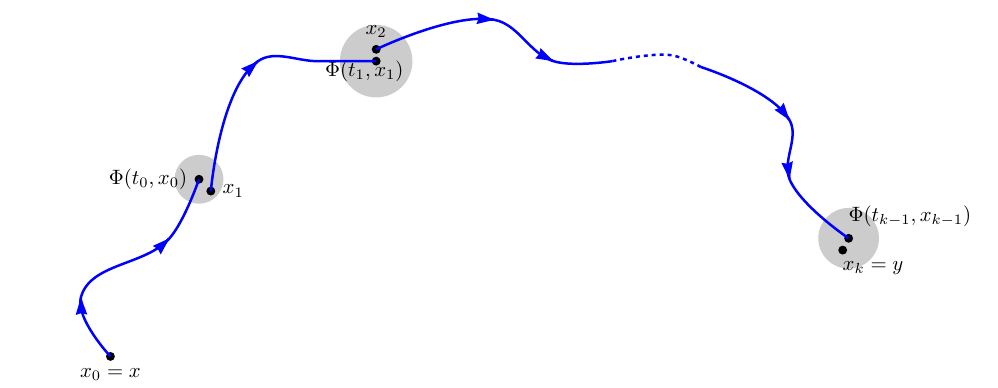}
\caption{An $(\varepsilon,T)$-chain (the depiction is of the continuous-time
case)}\label{fig:chain}
\end{figure}

In Figure~\ref{fig:chain}
we depict a chain.  We can then define an associated notion of recurrence.
\begin{definition}
Let $(\ts{X},\d)$ be a metric space and let $\Phi$ be a topological flow or
semiflow on $\ts{X}$\@.
\begin{compactenum}[(i)]
\item A point $x\in\ts{X}$ is \defn{chain recurrent} for $\Phi$ if, for each
$\varepsilon\in\mappings[0]{\ts{X}}{\realp}$ and $T\in\tdomainp$\@, there
exists an $(\varepsilon,T)$-chain from $x$ to itself.
\item Let $T\in\tdomainp$\@.  A point $x\in\ts{X}$ is \defn{$T$-chain
recurrent} for $\Phi$ if, for each
$\varepsilon\in\mappings[0]{\ts{X}}{\realp}$\@, there exists an
$(\varepsilon,T)$-chain from $x$ to itself.
\item Let $T\in\tdomainp$\@.  A point $x\in\ts{X}$ is \defn{exactly $T$-chain
recurrent} for $\Phi$ if, for each
$\varepsilon\in\mappings[0]{\ts{X}}{\realp}$\@, there exists an
$\varepsilon$-$T$-chain from $x$ to itself.
\item We denote by $\ChRec(\Phi)$ the set of chain recurrent points for
$\Phi$\@.
\item We denote by $\ChRec_{\ge{}T}(\Phi)$ the set of $T$-chain recurrent
points for $\Phi$\@.
\item We denote by $\ChRec_{=T}(\Phi)$ the set of exactly $T$-chain recurrent
points for $\Phi$\@.\oprocend
\end{compactenum}
\end{definition}

In the original definition of chain recurrence by \citet{CC:78}\@, constant
error functions are used.  For compact spaces, it is easy to see that the
definitions for constant and nonconstant error functions are equivalent.  For
noncompact spaces, however, they are not the same, as the following example
shows.
\begin{example}\label{eg:metric-dependent}
Let $\ts{X}=\real\times\realp$ be the upper half-plane.  We let $\d$ be the
standard Euclidean metric for $\ts{X}$ inherited from $\real^2$ and we let
$\subscr{\d}{h}$ be the so-called hyperbolic metric,~\ie~that metric whose
geodesics are arcs of half-circles in $\ts{X}$ with centres on the line
$\real\times\{0\}$\@.  Note that the topology defined by these two metrics is
the same.  On $\ts{X}$\@, we consider the continuous-time flow defined by
$\Phi(t,(x,y))=(x+t,y)$\@.  Let $\ChRec'(\Phi)$ be the chain recurrent set
where chains are defined using constant error functions and using the metric
$\d$\@.  Let $\subscr{\ChRec'}{h}(\Phi)$ be the chain recurrent set where
chains are defined using constant error functions and using the metric
$\subscr{\d}{h}$\@.  It is not difficult to show that
$\ChRec'(\Phi)=\emptyset$ and that $\subscr{\ChRec'}{h}(\Phi)=\ts{X}$\@.  The
idea in proving the second of these formulae is to observe that horizontal
distances in the hyperbolic metric become much larger than their Euclidean
counterparts as one gets close (in the Euclidean sense) to
$\real\times\{0\}$\@.  The details are given in
\cite[Example~2.7.13]{JMA/GSN:07}\@.\oprocend
\end{example}

The example illustrates the dangers of using constant error functions in
defining chains.  Specifically, chain recurrence becomes a metric-dependent
concept in this case.  There is not necessarily anything wrong with this, of
course; perhaps one is content to have notions that depend on metric,~\eg~on
normed vector spaces.  However, the Conley decomposition relates chain
recurrence to attracting and repelling sets, and these latter are purely
topological concepts.  Thus one cannot expect that the Conley decomposition
is valid when chains are defined using constant error functions.
Nonetheless, the use of constant error functions to define chains is common
in the literature, even for noncompact states spaces.

Let us enumerate some properties of the chain recurrent set.  We begin by
relating chain recurrence to one of the classical notions of recurrence, that
of a nonwandering point.
\begin{definition}
Let $(\ts{X},\sO)$ be a topological space and let $\Phi$ be a topological
flow or semiflow.  A point $x_0\in\ts{X}$ is \defn{nonwandering} for $\Phi$
if there exists $T\in\tdomainp$ such that, for each neighbourhood $\nbhd{U}$
of $x_0$\@, $\Phi_t(\nbhd{U})\cap\nbhd{U}\not=\emptyset$ for some
$t\in\tdomain_{\ge{}T}$\@.  We denote by $\NWnd(\Phi)$ the set forward
nonwandering points for $\Phi$\@.\oprocend
\end{definition}

The following characterisation of the set of nonwandering points will be
useful.
\begin{lemma}\label{lem:NWnd-1stcountable}
Let\/ $(\ts{X},\sO)$ be a first countable Hausdorff topological space and
let\/ $\Phi$ be a topological flow or semiflow.  Let\/ $x\in\NWnd(\Phi)$\@.
Then, for a neighbourhood\/ $\nbhd{U}$ of\/ $x$ and for\/ $T\in\tdomainp$\@,
there exists\/ $t\in\tdomain_{\ge{}T}$ such that\/
$\nbhd{U}\cap\Phi_t(\nbhd{U})\not=\emptyset$\@.
\begin{proof}
We prove the contrapositive.  Thus, we let $x\in\ts{X}$\@, we assume that
there exists a neighbourhood $\nbhd{U}$ of $x$ and $T\in\tdomainp$ such that
$\nbhd{U}\cap\Phi_t(\nbhd{U})=\emptyset$ for $t\in\tdomain_{\ge{}T}$\@, and
we prove that $x\not\in\NWnd(\Phi)$\@.

Note that the assumptions ensure that $x$ is not a periodic point.  We claim
that this implies that, for any $S\in\tdomainp$\@, there exists a
neighbourhood $\nbhd{V}$ of $x$ such that
$\nbhd{V}\cap\Phi_t(\nbhd{V})=\emptyset$ for all $t\in\interval[S,T]$\@.  We
prove this by the contrapositive.  Thus we let $x\in\ts{X}$ and
$S\in\tdomainp$\@, we assume that, for any neighbourhood $\nbhd{V}$ of $x$\@,
there exists $t\in\interval[S,T]$ such that
$\nbhd{V}\cap\Phi_t(\nbhd{V})\not=\emptyset$\@, and we prove that $x$ is a
periodic point.  By first countability of $(\ts{X},\sO)$\@, let
$\ifam{\nbhd{U}_j}_{j\in\integerp}$ be a neighbourhood base for $x$\@.  Our
hypotheses ensure that, for each $j\in\integerp$\@, there exists
$x_j\in\ts{X}$ and $t_j\in\interval[S,T]$ such that
\begin{equation*}
x_j\in\nbhd{U}_j,\ \Phi(t_j,x_j)\in
\nbhd{V}_j\cap\Phi_{t_j}(\nbhd{V}_j),\qquad j\in\integerp.
\end{equation*}
In particular, $\lim_{j\to\infty}x_j=x$\@.  Since the sequence of times
$\ifam{t_j}_{j\in\integerp}$ resides in the compact interval
$\interval[S,T]$\@, there exists a subsequence
$\ifam{t_{j_k}}_{k\in\integerp}$ that converges to some
$\tau\in\interval[S,T]$\@.  Then we have
\begin{equation*}
\Phi(\tau,x)=\lim_{j_k\to\infty}\Phi(t_{j_k},x_{j_k})=x.
\end{equation*}
This shows that $x$ is a periodic point.

Now, combining the first two paragraphs of the proof, for any
$S\in\tdomainp$\@, there exist neighbourhoods $\nbhd{U}$ and $\nbhd{V}$ of
$x$ such that
\begin{equation*}
(\nbhd{U}\cap\nbhd{V})\cap\Phi_t(\nbhd{U}\cap\nbhd{V})=\emptyset
\end{equation*}
for all $t\in\tdomain_{\ge{}S}$\@.  That is to say, $x\not\in\NWnd(\Phi)$\@.
\end{proof}
\end{lemma}

We now show that nonwandering points are chain recurrent.
\begin{proposition}
Let\/ $(\ts{X},\d)$ be a metric space and let\/ $\Phi$ be a topological flow
or semiflow.  Then\/ $\NWnd(\Phi)\subset\ChRec(\Phi)$\@.
\begin{proof}
Let $\varepsilon\in\mappings[0]{\ts{X}}{\realp}$ and let $T\in\tdomainp$\@.
Let $x\in\NWnd(\Phi)$\@.  Since $\Phi_T$ and $\varepsilon$ are continuous,
let $\delta\in\interval(0,{\tfrac{1}{2}\varepsilon(x)})$ be such that
\begin{equation*}
\d(x,y)<\delta\implies\d(\Phi(T,x),\Phi(T,y))<
\varepsilon(\Phi(T,x)).
\end{equation*}
and
\begin{equation*}
\d(x,y)<\delta\implies\snorm{\varepsilon(y)-\varepsilon(x)}<
\tfrac{1}{2}\varepsilon(x).
\end{equation*}
By Lemma~\ref{lem:NWnd-1stcountable}\@, noting that metric spaces are first
countable and Hausdorff, let $t_1\in\tdomain_{>2T}$ be such that
$\oball[\d{}]{\delta}{x}\cap\Phi_{t_1}(\oball[\d{}]{\delta}{x})\not=
\emptyset$\@, and let $y\in\oball[\d{}]{\delta}{x}$ be such that
$\Phi(t_1,y)\in\oball[\d{}]{\delta}{x}$\@.  Then
\begin{equation*}
\d(x,y)<\delta\implies\d(\Phi(T,x),\Phi(T,y))<
\varepsilon(\Phi(T,x)).
\end{equation*}
Also,
\begin{equation*}
\d(\Phi(t_1-T,\Phi(T,y)),x)=\d(\Phi(t_1,y),x)<\delta<
\tfrac{1}{2}\varepsilon(x).
\end{equation*}
Thus,
\begin{equation*}
\varepsilon(\Phi(t_1,y))-\varepsilon(x)>-\tfrac{1}{2}\varepsilon(x)
\implies\varepsilon(\Phi(t_1,y))>\tfrac{1}{2}\varepsilon(x),
\end{equation*}
and so
\begin{equation*}
\d(\Phi(t_1-T,\Phi(T,y)),x)<\varepsilon(\Phi(t_1,y)).
\end{equation*}
We conclude, then, that
\begin{equation*}
x,\Phi(T,y),x,\quad T,t_1-T
\end{equation*}
is an $(\varepsilon,T)$-chain, and so $x\in\ChRec(\Phi)$\@.
\end{proof}
\end{proposition}

The invariance and closedness of the chain recurrent set will be needed in
our results below.
\begin{proposition}\label{prop:ChRec-props}
Let\/ $(\ts{X},\d)$ be a metric space and let\/ $\Phi$ be a topological flow
(\resp~semiflow) on\/ $\ts{X}$\@.  Then the following statements hold:
\begin{compactenum}[(i)]
\item \label{pl:ChRec1} $\ChRec(\Phi)$ is invariant (\resp~forward-invariant)
for\/ $\Phi$\@;
\item \label{pl:ChRec2} $\ChRec(\Phi)$ is closed;
\end{compactenum}
\begin{proof}
\eqref{pl:ChRec1} Let $x_0\in\ChRec(\Phi)$ and let $t\in\tdomain$\@.  Let
$\varepsilon\in\mappings[0]{\ts{X}}{\realp}$ and let $T\in\tdomainp$\@.

Assume first that $t\in\tdomainp$\@.  Let $T'=T+t$\@.  Let $\nbhd{U}$ be a
neighbourhood of $\Phi(t,x_0)$\@, chosen so that
\begin{equation*}
y\in\nbhd{U}\implies\begin{cases}
\d(y,\Phi(t,x_0))<\tfrac{1}{2}\varepsilon(\Phi(t,x_0)),&\\
\varepsilon(y)>\tfrac{1}{2}\varepsilon(\Phi(t,x_0)).&\end{cases}
\end{equation*}
Let $\delta\in\realp$ be small enough that
\begin{equation*}
\d(x_0,x)<\delta\implies\Phi(t,x)\in\nbhd{U}.
\end{equation*}
Take $\eta\in\mappings[0]{\ts{X}}{\realp}$ to be
$\eta(x)=\min\{\varepsilon(x),\delta\}$\@, $x\in\ts{X}$\@.  Let
\begin{equation*}
x_0,x_1,\dots,x_n=x_0,\quad t_0,t_1,\dots,t_{n-1}
\end{equation*}
be an $(\eta,T')$-chain and denote
\begin{align*}
&y_0=\Phi(t,x_0),y_1=x_1,\dots,y_{n-1}=x_{n-1},y_n=\Phi(t,x_0),\\
&s_0=t_0-t,s_1=t_1,\dots,s_{n-2}=t_{n-2},s_{n-1}=t_{n-1}+t.
\end{align*}
Note that
\begin{equation*}
\d(\Phi(s_0,y_0),y_1)=\d(\Phi(t_0-t,\Phi(t,x_0)),x_1)=
\d(\Phi(t_0,x_0),x_1).
\end{equation*}
Thus we have
\begin{equation*}
\d(\Phi(s_0,y_0),y_1)<\eta(\Phi(t_0,x_0))\le
\varepsilon(\Phi(s_0,y_0)).
\end{equation*}
Clearly we have
\begin{equation*}
\d(\Phi(s_j,y_j),y_{j+1})<\varepsilon(\Phi(s_j,y_j)),\qquad
j\in\{1,\dots,n-2\}.
\end{equation*}
Finally, we have
\begin{equation*}
\d(\Phi(s_{n-1},y_{n-1}),y_n)=\d(\Phi(t_{n-1}+t,y_{n-1}),y_n)=
\d(\Phi(t,\Phi(t_{n-1},x_{n-1})),\Phi(t,x_0)).
\end{equation*}
From this and the various definitions we calculate
\begin{align*}
&\d(\Phi(t_{n-1},x_{n-1}),x_0)<\eta(\Phi(t_{n-1},x_{n-1}))<\delta\\
\implies&\Phi(t,\Phi(t_{n-1},x_{n-1}))\in\nbhd{U}\\
\implies&\begin{cases}
\d(\Phi(t,\Phi(t_{n-1},x_{n-1})),\Phi(t,x_0))<
\tfrac{1}{2}\varepsilon(\Phi(t,x_0)),&\\
\varepsilon(\Phi(t,\Phi(t_{n-1},x_{n-1})))>
\tfrac{1}{2}\varepsilon(\Phi(t,x_0))&\end{cases}\\
\implies&\d(\Phi(t,\Phi(t_{n-1},x_{n-1})),\Phi(t,x_0))<
\varepsilon(\Phi(t,\Phi(t_{n-1},x_{n-1})))\\
\implies&\d(\Phi(s_{n-1},y_{n-1}),y_n)<
\varepsilon(\Phi(t,\Phi(t_{n-1},x_{n-1}))).
\end{align*}
Therefore,
\begin{equation*}
y_0,y_1,\dots,y_n=y_0,\quad s_0,s_1,\dots,s_{n-1}
\end{equation*}
is an $(\varepsilon,T)$-chain from $\Phi(t,x)$ to $\Phi(t,x)$\@.

If $\Phi$ is a flow and if $t\in\tdomainn$\@, then let $T'=T-t$\@.  One can
then proceed exactly as in the preceding paragraph to deduce that
$\Phi(t,x_0)\in\ChRec(\Phi)$\@.

\eqref{pl:ChRec2} Let $x_0\in\closure(\ChRec(\Phi))$\@, and let
$\varepsilon\in\mappings[0]{\ts{X}}{\realp}$ and $T\in\tdomainp$\@.  Let
$\delta\in\realp$ be such that
\begin{equation*}
\d(x_0,x)<\delta\implies\begin{cases}\d(\Phi(T,x_0),\Phi(T,x))<
\tfrac{1}{2}\varepsilon(\Phi(T,x_0)),&\\
\varepsilon(x)>\tfrac{1}{2}\varepsilon(x_0).\end{cases}
\end{equation*}
Since $x_0\in\closure(\ChRec(\Phi))$\@, let $x\in\ChRec(\Phi)$ be
such that
\begin{equation*}
\d(x_0,x)<\tfrac{1}{2}\min\{\delta,\varepsilon(\Phi(T,x_0)),
\tfrac{1}{2}\varepsilon(x_0)\}.
\end{equation*}
Let $\eta\in\mappings[0]{\ts{X}}{\realp}$ be defined by
$\eta(x)=\frac{1}{2}\min\{\delta,\varepsilon(x)\}$\@.

First note that
\begin{equation*}
x_0,x_1=\Phi(T,x),\quad T
\end{equation*}
is an $(\varepsilon,T)$-chain from $x_0$ to $\Phi(T,x)$\@.

Now, since $x\in\ChRec(\Phi)$\@, let
\begin{equation*}
y_0=x,y_1,\dots,y_n=x,\quad s_0,s_1,\dots,s_{n-1}
\end{equation*}
be an $(\eta,2T)$-chain from $x$ to itself.  Since
\begin{equation*}
\Phi(s_0-T,\Phi(T,x))=\Phi(s_0,x),
\end{equation*}
it follows immediately that
\begin{equation*}
\Phi(T,x),y_1,\dots,y_{n-1},\quad s_0-T,s_1,\dots,s_{n-2}
\end{equation*}
is an $(\varepsilon,T)$-chain from $\Phi(T,x)$ to $y_{n-1}$\@.

Finally, we claim that
\begin{equation*}
y_{n-1},x_0,\quad s_{n-1}
\end{equation*}
is an $(\varepsilon,T)$-chain from $y_{n-1}$ to $x_0$\@.  First note that,
\begin{align*}
\d(\Phi(s_{n-1},y_{n-1}),x_0)\le&\;
\d(\Phi(s_{n-1},y_{n-1}),x)+\d(x,x_0)\\
<&\;\frac{\delta}{2}+\frac{\delta}{2}=\delta.
\end{align*}
Therefore,
\begin{equation*}
\varepsilon(\Phi(s_{n-1},y_{n-1}))>\tfrac{1}{2}\varepsilon(x_0),
\end{equation*}
and so
\begin{align*}
\d(\Phi(s_{n-1},y_{n-1}),x_0)\le&\;
\d(\Phi(s_{n-1},y_{n-1}),x)+\d(x,x_0)\\
<&\;\tfrac{1}{2}\varepsilon(\Phi(s_{n-1},y_{n-1}))+
\tfrac{1}{4}\varepsilon(x_0)\\
<&\;\varepsilon(\Phi(s_{n-1},y_{n-1})),
\end{align*}
giving our claim.

Putting the above three constructions together, we see that
\begin{equation*}
x_0,\Phi(T,x),y_1,\dots,y_{n-1},x_0,\quad
T,s_0-T,s_1,\dots,s_{n-1}
\end{equation*}
is an $(\varepsilon,T)$-chain from $x_0$ to itself, whence $x_0\in\ChRec(\Phi)$\@.
\end{proof}
\end{proposition}

\subsection{Chain equivalence}

On the set of chain recurrent points of a topological flow or semiflow, there
is an important equivalence relation.
\begin{definition}
Let $(\ts{X},\d)$ be a metric space and let $\Phi$ be a topological flow or
semiflow on $\ts{X}$\@.
\begin{compactenum}[(i)]
\item Points $x,y\in\ts{X}$ are \defn{chain equivalent} for $\Phi$ if, for
each $\varepsilon\in\mappings[0]{\ts{X}}{\realp}$ and for each
$T\in\tdomainp$\@, there exist $(\varepsilon,T)$-chains
\begin{equation*}
x_0=x,x_1,\dots,x_k=y,\quad t_0,t_1,\dots,t_{k-1}
\end{equation*}
and
\begin{equation*}
y_0=y,y_1,\dots,y_m=x,\quad s_0,s_1,\dots,s_{m-1}.
\end{equation*}
\item Let $T\in\tdomainp$\@.  Points $x,y\in\ts{X}$ are \defn{$T$-chain
equivalent} for $\Phi$ if, for each
$\varepsilon\in\mappings[0]{\ts{X}}{\realp}$\@, there exist
$(\varepsilon,T)$-chains
\begin{equation*}
x_0=x,x_1,\dots,x_k=y,\quad t_0,t_1,\dots,t_{k-1}
\end{equation*}
and
\begin{equation*}
y_0=y,y_1,\dots,y_m=x,\quad s_0,s_1,\dots,s_{m-1}.
\end{equation*}
\item Let $T\in\tdomainp$\@.  Points $x,y\in\ts{X}$ are \defn{exactly
$T$-chain equivalent} for $\Phi$ if, for each
$\varepsilon\in\mappings[0]{\ts{X}}{\realp}$\@, there exist
$\varepsilon$-$T$-chains
\begin{equation*}
x_0=x,x_1,\dots,x_k=y,\quad t_0,t_1,\dots,t_{k-1}
\end{equation*}
and
\begin{equation*}
y_0=y,y_1,\dots,y_m=x,\quad s_0,s_1,\dots,s_{m-1}.
\end{equation*}
\item A \defn{chain component} for $\Phi$ is a subset $C\subset\ts{X}$ such
that, if $x,y\in C$\@, then $x$ and $y$ are chain equivalent.
\item Let $T\in\tdomainp$\@.  A \defn{$T$-chain component} for $\Phi$ is a
subset $C\subset\ts{X}$ such that, if $x,y\in C$\@, then $x$ and $y$ are
$T$-chain equivalent.
\item Let $T\in\tdomainp$\@.  An \defn{exact $T$-chain component} for $\Phi$
is a subset $C\subset\ts{X}$ such that, if $x,y\in C$\@, then $x$ and $y$ are
exactly $T$-chain equivalent.\oprocend
\end{compactenum}
\end{definition}

Note that the relation of chain equivalence is not generally an equivalence
relation on $\ts{X}$\@; it is symmetric and transitive, but not generally
reflexive.  However, chain equivalence is an equivalence relation on
$\ChRec(\Phi)$\@.  Therefore, $\ChRec(\Phi)$ is a disjoint union of chain
components.  This fact features significantly in the Fundamental Theorem of
Dynamical Systems.

Let us prove a few essential properties of chain components.
\begin{proposition}\label{prop:chain-comp}
Let\/ $(\ts{X},\d)$ be a metric space and let\/ $\Phi$ be a topological flow
(\resp~semiflow) on\/ $\ts{X}$\@.  Then the following statements hold:
\begin{compactenum}[(i)]
\item \label{pl:chain-comp1} if\/ $C$ is a chain component for\/ $\Phi$, then
it is closed;
\item \label{pl:chain-comp2} if\/ $C$ is a chain component for\/ $\Phi$\@,
then it is invariant (\resp~forward-invariant).
\end{compactenum}
\begin{proof}
\eqref{pl:chain-comp1} Let $C\subset\ChRec(\Phi)$ be a chain component.  To
prove closedness of $C$\@, we shall show that, if $y\in\closure(C)$\@, then
$y\in C$\@.  To do this, we will let $x\in C$\@,
$\varepsilon\in\mappings[0]{\ts{X}}{\realp}$\@, and $T\in\tdomainp$ and
construct two $(\varepsilon,T)$-chains, one from $x$ to $y$ and one from $y$
to $x$\@.  Since $x\in\ChRec(\Phi)$\@, there is also an
$(\varepsilon,T)$-chain from $x$ to itself.  By concatenating chains, one
concludes that (1)~$y\in\ChRec(\Phi)$ and (2)~$y$ is chain equivalent to any
point in $C$\@.  This suffices to show that $y\in C$\@.

By \eqref{lem:general-epsilon-delta}\@, let
$\delta\in\mappings[0]{\ts{X}}{\realp}$ be such that
\begin{equation*}
\d(x_1,x_2)<\delta(x_1)\implies\d(\Phi(T,x_1),\Phi(T,x_2))<
\varepsilon(\Phi(T,x_1)).
\end{equation*}
Let $x'\in C$ be such that $\d(x',y)<\delta(y)$\@; this is possible since
$y\in\closure(C)$\@.  Since $x,x'\in C$\@, let
\begin{equation*}
x_0=x',x_1,\dots,x_k=x,\quad t_0,t_1,\dots,t_{k-1}
\end{equation*}
be an $(\varepsilon,2T)$-chain from $x'$ to $x$\@.  Note that
\begin{equation*}
\d(x',y)<\delta(y)\implies
\d(\Phi(T,x'),\Phi(T,y))<\varepsilon(\Phi(T,y)),
\end{equation*}
and from this we easily deduce that
\begin{equation*}
y_0=y,y_1=\Phi(T,x'),y_2=x_1,\dots,y_{k+1}=x_k=x,\quad
s_0=T,s_1=t_0-T,\dots,s_k=t_{k-1}
\end{equation*}
is an $(\varepsilon,T)$-chain from $y$ to $x$\@.

Now we find an $(\varepsilon,T)$-chain from $x$ to $y$\@.  By
Lemma~\ref{lem:Papprox}\@, let $\delta\in\mappings[0]{\ts{X}}{\realp}$ be
such that
\begin{equation*}
\d(x_1,x_2)<\delta(x_1)\implies
\frac{1}{2}\varepsilon(x_1)<\varepsilon(x_2)<\frac{3}{2}\varepsilon(x_1).
\end{equation*}
Let $x'\in C$ satisfy
$\d(y,x')<\max\{\delta(y),\frac{1}{6}\varepsilon(y)\}$\@, this being possible
since $y\in\closure(C)$\@.  Let
\begin{equation*}
x_0=x,x_1,\dots,x_k=x',\quad t_0,t_1,\dots,t_{k-1}
\end{equation*}
be a $(\frac{1}{2}\varepsilon,T)$-chain from $x$ to $x'$\@.  We have
\begin{align*}
\d(\Phi(t_{k-1},x_{k-1}),y)\le&\;
\d(\Phi(t_{k-1},x_{k-1}),x')+\d(y,x')\\
<&\;\frac{1}{2}\varepsilon(\Phi(t_{k-1},x_{k-1}))+
\frac{1}{6}\varepsilon(y)
\end{align*}
We have
\begin{equation*}
\d(x',y)<\delta(y)\implies\frac{1}{6}\varepsilon(y)<
\frac{1}{3}\varepsilon(x')
\end{equation*}
and
\begin{equation*}
\d(\Phi(t_{k-1},x_{k-1}),x')<\delta(\Phi(t_{k-1},x_{k-1}))\implies
\frac{1}{3}\varepsilon(x')<\frac{1}{2}\varepsilon(\Phi(t_{k-1},x_{k-1})).
\end{equation*}
Putting this all together, we have
\begin{equation*}
\d(\Phi(t_{k-1},x_{k-1}),y)<\varepsilon(\Phi(t_{k-1},x_{k-1})).
\end{equation*}
From this, we conclude that
\begin{equation*}
x_0=x,x_1,\dots,x_k=y,\quad t_0,t_1,\dots,t_{k-1}
\end{equation*}
is an $(\varepsilon,T)$-chain from $x$ to $y$\@.

\eqref{pl:chain-comp2} Let $C$ be a chain component, let $x\in C$\@, and let
$t\in\tdomainp$ (if $\Phi$ is a semiflow) and $t\in\tdomain$ (if $\Phi$ is a
flow).  As in the preceding part of the proof, to show that
$\Phi(t,x)\in C$\@, it suffices to show that, for any
$\varepsilon\in\mappings[0]{\ts{X}}{\realp}$ and $T\in\tdomainp$\@, there
exist two $(\varepsilon,T)$-chains, one from $x$ to $\Phi(t,x)$ and one from
$\Phi(t,x)$ to $x$\@.

Since $C\subset\ChRec(\Phi)$\@, let
\begin{equation*}
x_0=x,x_1,\dots,x_k=x,\quad t_0,t_1,\dots,t_{k-1}
\end{equation*}
be an $(\varepsilon,T+\snorm{t})$-chain from $x$ to itself.  Then it is
evident that
\begin{equation*}
\Phi(t,x),x_1,\dots,x_k=x,\quad t_0-t,t_1,\dots,t_{k-1}
\end{equation*}
is an $(\varepsilon,T)$-chain from $\Phi(t,x)$ to $x$\@.

To construct an $(\varepsilon,T)$-chain from $x$ to $\Phi(t,x)$\@, let
$\delta\in\mappings[0]{\ts{X}}{\realp}$ be such that
\begin{equation*}
\d(x_1,x_2)<\delta(x_1)\implies
\d(\Phi(t,x_1),\Phi(t,x_2))<\varepsilon(x_1).
\end{equation*}
Let
\begin{equation*}
x_0=x,x_1,\dots,x_k=x,\quad t_0,t_1,\dots,t_{k-1}
\end{equation*}
be a $(\min\{\varepsilon,\delta\},T+\snorm{t})$-chain from $x$ to itself;
this is possible since $x\in C\subset\ChRec(\Phi)$\@.  Since
\begin{equation*}
\d(\Phi(t_{k-1},x_{k-1}),x)<\delta(\Phi(t_{k-1},x_{k-1})),
\end{equation*}
we have
\begin{align*}
\d(\Phi(\Phi(t_{k-1}+t),x_{k-1}),\Phi(t,x))=&\;
\d(\Phi(t,\Phi(t_{k-1},x_{k-1})),\Phi(t,x))\\
<&\;\varepsilon(\Phi(t_{k-1},x_{k-1})).
\end{align*}
Thus we conclude that
\begin{equation*}
x_0=0,x_1,\dots,x_{k-1},\Phi(t,x),\quad t_0,t_1,\dots,t_{k-2},t_{k-1}+t
\end{equation*}
is an $(\varepsilon,T)$-chain from $x$ to $\Phi(t,x)$\@.
\end{proof}
\end{proposition}

\subsection{Alternative characterisations of chain equivalence}

Our objective in this section is to give two alternative characterisations of
chain equivalence, as per~\cite{MH:95}\@.  There are two principal
simplifications we consider.  First we shall show that, for chain
equivalence, it is possible to fix $T$\@, provided that $T$ possesses a
multiplicative inverse in $\tdomain$\@.  Next we shall show that the precise
choices of the times $t_0,t_1,\dots,t_{k-1}$ can also be fixed to be the same
time $T$\@, provided again that $T$ possesses a multiplicative inverse in
$\tdomain$\@.  In proving that these simplifications can be made without loss
of generality, we make extensive use of the technical results about error
functions from Section~\ref{subsec:error-functions}\@.

The result we prove is the following.
\begin{theorem}\label{the:chain-equiv}
Let\/ $(\ts{X},\d)$ be a metric space and let\/ $\Phi$ be a topological
flow or semiflow on\/ $\ts{X}$\@.  Then the following statements are
equivalent for\/ $x,y\in\ts{X}$ and for $T_0\in\tdomainp$ possessing a
multiplicative inverse in\/ $\tdomain$\@:
\begin{compactenum}[(i)]
\item \label{pl:chain-equiv1} $x$ and\/ $y$ are chain equivalent for\/
$\Phi$\@;
\item \label{pl:chain-equiv2} $x$ and\/ $y$ are\/ $T_0$-chain equivalent
for\/ $\Phi$;
\item \label{pl:chain-equiv3} $x$ and\/ $y$ are\/ exactly $T_0$-chain
equivalent for\/ $\Phi$\@.
\end{compactenum}
\begin{proof}
\eqref{pl:chain-equiv1}$\implies$\eqref{pl:chain-equiv2} This is clear from
the definitions.

\eqref{pl:chain-equiv2}$\implies$\eqref{pl:chain-equiv3} We break the proof
into two parts, first for the discrete-time case then for the continuous-time
case.

\paragraph{The discrete-time case}

The discrete-time case is almost immediate.  The requirement that $T_0$
possess a multiplicative inverse in $\tdomain$ means that $T_0=1$ in the
discrete-time case.  In this case, however, every $(\varepsilon,T)$-chain
gives rise to a $\varepsilon$-$1$-chain, simply by adding zero jumps at times
for which there is not already a jump.

\paragraph{The continuous-time case}

In this case, the requirement that $T_0\in\tdomainp$ possess a multiplicative
inverse in $\tdomain$ places no restriction on $T_0$\@.  Let $x$ and $y$
satisfy the hypotheses of~\eqref{pl:chain-equiv2}\@.  We claim that this
implies that, for each\/ $\varepsilon\in\mappings[0]{\ts{X}}{\realp}$ and for
each $\tau\in\realp$\@, there exists an $(\varepsilon,T_0)$ chain
\begin{equation*}
x_0=x,x_1,\dots,x_k=y,\quad t_0,t_1,\dots,t_{k-1}
\end{equation*}
and $N\in\integerp$ such that $t_j\in\interval[{T_0},{2T_0})$\@,
$j\in\{0,1,\dots,k-1\}$\@, and $\snorm{\sum_{j=1}^kt_j-NT_0}<\tau$\@.

To prove this, let $\varepsilon\in\mappings[0]{\ts{X}}{\realp}$ and
$\tau\in\realp$\@.  Note that the hypotheses ensure the existence of
$(\frac{1}{4}\varepsilon,T_0)$-chains
\begin{equation*}
u_0=x,u_1,\dots,u_l=y,\quad r_0,r_1,\dots,r_{l-1}
\end{equation*}
and
\begin{equation*}
v_0=y,v_1,\dots,v_m=x,\quad s_0,s_1,\dots,s_{m-1}.
\end{equation*}
First of all, we can assume that $r_a\in\interval[{T_0},{2T_0})$\@,
$a\in\{0,1,\dots,l-1\}$\@, and $s_b\in\interval[{T_0},{2T_0})$\@,
$b\in\{0,1,\dots,m-1\}$\@.  Indeed, if this is not so, then we can add jump
times with zero jumps to ensure that these conditions are met.  Let us
abbreviate
\begin{equation*}
R=r_0+r_1+\dots+r_{l-1},\quad S=s_0+s_1+\dots+s_{m-1}.
\end{equation*}
By Lemma~\ref{lem:Papprox}\@, let $\delta\in\mappings[0]{\ts{X}}{\realp}$ be
such that
\begin{equation*}
\d(z_1,z_2)<\delta(z_1)\implies
\frac{1}{2}\varepsilon(z_1)<\varepsilon(z_2).
\end{equation*}
Let $\sigma\in\realp$ be such that
\begin{equation*}
\snorm{r'_0-r_0}<\sigma\implies
\d(\Phi(r_0,u_0),\Phi(r'_0,u_0))<
\min\{\tfrac{1}{4}\varepsilon(\Phi(r_0,u_0),\delta(\Phi(r_0,u_0))\}.
\end{equation*}
Choose $r'_0\in\interval[{r_0},{r_0+\sigma})$ such that $R+S-r_0+r'_0$ is
irrational and such that $r'_0\in\interval[{T_0},{2T_0})$\@.  Then, since
\begin{equation*}
\d(\Phi(r_0,u_0),\Phi(r'_0,u_0))<\delta(\Phi(r_0,u_0)),
\end{equation*}
we have
\begin{align*}
\d(\Phi(r'_0,u_0),u_1)\le&\;\d(\Phi(r'_0,u_0),\Phi(r_0,u_0))+
\d(\Phi(r_0,u_0),u_1)\\
<&\;\frac{1}{4}\varepsilon(\Phi(r_0,u_0))+
\frac{1}{4}\varepsilon(\Phi(r_0,u_0))\\
<&\;\varepsilon(\Phi(r'_0,u_0)),
\end{align*}
and we deduce that
\begin{equation*}
u_0=x,u_1,\dots,u_l=y,\quad r'_0,r_1,\dots,r_{l-1}
\end{equation*}
is an $(\varepsilon,T_0)$-chain from $x$ to $y$\@.  Abbreviate
\begin{equation*}
R'=r'_0+r_1+\dots+r_{l-1}.
\end{equation*}
By irrationality of $R'+S$\@, there exist $M,N\in\integerp$ such that
\begin{equation*}
\snorm{R'+M(R'+S)-NT_0}<\tau.
\end{equation*}
Now build an $(\varepsilon,T_0)$ chain as follows:
\begin{multline*}
u_0=x,u_1,\dots,u_l=y=\underbrace{v_0,v_1,\dots,v_m,u_0,u_1,\dots,u_l,\dots,
v_0,v_1,\dots,v_m,u_0,u_1,\dots,u_l}_{v_0,v_1,\dots,v_m,u_0,u_1,\dots,u_l\
\textrm{repeated}\ M\ \textrm{times}},\\
r'_0,r_1,\dots,r_{l-1},\underbrace{s_0,s_1,\dots,s_{m-1},
r'_0,r_1,\dots,r_{l-1},\dots,s_0,s_1,\dots,s_{m-1},
r'_0,r_1,\dots,r_{l-1}}_{s_0,s_1,\dots,s_{m-1},r'_0,r_1,\dots,r_{l-1}\
\textrm{repeated}\ M\ \textrm{times}}.
\end{multline*}
This is the desired $(\varepsilon,T_0)$-chain whose total duration is within
$\tau$ of $NT_0$\@.

Now, with this construction, we prove the desired implication.  We let
$\varepsilon\in\mappings[0]{\ts{X}}{\realp}$\@.  By
Lemma~\ref{lem:Papprox}\@, let $\delta_1\in\mappings[0]{\ts{X}}{\realp}$ be
such that
\begin{equation}\label{eq:chain-equiv2}
\d(z_1,z_2)<\delta_1(z_1)\implies\frac{1}{2}\varepsilon(z_1)<
\varepsilon(z_2)<\frac{3}{2}\varepsilon(z_1).
\end{equation}
By Lemma~\ref{lem:PPhiestimate}\@, let
$\delta_2\in\mappings[0]{\ts{X}}{\realp}$ be such that
\begin{multline}\label{eq:chain-equiv1}
\d(z_1,z_2)<\delta_2(z_1)\implies\d(\Phi(t,z_1),\Phi(t,z_2))\\
<\min\left\{\frac{1}{4}\varepsilon(\Phi(t,z_1)),
\frac{1}{2}\delta_1(\Phi(t,z_1))\right\},\qquad t\in\interval[0,{3T_0}].
\end{multline}
Let
\begin{equation}\label{eq:chain-equiv3}
u_0=x,u_1,\dots,u_l=y,\quad r_0,r_1,\dots,r_{l-1}
\end{equation}
and
\begin{equation}\label{eq:chain-equiv4}
v_0=y,v_1,\dots,v_m=x,\quad s_0,s_1,\dots,s_{m-1}
\end{equation}
be $(\min\{\frac{1}{32}\varepsilon,\frac{1}{4}\delta_2\},T_0)$-chains.  As we
saw above, we can suppose that $r_a\in\interval[{T_0},{2T_0})$\@,
$a\in\{0,1,\dots,l-1\}$\@, and $s_b\in\interval[{T_0},{2T_0})$\@,
$b\in\{0,1,\dots,m-1\}$\@.  Our construction above gives rise, for every
$\tau'\in\realp$\@, to $N\in\integerp$ and a
$(\min\{\frac{1}{8}\varepsilon,\delta_2\},T_0)$-chain
\begin{equation}\label{eq:chain-equiv5}
x_0=x,x_1,\dots,x_k=y,\quad t_0,t_1,\dots,t_{k-1}
\end{equation}
satisfying $\snorm{T-N}<\tau'$\@, where $T=\sum_{j=0}^{k-1}t_j$\@.
Furthermore, we can see from our constructions that $t_{k-1}=r_{l-1}$ and
$x_{k-1}=u_{l-1}$\@, independently of $\tau'$\@.  We now specify a suitable
value of $\tau'$\@.

To do so, denote
\begin{align*}
K_a=&\;\setdef{\Phi(t,u_a)}{t\in\interval[0,{r_a}]},\qquad
a\in\{0,1,\dots,l-1\},\\
L_b=&\;\setdef{\Phi(t,v_b)}{t\in\interval[0,{s_b}]},\qquad
a\in\{0,1,\dots,m-1\},\\
K=&\;\left(\bigcup_{a=0}^{l-1}K_a\right)\cup\left(\bigcup_{b=0}^{m-1}L_b\right).
\end{align*}
Note that $K_a$\@, $a\in\{0,1,\dots,l-1\}$\@, and $L_b$\@,
$b\in\{0,1,\dots,m-1\}$\@, are compact, and thus so too is $K$\@.
Importantly, note that $K$ depends only on the chains~\eqref{eq:chain-equiv3}
and~\eqref{eq:chain-equiv4}\@,~\ie~not on the choice of any $\tau'$\@.  Let
\begin{equation*}
\alpha=\inf\asetdef{\min\left\{\frac{1}{8}\varepsilon(z),
\frac{1}{2}\delta_1(z),\delta_2(z)\right\}}{z\in K},
\end{equation*}
noting that $\alpha>0$\@.  We claim that we can choose
$\tau\in\interval(0,{T_0})$ such that
\begin{equation}\label{eq:chain-equiv6}
\snorm{s}<\tau\implies\d(\Phi(s,x),x)<\alpha\qquad x\in K.
\end{equation}
Indeed, the function
\begin{equation*}
\interval[-1,1]\times K\ni(t,z)\mapsto\Phi(t,z)
\end{equation*}
is uniformly continuous (its domain being compact), and so there exists $\eta\in\realp$ such that
\begin{equation*}
(\snorm{\tau_1-\tau_2}<\eta,\ \d(z_1,z_2)<\eta)\implies
\d(\Phi(\tau_1,z_2),\Phi(\tau_2,z_2))<\alpha.
\end{equation*}
Applying this formula to $\tau_1=s$\@, $\tau_2=0$\@, and $z_1=z_2=z$ gives
the claim.  Thus we choose $\tau\in\interval(0,{T_0})$ so
that~\eqref{eq:chain-equiv6} holds and choose the
chain~\eqref{eq:chain-equiv5} so that $\snorm{T-NT_0}<\tau$\@, where
$T=\sum_{j=0}^{k-1}t_j$\@.  Note that we can, without loss of generality,
assume that $T>NT_0$\@, and so $T-NT_0\in\interval[0,\tau)$\@.

We now build points $y_0,y_1,\dots,y_N\in\ts{X}$ such that
\begin{compactenum}
\item $y_0=x$\@,
\item $y_N=y$\@, and
\item $\d(\Phi(T_0,y_j),y_{j+1})<\varepsilon(\Phi(T_0,y_j))$\@,
$j\in\{0,1,\dots,N-1\}$\@.
\end{compactenum}
We take $y_0=x_0=x$\@.  If $t_0=t_1=\dots=t_j=T_0$ for some
$j\in\{0,1,\dots,j-2\}$\@, then we can take $y_a=x_a$\@,
$a\in\{0,1,\dots,j\}$\@.  Our constructions then begin with the next jump.
We can, therefore, simplify life notationally, and not lose generality, by
simply supposing that $t_0\in\interval({T_0},{2T_0})$\@.  (We recommend that
a reader draw pictures of trajectories with times labelled to make sense of
the constructions that follow.)

We take $y_1=\Phi(T_0,y_0)$\@, whereupon
\begin{equation*}
\d(\Phi(T_0,y_0),y_1)=0<\varepsilon(\Phi(T_0,y_0)).
\end{equation*}

To define $y_2$\@, we first note that
\begin{equation*}
t_0\in\interval({T_0},{2T_0})\implies t_0-T_0\in\interval(0,{T_0}).
\end{equation*}
Thus we can follow the trajectory through $\Phi(T_0,y_0)$ for time $t_0-T_0$
and then jump to $x_1$ and follow the trajectory through $x_1$ for time
$T_0-(t_0-T_0)=2T_0-t_0$\@.  Since $2T_0-t_0\in\interval(0,{T_0})$\@, we have
$t_1>2T_0-t_0$\@, and so we do not need to jump from the trajectory through
$x_1$ and so we can define $y_2=\Phi(2T_0-t_0,x_1)$\@.

Since
\begin{equation*}
\d(\Phi(t_0,x_0),x_1)<\delta_2(\Phi(t_0,x_0))
\end{equation*}
and $2T_0-t_0\le2T_0$\@, by~\eqref{eq:chain-equiv1} we have
\begin{equation*}
\d(\Phi(2T_0-t_0,\Phi(t_0,x_0)),\Phi(2T_0-t_0,x_1))<
\frac{1}{4}\varepsilon(\Phi(2T_0-t_0,\Phi(t_0,x_0))).
\end{equation*}
Therefore,
\begin{align*}
\d(\Phi(T_0,y_1),y_2)=&\;\d(\Phi(2T_0,y_0),\Phi(2T_0-t_0,x_1)\\
=&\;\d(\Phi(2T_0-t_0,\Phi(t_0,x_0)),\Phi(2T_0-t_0,x_1))\\
<&\;\frac{1}{4}\varepsilon(\Phi(2T_0-t_0,\Phi(t_0,x_0)))\\
=&\;\frac{1}{4}\varepsilon(\Phi(2T_0,x_0))<\varepsilon(\Phi(T_0,y_0)).
\end{align*}

Note that $y_2$ lies on the trajectory through $x_1$ of time duration
$t_1$\@.  Specifically, $y_2=\Phi(s,x_1)$ with $s=2T_0-t_0$\@.  To define
$y_3$\@, we need to follow the chain for time $T_0$\@, and we need to
bookkeep if and when we need to jump to $x_2$\@.  There are three cases to
consider.
\begin{compactenum}
\item $t_1-s>T_0$\@: In this case we do not need to jump, and can
immediately define $y_3=\Phi(T_0,y_2)$\@.  In this case we have
\begin{equation*}
\d(\Phi(T_0,y_2),y_3)=0<\varepsilon(\Phi(T_0,y_2)).
\end{equation*}

\item $t_1-s=T_0$\@: In this case, we jump at the end of the time interval
of duration $T_0$\@.  That is, we take $y_3=x_2$\@.  Note that
\begin{equation*}
\Phi(T_0,y_2)=\Phi(t_1,x_1).
\end{equation*}
Thus we have
\begin{align*}
\d(\Phi(T_0,y_2),y_3)=&\;\d(\Phi(t_1,x_1),x_2)<
\frac{1}{4}\varepsilon(\Phi(t_1,x_1))\\
=&\;\frac{1}{4}\varepsilon(\Phi(T_0,y_2))<\varepsilon(\Phi(T_0,y_2)).
\end{align*}

\item $t_1-s<T_0$\@: In this case, we follow the trajectory through $y_2$ for
time $t_1-s$\@, jump to $x_2$\@, and then take $y_3=\Phi(T_0-(t_1-s),x_2)$\@.
Note that
\begin{equation*}
\Phi(T_0,y_2)=\Phi(T_0,\Phi(s,x_1))=\Phi(T_0+s,x_1).
\end{equation*}
Since
\begin{equation*}
\d(\Phi(t_1,x_1),x_2)<\delta_2(\Phi(t_1,x_1))
\end{equation*}
and since $T_0-(t_1-s)\le2T_0$\@, by~\eqref{eq:chain-equiv1} we have
\begin{equation*}
\d(\Phi(T_0-(t_1-s),\Phi(t_1,x_1)),\Phi(T_0-(t_1-s),x_2))<
\frac{1}{4}\varepsilon(\Phi(T_0-(t_1-s),\Phi(t_1,x_1))).
\end{equation*}
Therefore,
\begin{align*}
\d(\Phi(T_0,y_2),y_3)=&\;\d(\Phi(T_0+s,x_1),\Phi(T_0-(t_1-s),x_2))\\
=&\;\d(\Phi(T_0-(t_1-s),\Phi(t_1,x_1)),\Phi(T_0-(t_1-s),x_2))\\
<&\;\frac{1}{4}\varepsilon(\Phi(T_0-(t_1-s),\Phi(t_1,x_1)))\\
=&\;\frac{1}{4}\varepsilon(\Phi(T_0,y_2))<\varepsilon(\Phi(T_0,y_2)).
\end{align*}
\end{compactenum}

Now one can proceed as in the construction of $y_3$ to define
$y_4,\dots,y_{N-1}$\@.  Note that, having done this, we have followed the
chain~\eqref{eq:chain-equiv5} for time duration $(N-1)T_0$\@.  Thus the time
remaining along the chain is
$T-(N-1)T_0=T_0+(T-NT_0)\in\interval[{T_0},{T_0+\tau})$\@.  That is to say,
the time remaining in the chain is within $\tau$ of $T_0$\@.

We take $y_N=y$ and consider two cases.
\begin{compactenum}
\item $y_{N-1}=\Phi(s,x_{k-1})$ for some $s\in\interval[0,{t_{k-1}})$\@: In
this case, the time left to travel along the chain~\eqref{eq:chain-equiv5} is
$t_{k-1}-s=T_0+(T-NT_0)$\@.  Thus
\begin{equation*}
\Phi(T_0,y_{N-1})=\Phi(T_0+s,x_{k-1})=\Phi(t_{k-1}-(T-NT_0),x_{k-1}).
\end{equation*}
Since $\snorm{NT_0-T}<\tau$\@, by~\eqref{eq:chain-equiv6}
and~\eqref{eq:chain-equiv2} we have
\begin{align*}
&\Phi(T_0,y_{N-1})=\Phi(NT_0-T,\Phi(t_{k-1},x_{k-1}))\\
\implies&\d(\Phi(T_0,y_{N-1}),\Phi(t_{k-1},x_{k-1}))<\alpha\le
\delta_1(\Phi(t_{k-1},x_{k-1}))\\
\implies&\frac{1}{2}\varepsilon(\Phi(t_{k-1},x_{k-1}))<
\varepsilon(\Phi(T_0,y_{N-1}))<
\frac{3}{2}\varepsilon(\Phi(t_{k-1},x_{k-1})).
\end{align*}
Similarly, since $\snorm{T-NT_0}<\tau$ and $t_{k-1}\le2$\@,
by~\eqref{eq:chain-equiv6} and~\eqref{eq:chain-equiv1} we have
\begin{align*}
&\d(x_{k-1},\Phi(NT_0-T,x_{k-1}))<\alpha\le\delta_2(x_{k-1})\\
\implies&\d(\Phi(t_{k-1},x_{k-1}),\Phi(t_{k-1},\Phi(NT_0-T,x_{k-1})))<
\frac{1}{4}\varepsilon(\Phi(t_{k-1},x_{k-1})).
\end{align*}
Therefore,
\begin{align*}
\d(\Phi(T_0,y_{N-1}),y_N)=&\;\d(\Phi(t_{k-1}-(T-NT_0),x_{k-1}),y)\\
\le&\;\d(\Phi(t_{k-1}-(T-NT_0),x_{k-1}),\Phi(t_{k-1},x_{k-1})\\
&\;+\d(\Phi(t_{k-1},x_{k-1}),y)\\
=&\;\d(\Phi(t_{k-1},\Phi(NT_0-T,x_{k-1})),\Phi(t_{k-1},x_{k-1}))\\
&\;+\d(\Phi(t_{k-1},x_{k-1}),y)\\
<&\;\frac{1}{4}\varepsilon(\Phi(t_{k-1},x_{k-1}))+
\frac{1}{4}\varepsilon(\Phi(t_{k-1},x_{k-1}))\\
=&\;\frac{1}{2}\varepsilon(\Phi(t_{k-1},x_{k-1}))<
\varepsilon(\Phi(T_0,y_{N-1})).
\end{align*}

\item $y_{N-1}=\Phi(s,x_{k-2})$ for some $s\in\interval[0,{t_{k-2}})$\@:
In this case, to get to $y_N$\@, we must jump to the trajectory through
$x_{k-1}$ after time $t_{k-2}-s$\@, and then follow this trajectory for time
$t_{k-1}-(T-NT_0)$\@.  Note that
\begin{equation*}
t_{k-1}-(T-NT_0)+t_{k-2}-s=T_0\implies T_0+s=t_{k-1}+t_{k-2}-(T-NT_0).
\end{equation*}
Also note that
\begin{equation*}
(t_{k-1}\in\interval[{T_0},{2T_0}),\ T-NT_0\in\interval[0,\tau))
\implies t_{k-1}-(T-NT_0)\le2T_0.
\end{equation*}
Since
\begin{equation*}
\d(\Phi(t_{k-2},x_{k-2}),x_{k-1})<\delta_2(\Phi(t_{k-2},x_{k-2})),
\end{equation*}
and $t_{k-1}-(T-NT_0)\le2T_0$\@, by~\eqref{eq:chain-equiv1} we have
\begin{align}\notag
\frac{1}{4}\varepsilon(\Phi(T_0+&s,x_{k-2}))
=\frac{1}{4}\varepsilon(\Phi(t_{k-1}+t_{k-2}-(T-NT_0),x_{k-2}))\\\notag
>&\;\d(\Phi(t_{k-1}+t_{k-2}-(T-NT_0),x_{k-2}),\Phi(t_{k-1}-(T-NT_0),x_{k-1}))\\
\label{eq:chain-equiv7}
=&\;\d(\Phi(T_0+s,x_{k-2}),\Phi(t_{k-1}-(T-NT_0),x_{k-1})).
\end{align}
Since $\snorm{NT_0-T}<\tau$\@, by~\eqref{eq:chain-equiv6} we have
\begin{align}\notag
&\Phi(t_{k-1}-(T-NT_0),x_{k-1})=\Phi(NT_0-T,\Phi(t_{k-1},x_{k-1}))\\
\label{eq:chain-equiv8}
\implies&\d(\Phi(t_{k-1}-(NT_0-T),x_{k-1})),\Phi(t_{k-1},x_{k-1}))<
\alpha\le\frac{1}{8}\varepsilon(\Phi(t_{k-1},x_{k-1})).
\end{align}
Using the definition of $\alpha$\@, the same argument gives
\begin{equation*}
\d(\Phi(t_{k-1}-(T-NT_0),x_{k-1}),\Phi(t_{k-1},x_{k-1})<
\alpha\le\frac{1}{2}\delta_1(\Phi(t_{k-1},x_{k-1})).
\end{equation*}
Since
\begin{equation*}
\d(\Phi(t_{k-2},x_{k-2}),x_{k-1})<\delta_2(\Phi(t_{k-2},x_{k-2}))
\end{equation*}
and $t_{k-1}-(T-NT_0)\le2T_0$\@, by~\eqref{eq:chain-equiv1} we have
\begin{multline*}
\d(\Phi(t_{k-1}-(T-NT_0),\Phi(t_{k-2},x_{k-2})),
\Phi(t_{k-1}-(T-NT_0),x_{k-1}))\\
<\frac{1}{2}\delta_1(\Phi(t_{k-1}-(T-NT_0),\Phi(t_{k-2},x_{k-2})))=
\frac{1}{2}\delta_1(T_0+s,x_{k-2}).
\end{multline*}
Putting the above together, we have
\begin{align*}
\d(\Phi(t_{k-1}+&t_{k-2}-(T-NT_0),x_{k-2}),\Phi(t_{k-1},x_{k-1}))\\
=&\;\d(\Phi(t_{k-1}-(T-NT_0),\Phi(t_{k-2},x_{k-2})),
\Phi(t_{k-1},x_{k-1})\\
\le&\;\d(\Phi(t_{k-1}-(T-NT_0),x_{k-1}),\Phi(t_{k-1},x_{k-1})\\
&\;+\d(\Phi(t_{k-1}-(T-NT_0),\Phi(t_{k-2},x_{k-2})),
\Phi(t_{k-1}-(T-NT_0),x_{k-1}))\\
<&\;\frac{1}{2}\delta_1(\Phi(t_{k-1},x_{k-1}))+
\frac{1}{2}\delta_1(T_0+s,x_{k-2})\\
\le&\;\max\{\delta_1(\Phi(t_{k-1},x_{k-1})),\delta_1(\Phi(T_0+s,x_{k-2}))\},
\end{align*}
using the standard relation $\dnorm{\cdot}_1\le n\dnorm{\cdot}_\infty$
between the $1$- and $\infty$-norms for $\real^n$\@.  Keeping in mind that
$T_0+s=t_{k-1}+t_{k-2}-(T-NT_0)$\@, we now consider two cases.
\begin{compactenum}[(a)]
\item
$\delta_1(\Phi(T_0+s,x_{k-2}))\le\delta_1(\Phi(t_{k-1},x_{k-1}))$\@: In
this case,
\begin{equation*}
\max\{\delta_1(\Phi(t_{k-1},x_{k-1})),\delta_1(\Phi(T_0+s,x_{k-2}))\}=
\delta_1(\Phi(t_{k-1},x_{k-1})).
\end{equation*}
By~\eqref{eq:chain-equiv2} we have
\begin{align*}
&\d(\Phi(T_0+s,x_{k-2}),\Phi(t_{k-1},x_{k-1}))<
\delta_1(\Phi(t_{k-1},x_{k-1}))\\
\implies&\frac{1}{2}\varepsilon(\Phi(t_{k-1},x_{k-1}))<
\varepsilon(\Phi(T_0+s,x_{k-2}))<
\frac{3}{2}\varepsilon(\Phi(t_{k-1},x_{k-1}))
\end{align*}

\item $\delta_1(T_0+s,y_{N-1})>\delta_1(\Phi(t_{k-1},x_{k-1}))$\@: In this
case,
\begin{equation*}
\max\{\delta_1(\Phi(t_{k-1},x_{k-1})),\delta_1(\Phi(T_0+s,y_{N-1}))\}=
\delta_1(\Phi(T_0+s,x_{k-2}))
\end{equation*}
and, again by~\eqref{eq:chain-equiv2}\@, we have
\begin{align*}
&\d(\Phi(T_0+s,x_{k-2}),\Phi(t_{k-1},x_{k-1}))<
\delta_1(\Phi(T_0+s,x_{k-s}))\\
\implies&\frac{1}{2}\varepsilon(\Phi(T_0+s,x_{k-2}))<
\varepsilon(\Phi(t_{k-1},x_{k-1}))<
\frac{3}{2}\varepsilon(\Phi(T_0+s,x_{k-2})).
\end{align*}
\end{compactenum}
Therefore, considering both cases together,
\begin{equation}\label{eq:chain-equiv9}
\varepsilon(\Phi(t_{k-1},x_{k-1}))<2\varepsilon(\Phi(T_0+s,x_{k-2})).
\end{equation}
Finally, we assemble~\eqref{eq:chain-equiv7}\@,~\eqref{eq:chain-equiv8}\@,
and~\eqref{eq:chain-equiv9}\@:
\begin{align*}
\d(\Phi(T_0,y_{N-1}),y_N)=&\;\d(\Phi(T_0+s,x_{k-2}),y)\\
\le&\;\d(\Phi(T_0+s,x_{k-2}),\Phi(t_{k-1}-(T-NT_0),x_{k-1}))\\
&\;+\d(\Phi(t_{k-1}-(T-NT_0),x_{k-1}),\Phi(t_{k-1},x_{k-1})))\\
&\;+\d(\Phi(t_{k-1},x_{k-1}),y)\\
<&\;\frac{1}{4}\varepsilon(\Phi(T_0+s,x_{k-2}))
+\frac{1}{8}\varepsilon(\Phi(t_{k-1},x_{k-1}))\\
&\;+\frac{1}{8}\varepsilon(\Phi(t_{k-1},x_{k-1}))\\
<&\;\frac{3}{4}\varepsilon(\Phi(T_0+s,x_{k-2}))<
\varepsilon(\Phi(T_0,y_{N-1})).
\end{align*}
\end{compactenum}
Thus the points $y_0,y_1,\dots,y_N$ have the desired properties.

\eqref{pl:chain-equiv3}$\implies$\eqref{pl:chain-equiv1} Under the stated
hypotheses, let $\varepsilon\in\mappings[0]{\ts{X}}{\realp}$ and let
$T\in\tdomainp$\@.  Without loss of generality, suppose that $T=MT_0$ for
some $M\in\integerp$\@.  By Lemma~\ref{lem:Papprox}\@, let
$\varepsilon_{2M-1}\in\mappings[0]{\ts{X}}{\realp}$ be such that
$\varepsilon_{2M-1}<\frac{1}{2}\varepsilon$ and such that
\begin{equation}\label{eq:chain-equiv10}
\d(z_1,z_2)<\varepsilon_{2M-1}(z_1)\implies
\varepsilon(z_2)<\tfrac{3}{2}\varepsilon(z_1).
\end{equation}
Then, by Lemma~\ref{lem:general-epsilon-delta}\@, let
$\eta_{2M-2}\in\mappings[0]{\ts{X}}{\realp}$ be such that
\begin{equation}\label{eq:chain-equiv11}
\d(z_1,z_2)<\eta_{2M-2}(z_1)\implies\d(\Phi(T_0,z_1),\Phi(T_0,z_2))<
\varepsilon_{2M-1}(\Phi(T_0,z_1)).
\end{equation}
Then, using Lemmata~\ref{lem:Papprox} and~\ref{lem:general-epsilon-delta}\@,
recursively define
$\varepsilon_{2M-1},\dots,\varepsilon_2\in\mappings[0]{\ts{X}}{\realp}$ and
$\eta_{2M-2},\dots,\eta_1\in\mappings[0]{\ts{X}}{\realp}$ such that
$\varepsilon_j<\frac{1}{2}\eta_j$\@, $j\in\{2,\dots,2M-2\}$\@,
$x\in\ts{X}$\@, and such that
\begin{equation}\label{eq:chain-equiv12}
\d(z_1,z_2)<\varepsilon_j(z_1)\implies\eta_j(z_2)<\tfrac{3}{2}\eta_j(z_1),
\qquad j\in\{2,\dots,2M-1\},
\end{equation}
and
\begin{equation}\label{eq:chain-equiv13}
\d(z_1,z_2)<\eta_j(z_1)\implies\d(\Phi(T_0,z_1),\Phi(T_0,z_2))<
\varepsilon_{j+1}(\Phi(T_0,z_1)),\qquad j\in\{1,\dots,2M-2\}.
\end{equation}
Define $\delta\in\mappings[0]{\ts{X}}{\realp}$ by
\begin{equation*}
\delta=\min\{\eta_1,\tfrac{1}{3}\eta_2,\dots,
\tfrac{1}{3}\eta_{N-1},\tfrac{1}{3}\varepsilon\}.
\end{equation*}

Now, by hypothesis, let
\begin{equation*}
x_0=x,x_1,\dots,x_k=y,\quad T_0,T_0,\dots,T_0
\end{equation*}
and
\begin{equation*}
y_0=y,y_1,\dots,y_m=x,\quad T_0,T_0,\dots,T_0
\end{equation*}
be $\varepsilon$-$T_0$-chains,~\ie
\begin{equation}\label{eq:chain-equiv14}
\d(\Phi(T_0,x_j),x_{j+1})<\delta(\Phi(T_0,x_j)),\qquad j\in\{1,\dots,k\},
\end{equation}
and
\begin{equation}\label{eq:chain-equiv15}
\d(\Phi(T_0,y_l),y_{l+1})<\delta(\Phi(T_0,y_l)),\qquad l\in\{1,\dots,m\}.
\end{equation}
From the two finite sequences $x_0,x_1,\dots,x_k$ and $y_0,y_1,\dots,y_m$\@,
build another finite sequence
\begin{multline*}
x_0=x,x_1,\dots,x_{k-1},\\
\underbrace{y_0,y_1,\dots,y_{m-1},x_0,x_1,\dots,x_{k-1},
\dots,y_0,y_1,\dots,y_{m-1},x_0,x_1,\dots,x_{k-1}}_{y_0,y_1,\dots,y_{m-1},
x_0,x_1,\dots,x_{k-1}\ \textrm{repeated}\ M\ \textrm{times}},y,
\end{multline*}
this of length $k+M(k+m)+1$\@.  Let us use this sequence to assemble an
$\varepsilon$-$T_0$-chain
\begin{equation*}
z_0,z_1,\dots,z_{k+M(k+m)},\quad T_0,T_0,\dots,T_0,
\end{equation*}
\ie
\begin{equation*}
\d(\Phi(T_0,z_j),z_{j+1})<\delta(\Phi(T_0,z_j)),\qquad
j\in\{0,1,\dots,k+M(k+m)-1\}.
\end{equation*}
By the Euclidean Algorithm, let $N',R'\in\integerp$ be such that
\begin{equation*}
k+M(k+m)=N'M+R',\quad R'\in\{0,1,\dots,M-1\}.
\end{equation*}
and take $N=N'-1$ and $R=R'+M$ so that
\begin{multline*}
k+m(k+m)=N'M+R'=(N+1)M+R-M=NM+R,\\R\in\{M,M+1,\dots,2M-1\}.
\end{multline*}
Define
\begin{equation*}
z'_j=\begin{cases}z_{jM},&j\in\{0,1,\dots,N-1\},\\y,&j=N\end{cases}
\end{equation*}
and
\begin{equation*}
t_j=\begin{cases}T,&j\in\{0,1,\dots,N-2\},\\R,&j=N-1.\end{cases}
\end{equation*}
We claim that
\begin{equation*}
z'_0,z'_1,\dots,z'_N,\quad t_0,t_1,\dots,t_{N-1}
\end{equation*}
is an $(\varepsilon,T)$-chain from $x$ to $y$\@.

To establish this, we first claim that, for $j\in\{0,1,\dots,N\}$ and for
$l\in\{1,\dots,R\}$\@, we have
\begin{equation}\label{eq:chain-equiv16}
\d(\Phi(lT_0,z_{jM}),z_{jM+l})<\eta_l(\Phi(lT_0,z_{jM})).
\end{equation}
For $l=1$\@, we have
\begin{equation*}
\d(\Phi(T_0,z_{jM}),z_{jM+1})<\delta(\Phi(T_0,z_{jM}))<
\eta_1(\Phi(T_0,z_{jM})),
\end{equation*}
verifying the claim in this case.  To argue inductively, suppose that
\begin{equation*}
\d(\Phi(T_0,z_{jM}),z_{jM+l})<\eta_l(\Phi(lT_0,z_{jM}))
\end{equation*}
for some $l\in\{2,\dots,R-1\}$\@.  From this inequality and
by~\eqref{eq:chain-equiv13}\@, we have
\begin{equation*}
\d(\Phi((l+1)T_0,z_{jM}),\Phi(T_0,z_{jM+l}))<
\varepsilon_{l+1}(\Phi((l+1)T_0,z_{jM}))
<\tfrac{1}{2}\eta_{l+1}(\Phi((l+1)T_0,z_{jM})).
\end{equation*}
Similarly, by~\eqref{eq:chain-equiv12} we have
\begin{equation*}
\eta_{l+1}(\Phi(T_0,z_{jM+l}))<
\tfrac{3}{2}\eta_{l+1}(\Phi((l+1)T_0,z_{jM})).
\end{equation*}
Thus, using~\eqref{eq:chain-equiv14}\@,~\eqref{eq:chain-equiv15}\@, and the
definition of $\delta$\@,
\begin{align*}
\d(\Phi(T_0,z_{jM+l}),z_{jM+l+1})<&\;\delta(\Phi(T_0,z_{jM+l}))
<\tfrac{1}{3}\eta_{l+1}(\Phi(T_0,z_{jM+l}))\\
<&\;\tfrac{1}{2}\eta_{l+1}(\Phi((l+1)T_0,z_{jM})).
\end{align*}
Putting this together,
\begin{align*}
\d(\Phi((l+1)T_0,z_{jM}),z_{jM+l+1})
\le&\;\d(\Phi((l+1)T_0,z_{jM}),\Phi(T_0,z_{jM+l}))\\
&\;+\d(\Phi(T_0,z_{jM}),z_{jM+l+1})\\
<&\;\eta_{l+1}(\Phi((l+1)T_0,z_{jM})).
\end{align*}
This establishes~\eqref{eq:chain-equiv16} by induction.

In the particular case of $l=N$\@,~\eqref{eq:chain-equiv16} gives
\begin{equation*}
\d(\Phi(NT_0,z'_j),z'_{j+l})<\eta_N(\Phi(NT_0,z'_j))<
\varepsilon(\Phi(NT_0,z'_j)).
\end{equation*}
If $l=R$ and $j=N-1$\@,~\eqref{eq:chain-equiv16} gives
\begin{equation*}
\d(\Phi(RT_0,z'_{N-1}),z'_N)<\eta_{N-1}(\Phi(RT_0,z'_{N-1}))<
\varepsilon(\Phi(RT_0,z'_{N-1})).
\end{equation*}

Of course, the same argument gives, for
$\varepsilon\in\mappings[0]{\ts{X}}{\realp}$ and $T\in\tdomainp$\@, an
$(\varepsilon,T)$-chain from $y$ to $x$\@, which gives this part of the
result.
\end{proof}
\end{theorem}

The following corollary is immediate since a point is chain recurrent if and
only if it is chain equivalent to itself.
\begin{corollary}\label{cor:chan-recurrent}
Let\/ $(\ts{X},\d)$ be a metric space and let\/ $\Phi$ be a topological
flow or semiflow on\/ $\ts{X}$\@.  Then, for\/ $T_0\in\tdomainp$ possessing a
multiplicative inverse in\/ $\tdomain$\@,
\begin{equation*}
\ChRec(\Phi)=\ChRec_{\ge{}T_0}(\Phi)=\ChRec_{=T_0}(\Phi)=\ChRec(\Phidisc{T_0}).
\end{equation*}
\end{corollary}

Note that the theorem and the corollary indicate why, when working with
discrete-time flows and semiflows, one can ignore the switching times and
simply take them to be $1$\@.  This leads to the following simplified
notion of a chain in these cases.  Indeed, this is the usual definition of a
chain for discrete-time flows and semiflows.
\begin{definition}\label{def:vareps-chain}
Let $(\ts{X},\d)$ be a metric space and let $\Phi$ be a discrete-time
topological flow or semiflow on $\ts{X}$ with\/ $\phi=\Phi_1$\@.  For
$x,y\in\ts{X}$ and $\varepsilon\in\mappings[0]{\ts{X}}{\realp}$\@, an
\defn{$\varepsilon$-chain} for $\Phi$ from $x$ to $y$ is a finite sequence
\begin{equation*}
x_0,x_1,\dots,x_k,
\end{equation*}
with
\begin{compactenum}[(i)]
\item $x_0,x_1,\dots,x_k\in\ts{X}$\@,
\item $x_0=x$ and $x_k=y$\@, and
\item $\d(\phi(x_j),x_{j+1})<\varepsilon(\phi(x_j))$\@,
$j\in\{0,1,\dots,k-1\}$\@.\oprocend
\end{compactenum}
\end{definition}

Staying with discrete-time flows and semiflows for a moment, let us introduce
some notation for these that will facilitate a comparison of the chain
recurrent set of a mapping with the chain recurrent set of its iterate
mappings.  Let $\Phi$ be a discrete-time flow or semiflow on a topological
space $(\ts{X},\sO)$ and let $k\in\integerp$\@.  Then define the
discrete-time flow or semiflow $\Phi^k$ on $\ts{X}$ by
\begin{equation*}
\mapdef{\Phi^k}{\tdomain\times\ts{X}}{\ts{X}}{(j,x)}{\Phi(jk,x).}
\end{equation*}
With this notation, we have the following result, which is essentially
corollary to Theorem~\ref{the:chain-equiv}\@.
\begin{corollary}\label{cor:Phi=Phik}
Let\/ $(\ts{X},\d)$ be a metric space and let\/ $\Phi$ be a topological flow
or semiflow on\/ $\ts{X}$\@.  Then\/ $\ChRec(\Phi)=\ChRec(\Phi^k)$ for
every\/ $k\in\integerp$\@.
\begin{proof}
If $x\in\ChRec(\Phi^k)$\@, then an elementary argument like that given in the
discrete-time case of the
implication~\eqref{pl:chain-equiv2}$\implies$\eqref{pl:chain-equiv3} from the
proof of Theorem~\ref{the:chain-equiv} shows that $x\in\ChRec(\Phi)$\@.

Now suppose that $x\in\ChRec(\Phi)$\@.  Let
$\varepsilon\in\mappings[0]{\ts{X}}{\realp}$ and, using an inductive proof
like that from the proof of the
implication~\eqref{pl:chain-equiv3}$\implies$\eqref{pl:chain-equiv1} from
Theorem~\ref{the:chain-equiv}\@, let $\delta\in\mappings[0]{\ts{X}}{\realp}$
be such that, if
\begin{equation*}
y_0,y_1,\dots,y_k
\end{equation*}
is a $\delta$-chain for $\Phi$\@, then
\begin{equation*}
\d(y_k,\Phi^k(y_0))<\varepsilon(\Phi^k(y_0)).
\end{equation*}
By Theorem~\ref{the:chain-equiv}\@, suppose that
\begin{equation*}
x_0,x_1,\dots,x_m
\end{equation*}
is a $\delta$-chain from $x$ to $x$ for $\Phi$\@.  Then
\begin{equation*}
\underbrace{x_0,x_1,\dots,x_{m-1},x_0,x_1,\dots,x_m,\dots,
x_0,x_1,\dots,x_m}_{k\ \textrm{times}}
\end{equation*}
is a $\delta$-chain from $x$ to $x$ for $\Phi$\@, and the choice of $\delta$
ensures that
\begin{equation*}
x_0,x_k,\dots,x_{mk}
\end{equation*}
is a $\varepsilon$-chain from $x$ to $x$ for $\Phi^k$\@.  Thus, by
Theorem~\ref{the:chain-equiv}\@, $x\in\ChRec(\Phi^k)$\@.
\end{proof}
\end{corollary}

\section{The Conley decomposition}

Now we turn to the first part of the Fundamental Theorem of Dynamical
Systems, the so-called Conley decomposition.  This gives a decomposition of
the state space for a flow or semiflow into chain recurrent dynamics on the
chain recurrent set and gradient-like dynamics off the chain recurrent set,
as (roughly) one flows from repelling sets to attracting sets.  We first
spend some time understanding the relationship between chains and trapping
regions, as this is essential for any sort of understanding of how the
Fundamental Theorem of Dynamical Systems works.  After this, we prove the
decomposition theorem.

\subsection{Chains and trapping regions}\label{subsec:chain-trap}

We begin with some notation.  Let $(\ts{X},\d)$ be a metric space and let
$S\subset\ts{X}$\@.  The function
\begin{equation*}
\mapdef{\dist_S}{\ts{X}}{\realnn}{x}
{\inf\setdef{\d(x,y)}{y\in S}}
\end{equation*}
is the \defn{distance function to\/ $S$}\@.  We shall also denote
$\dist(x,S)=\dist_S(x)$\@.  Evidently, if $S$ is closed and $x\not\in S$\@,
$\dist_S(x)\in\realp$\@.  We use this notation in the proof of the following
lemma that shows how trapping regions give natural error functions with
useful properties.
\begin{lemma}\label{lem:trap-erf}
Let\/ $(\ts{X},\d)$ be a metric space and let\/ $\Phi$ be a topological flow
or semiflow on\/ $\ts{X}$\@.  For an open set\/ $\nbhd{T}\subset\ts{X}$ and
for\/ $T\in\tdomainp$\@, the following statements are equivalent:
\begin{compactenum}[(i)]
\item \label{pl:trap-erf1} $\nbhd{T}$ is a trapping region with\/ $\closure(\Phi(\tdomain_{\ge{}T}\times\nbhd{T}))\subset\nbhd{T}$\@;
\item \label{pl:trap-erf2} there exists\/
$\varepsilon\in\mappings[0]{\ts{X}}{\interval(0,1]}$ such that
\begin{equation*}
\oball[\d{}]{\varepsilon(\Phi(t,x))}{\Phi(t,x)}\subset\nbhd{T},
\qquad t\in\tdomain_{\ge{}T},\ x\in\nbhd{T}.
\end{equation*}
\end{compactenum}
\begin{proof}
\eqref{pl:trap-erf1}$\implies$\eqref{pl:trap-erf2} Let
\begin{equation*}
\varepsilon'(x)=
\frac{1}{2}\left(\dist_{\closure(\Phi(\tdomain_{\ge{}T}\times\nbhd{T}))}(x)+
\dist_{\ts{X}\setminus\nbhd{T}}(x)\right),
\end{equation*}
noting that $\varepsilon'$ is continuous~\cite[Theorem~3.16]{CDA/KCB:06}\@.
Also, $\varepsilon'$ takes values in $\realp$ since, if
$x\in\closure(\Phi(\tdomain_{\ge{}T}\times\nbhd{T}))$\@, then
$x\not\in\ts{X}\setminus\nbhd{T}$\@.  For $x\in\nbhd{T}$ and
$t\in\tdomain_{\ge{}T}$\@, we have
\begin{equation*}
\Phi(t,x)\in\Phi(\tdomain_{\ge{}T}\times\nbhd{T})\implies
\dist_{\closure(\Phi(\tdomain_{\ge{}T}\times\nbhd{T}))}(\Phi(t,x))=0.
\end{equation*}
Thus, if $y\in\oball[\d{}]{\varepsilon'(\Phi(t,x))}{\Phi(t,x)}$\@, we have
\begin{equation*}
\d(\Phi(t,x),y)<\varepsilon'(\Phi(t,x))=
\frac{1}{2}\dist_{\ts{X}\setminus\nbhd{T}}(\Phi(t,x)).
\end{equation*}
Also,
\begin{align*}
\dist_{\ts{X}\setminus\nbhd{T}}(\Phi(t,x))=&\;
\inf\setdef{\d(\Phi(t,x),z)}{z\in\ts{X}\setminus\nbhd{T}}\\
\le&\;\inf\setdef{\d(\Phi(t,x),y)+\d(z,y)}{x\in\ts{X}\setminus\nbhd{T}}\\
=&\;\d(\Phi(t,x),y)+\dist_{\ts{X}\setminus\nbhd{T}}(y).
\end{align*}
Putting this together,
\begin{equation*}
2\d(\Phi(t,x),y)<\dist_{\ts{X}\setminus\nbhd{T}}(\Phi(t,x))\le
\d(\Phi(t,x),y)+\dist_{\ts{X}\setminus\nbhd{T}}(y).
\end{equation*}
Then
\begin{equation*}
\dist_{\ts{X}\setminus\nbhd{T}}(y)>\d(y,\Phi(t,x))\ge0
\implies y\in\nbhd{T},
\end{equation*}
showing that
$\oball[\d{}]{\varepsilon'(\Phi(t,x))}{\Phi(t,x)}\subset\nbhd{T}$\@.  Taking
\begin{equation*}
\varepsilon(x)=\min\{\varepsilon'(x),1\}
\end{equation*}
gives this part of the lemma.

\eqref{pl:trap-erf2}$\implies$\eqref{pl:trap-erf1} Let $\varepsilon$ be as
stated.  Let $y\in\closure(\Phi(\tdomain_{\ge{}T}\times\nbhd{T}))$ and let
$\ifam{(t_j,x_j)}_{j\in\integerp}$ be a sequence in
$\tdomain_{\ge{}T}\times\nbhd{T}$ for which
$y=\lim_{j\to\infty}\Phi(t_j,x_j)$\@.  By continuity of $\varepsilon$\@, let
$N\in\integerp$ be sufficiently large that
$\varepsilon(\Phi(t_j,x_j))\ge\frac{1}{2}\varepsilon(y)$\@, $j\ge N$\@.  Also
suppose that $N$ is large enough that
$\d(\Phi(t_j,x_j),y)<\frac{1}{4}\varepsilon(y)$ for $j\ge N$\@.  Let
$z\in\oball[\d{}]{\frac{1}{4}\varepsilon(y)}{y}$\@.  Then
\begin{align*}
\d(\Phi(t_N,x_N),z)\le&\;\d(\Phi(t_N,x_N),y)+\d(y,z)\\
<&\;\frac{1}{4}\varepsilon(y)+\frac{1}{4}\varepsilon(y)
\le\varepsilon(\Phi(t_N,x_N)).
\end{align*}
Thus
\begin{equation*}
\oball[\d{}]{\tfrac{1}{4}\varepsilon(y)}{y}\subset
\oball[\d{}]{\varepsilon(\Phi(t_N,x_N))}{\Phi(t_N,x_N)}\subset\nbhd{T}
\end{equation*}
and so $y\in\nbhd{T}$\@.
\end{proof}
\end{lemma}

The following lemma provides an essential conceptual step in understanding
the Fundamental Theorem of Dynamical Systems.  It shows how chains and the
chain recurrent set give rise to natural trapping regions.
\begin{lemma}\label{lem:chain->trap}
Let\/ $(\ts{X},\d)$ be a metric space and let\/ $\Phi$ be a topological flow
or semiflow on\/ $\ts{X}$\@.  Let\/ $x\in\ts{X}\setminus\ChRec(\Phi)$\@, and
let\/ $\varepsilon\in\mappings[0]{\ts{X}}{\realp}$ and\/ $T\in\tdomainp$ be
such that there is no\/ $(\varepsilon,T)$-chain from\/ $x$ to itself.  Let\/
$m\in\integerp$\@.  Then
\begin{equation*}
\nbhd{T}=\setdef{y\in\ts{X}}{\textrm{there exists an $(\varepsilon,T)$-chain
of length\/ $k\ge m$ from\/ $x$ to\/ $y$}}
\end{equation*}
is an open trapping region satisfying
\begin{compactenum}[(i)]
\item \label{pl:!ChRec-trap1}
$\Phi(\tdomain_{\ge{}T}\times\nbhd{T})\subset\nbhd{T}$ and
\item \label{pl:!ChRec-trap2} $x\not\in\nbhd{T}$\@.
\end{compactenum}
\begin{proof}
Clearly, $x\not\in\nbhd{T}$ by definition of $\varepsilon$ and $T$\@.

To prove openness of $\nbhd{T}$\@, let $y\in\nbhd{T}$\@, let $k\ge m$ and let
\begin{equation*}
x_0=x,x_1,\dots,x_k=y,\quad t_0,t_1,\dots,t_{k-1}
\end{equation*}
be an $(\varepsilon,T)$-chain from $x$ to $y\in\nbhd{T}$\@.  We claim that, if
\begin{equation*}
r=\varepsilon(\Phi(t_{k-1},x_{k-1}))-\d(\Phi(t_{k-1},x_{k-1}),y),
\end{equation*}
then $\oball[\d{}]{r}{y}\subset\nbhd{T}$\@.  Indeed, let
$z\in\oball[\d{}]{r}{y}$\@.  Then
\begin{align*}
\d(\Phi(t_{k-1},x_{k-1}),z)\le&\;
\d(\Phi(t_{k-1},x_{k-1}),y)+\d(y,z)\\
<&\;\d(\Phi(t_{k-1},x_{k-1}),y)+
\varepsilon(\Phi(t_{k-1},x_{k-1}))-\d(\Phi(t_{k-1},x_{k-1}),y)\\
=&\;\varepsilon(\Phi(t_{k-1},x_{k-1})).
\end{align*}
Thus we have
\begin{equation*}
\oball[\d{}]{r}{y}\subset
\oball[\d{}]{\varepsilon(\Phi(t_{k-1},x_{k-1}))}{\Phi(t_{k-1},x_{k-1})}.
\end{equation*}
Since the ball on the right consists of points $z$ for which there is an
$(\varepsilon,T)$-chain from $x$ to $z$\@, we obtain
$\oball[\d{}]{r}{y}\subset\nbhd{T}$\@, giving the desired openness.

Now we show that
$\closure(\Phi(\tdomain_{\ge{}T}\times\nbhd{T}))\subset\nbhd{T}$\@.  By
Lemma~\ref{lem:Papprox}\@, let $\delta\in\mappings[0]{\ts{X}}{\realp}$ be
such that
\begin{equation*}
\d(z_1,z_2)<\delta(z_1)\implies\frac{1}{2}\varepsilon(z_1)<
\varepsilon(z_2).
\end{equation*}
Without loss of generality, we can assume that
$\delta\le\frac{1}{2}\varepsilon$\@.  Let
$y\in\closure(\Phi(\tdomain_{\ge{}T}\times\nbhd{T}))$\@.  Then
\begin{equation*}
\oball[\d{}]{\delta(y)}{y}\cap\Phi(\tdomain_{\ge{}T}\times\nbhd{T})
\not=\emptyset,
\end{equation*}
and so there exists $t\in\tdomain_{\ge{}T}$ and $z\in\nbhd{T}$ such that
$\Phi(t,z)\in\oball[\d{}]{\delta(y)}{y}$\@.  By definition of $\nbhd{T}$\@,
let $k\ge m$ and let
\begin{equation*}
x_0=x,x_1,\dots,x_k=z,\quad t_0,t_1,\dots,t_{k-1}
\end{equation*}
be an $(\epsilon,T)$-chain from $x$ to $z$\@.  We have
\begin{equation*}
\d(\Phi(t,z),y)<\delta(y)\le\frac{1}{2}\varepsilon(y)<
\varepsilon(\Phi(t,z)),
\end{equation*}
from which we may conclude that
\begin{equation*}
x_0=x,x_1,\dots,x_k=z,y,\quad t_0,t_1,\dots,t_{k-1},t
\end{equation*}
is an $(\varepsilon,T)$-chain from $x$ to $y$\@, and so $y\in\nbhd{T}$\@.
This gives
$\closure(\Phi(\tdomain_{\ge{}T}\times\nbhd{T}))\subset\nbhd{T}$\@, as
desired.
\end{proof}
\end{lemma}

\subsection{The decomposition theorem}

With the understanding of the connection between chains and trapping regions
from the preceding section, we can prove the Conley decomposition.  The
following lemma captures an essential part of the theorem.  Typically the
proof of this lemma is given separately for the continuous-time case and the
discrete-time case (where the notion of trapping region is taken to be our
notion of strong trapping region as in
Remark~\ref{rem:special-trap}--\ref{enum:special-trap1} and in
Definition~\ref{def:strong-trap} below).  Our proof works for both cases, and
has the additional benefit of being simpler than the already simple proofs in
each of the separate cases.
\begin{lemma}\label{lem:ChReccapbasin}
Let\/ $(\ts{X},\d)$ be a metric space and let\/ $\Phi$ be a topological flow
or semiflow on\/ $\ts{X}$\@.  If\/ $\nbhd{T}$ is a trapping region with\/ $A$
its corresponding attracting set, then
\begin{equation*}
\ChRec(\Phi)\cap\Orb^-(\nbhd{T})=A.
\end{equation*}
\begin{proof}
Let $T\in\tdomainp$ be such that
$\closure(\Phi(\tdomain_{\ge{}T}\times\nbhd{T}))\subset\nbhd{T}$ and, by
Lemma~\ref{lem:trap-erf}\@, let $\varepsilon\in\mappings[0]{\ts{X}}{\realp}$
be such that
\begin{equation*}
\oball[\d{}]{\varepsilon(\Phi(t,x))}{\Phi(t,x)}\subset\nbhd{T},
\qquad t\in\tdomain_{\ge{}T},\ x\in\nbhd{T}.
\end{equation*}
Let $t\in\tdomain_{\ge{}T}$ and let $m\in\integerp$\@.  Let
$x\in\ChRec(\Phi)\cap\nbhd{T}$\@.  Let
\begin{equation*}
x_0=x,x_1,\dots,x_k=x,\quad t_0,t_1,\dots,t_{k-1}
\end{equation*}
be a $(\min\{\frac{1}{m},\varepsilon\},t)$-chain from $x$ to itself.  We have
\begin{equation*}
\d(\Phi(t_0,x_0),x_1)<\min
\left\{\frac{1}{m},\varepsilon(\Phi(t,x_0))\right\}
\le\varepsilon(\Phi(t_0,x_0)),
\end{equation*}
which implies that $x_1\in\nbhd{T}$ by definition of $\varepsilon$\@.  Now
suppose that $x_j\in\nbhd{T}$ for $j\in\{0,1,\dots,k-2\}$\@, and note that,
as above,
\begin{equation*}
\d(\Phi(t_j,x_j),x_{j+1})<\varepsilon(\Phi(t_j,x_j)).
\end{equation*}
Thus $x_{j+1}\in\nbhd{T}$ by definition of $\varepsilon$\@.  Thus
$x_1,\dots,x_{k-1}\in\nbhd{T}$\@.  Now we have
\begin{equation*}
\d(\Phi(t_{k-1},x_{k-1}),x)<
\min\left\{\frac{1}{m},\varepsilon(\Phi(t_{k-1},x_{k-1}))\right\}
\le\frac{1}{m}.
\end{equation*}
Since $t_{k-1}\ge T+t$ and $x_{k-1}\in\nbhd{T}$\@, we have
\begin{equation*}
\Phi(t_{k-1},x_{k-1})\in\Phi(\tdomain_{\ge{}T+t}\times\nbhd{T}),
\end{equation*}
which implies that
\begin{equation*}
\dist_{\Phi(\tdomain_{\ge{}T+t}\times\nbhd{T})}(x)\le
\d(x,\Phi(t_{k-1},x_{k-1}))<\frac{1}{m}.
\end{equation*}
As this construction can be made for any $m\in\integerp$\@, we conclude that
$x\in\closure(\Phi(\tdomain_{\ge{}T+t}\times\nbhd{T}))$\@.  Therefore, as
this holds for every $t\in\tdomain_{\ge{}T}$\@, we have
\begin{equation*}
x\in\bigcap_{t\in\tdomainnn}\Phi(\tdomain_{\ge{}T+t}\times\nbhd{T})=A.
\end{equation*}

Let $x\in\ChRec(\Phi)\cap\Orb^-(\nbhd{T})$\@.  Let $(\varepsilon,T)$ be as in
part~\eqref{pl:trap-erf2} of Lemma~\ref{lem:trap-erf}\@.  Suppose that
$\Phi(t,x)\in\nbhd{T}$ for $t\in\tdomainp$\@.  Since
$x\in\ChRec(\Phi)$\@, let
\begin{equation*}
x_0=x,x_1,\dots,x_k=x,\quad t_0,t_1,\dots,t_{k-1}
\end{equation*}
be an $(\varepsilon,T+t)$-chain from $x$ to $x$\@.  Since
$\Phi_t(x)\in\nbhd{T}$\@,
$\Phi_{s+t}(x)\in\closure(\nbhd{T})\subset\nbhd{T}$ for every
$s\in\tdomain_{\ge{}T}$\@; in particular, $\Phi(t_0,x_0)\in\nbhd{T}$\@.  By
definition of $\varepsilon$\@, $x_1\in\nbhd{T}$\@.  Thus
\begin{equation*}
x_1,x_2,\dots,x_k=x,\quad t_1,t_2,\dots,t_{k-1}
\end{equation*}
is an $(\varepsilon,T)$ chain from $x_1\in\nbhd{T}$ to $x$\@.  As we argued
in the first part of the proof, this implies that $x_k=x\in\nbhd{T}$\@.  Thus
$x\in\ChRec(\Phi)\cap\nbhd{T}=A$\@, again from the first part of the proof.
\end{proof}
\end{lemma}

Now we can state the decomposition theorem.  In the statement of the result,
we denote by $\sT(\Phi)$ the set of trapping regions for a topological flow
or semiflow.
\begin{theorem}\label{the:conley-decomp}
Let\/ $(\ts{X},\d)$ be a metric space and let\/ $\Phi$ be a topological flow
or semiflow on\/ $\ts{X}$\@.  Then
\begin{equation*}
\ts{X}\setminus\ChRec(\Phi)=
\bigcup_{\nbhd{T}\in\sT(\Phi)}\Orb^-(\nbhd{T})\setminus A_{\nbhd{T}}.
\end{equation*}
\begin{proof}
Let $x\not\in\ChRec(\Phi)$\@.  Let
$\varepsilon\in\mappings[0]{\ts{X}}{\realp}$ and $T\in\tdomainp$ be such that
there is no $(\varepsilon,T)$-chain from $x$ to itself.  Let $\nbhd{T}$ be
the associated trapping region with $m=1$ as in Lemma~\ref{lem:chain->trap}
and let $A$ be the attracting set associated with $\nbhd{T}$\@.  Clearly
$x\not\in A\subset\nbhd{T}$\@.  Since
\begin{equation*}
x_0=x,x_1=\Phi_T(x),\quad T
\end{equation*}
is an $(\varepsilon,T)$ chain, $\Phi_T(x)\in\nbhd{T}$\@.  Therefore,
$x\in\Orb^-(\nbhd{T})$\@.  Thus
\begin{equation*}
\ts{X}\setminus\ChRec(\Phi)\subset\Orb^-(\nbhd{T})\setminus A
\subset\bigcup_{\nbhd{T}'\in\sT(\Phi)}\Orb^-(\nbhd{T}')\setminus A_{\nbhd{T}'}.
\end{equation*}

Now, let $x\in\Orb^-(\nbhd{T})\setminus A$ for some $\nbhd{T}\in\sT(\Phi)$
with $A$ its attracting set.  By Lemma~\ref{lem:ChReccapbasin}\@,
$x\not\in\ChRec(\Phi)$ and so
\begin{equation*}\eqqed
\bigcup_{\nbhd{T}'\in\sT(\Phi)}\Orb^-(\nbhd{T}')\setminus
A_{\nbhd{T}'}\subset
\ts{X}\setminus\ChRec(\Phi),
\end{equation*}
\end{proof}
\end{theorem}

There is an insightful rephrasing of this decomposition for flows.
\begin{corollary}
Let\/ $(\ts{X},\d)$ be a metric space and let\/ $\Phi$ be a topological flow
or semiflow on\/ $\ts{X}$\@.  Then
\begin{equation*}
\ChRec(\Phi)=
\bigcap_{\nbhd{T}\in\sT(\Phi)}\setdef{A\cup R}
{\textrm{$(A,R)$ is the attracting-repelling pair for\/
$\nbhd{T}\in\sT(\Phi)$}}.
\end{equation*}
\begin{proof}
From the theorem,
\begin{equation*}
\ChRec(\Phi)=\bigcap_{\nbhd{T}\in\sT(\Phi)}
\ts{X}\setminus(\Orb^-(\nbhd{T})\setminus A_{\nbhd{T}})=
\bigcap_{\nbhd{T}\in\sT(\Phi)}(\ts{X}\setminus\Orb^-(\nbhd{T}))\cup A_{\nbhd{T}}.
\end{equation*}
The corollary now follows from Lemma~\ref{lem:repelling-negorb}\@.
\end{proof}
\end{corollary}

\section{Complete Lyapunov functions}

In this section we give the second part of the Fundamental Theorem of
Dynamical Systems, namely the existence of a complete Lyapunov function.  Our
construction comes with a few steps.  First we work with discrete-time flows
and semiflows, with the final results being valid for flows and semiflows on
separable metric spaces.  In our constructions, we make use of strong
trapping regions (as in
Remark~\ref{rem:special-trap}--\ref{enum:special-trap1} and
Definition~\ref{def:strong-trap} below), following~\cite{MH:98}\@.  Despite
using strong trapping regions in place of the trapping regions from the
decomposition theorem, the conclusions refer only to chain notions which are
themselves not concerned with whether the trapping regions are strong or not.
After these constructions are complete, we show that they imply the existence
of complete Lyapunov functions in the continuous-time case for separable
metric spaces.

\subsection{Definitions}

We begin by giving definitions and elementary results around the notion of a
complete Lyapunov function.
\begin{definition}
Let $(\ts{X},\d)$ be a metric space and let $\Phi$ be a topological flow or
semiflow on $\ts{X}$\@.  A \defn{complete Lyapunov function} for $\Phi$ is a
continuous function $\map{L}{\ts{X}}{\interval[0,1]}$ such that:
\begin{compactenum}[(i)]
\item $L\scirc\Phi(t_2,x)\le L\scirc\Phi(t_1,x)$ for $x\in\ts{X}$ and
$t_1,t_2\in\tdomainnn$ satisfying $t_1<t_2$\@;
\item $L\scirc\Phi(t,x)<L(x)$ for $x\in\ts{X}\setminus\ChRec(\Phi)$ and
$t\in\tdomainp$\@;
\item $L\scirc\Phi(t,x)=L(x)$ for $x\in\ChRec(\Phi)$ and $t\in\tdomainnn$\@;
\item points $x,y\in\ChRec(\Phi)$ are chain equivalent if and only if
$L(x)=L(y)$\@.\oprocend
\end{compactenum}
\end{definition}

Thus we see that a complete Lyapunov function (1)~distinguishes the forward
chain recurrent set, (2)~additionally distinguishes the chain components, and
(3)~tells us something about the flow or semiflow between the chain
components.

Complete Lyapunov functions for flows have particular properties that are
sometimes also casually (and falsely) ascribed to complete Lyapunov functions
for semiflows.  The following lemma gives two particular properties of
complete Lyapunov functions for flows.
\begin{lemma}\label{lem:clf-flow}
Let\/ $(\ts{X},\d)$ be a metric space and let\/ $\Phi$ be a topological flow
on\/ $\ts{X}$ with\/ $L$ a complete Lyapunov function for\/ $\Phi$\@.  Then
the following statements hold:
\begin{compactenum}[(i)]
\item \label{pl:clf-flow1} if\/ $x\in\ts{X}\setminus\ChRec(\Phi)$\@, then the
function $t\mapsto L\scirc\Phi(t,x)$ is strictly decreasing on\/
$\tdomain$\@;
\item \label{pl:clf-flow2} if\/ $x\in\ChRec(\Phi)$\@, then\/
$L\scirc\Phi(t,x)=L(x)$ for\/ $t\in\tdomain$\@.
\end{compactenum}
\begin{proof}
\eqref{pl:clf-flow1} Note that $\ChRec(\Phi)$ is invariant by
Proposition~\ref{prop:ChRec-props}\eqref{pl:ChRec1}\@.  By
Lemma~\ref{lem:comp-inv}\@, $\ts{X}\setminus\ChRec(\Phi)$ is also invariant.
Therefore, if $x\in\ts{X}\setminus\ChRec(\Phi)$ and if $t_1,t_2\in\tdomain$
satisfy $t_1<t_2$\@, we have
\begin{equation*}
L\scirc\Phi(t_2,x)=L\scirc\Phi(t_2-t_1,\Phi(t_1,x))<
L\scirc\Phi(t_1,x),
\end{equation*}
since $\Phi(t_1,x)\in\ts{X}\setminus\ChRec(\Phi)$\@.

\eqref{pl:clf-flow2} We need only prove this for $t\in\tdomainn$\@, so let
$x\in\ChRec(\Phi)$ and $t\in\tdomainn$\@.  Then $\Phi(t,x)\in\ChRec(\Phi)$ by
Proposition~\ref{prop:ChRec-props}\eqref{pl:ChRec1}\@.  By definition,
\begin{equation*}\eqqed
L\scirc\Phi(t,x)=L\scirc\Phi(-t,\Phi(t,x))=L(x).
\end{equation*}
\end{proof}
\end{lemma}

We shall show that every topological flow or semiflow on a separable metric
space possesses a complete Lyapunov function.  To do so will require some
work.  First we will establish the result in the discrete-time case, and then
use this result to establish the result in the continuous-time case.  That
the continuous- and discrete-time cases should be related should come as no
surprise, given that the chain recurrent set and its chain components are
essentially determined by the mapping $\Phi_1$\@, given
Theorem~\ref{the:chain-equiv}\@.

\subsection{Some constructions particular to the discrete-time case}

In this section and the two sections following it, we work with a
discrete-time flow or semiflow $\Phi$ on a metric space $(\ts{X},\d)$\@, and
we denote $\phi=\Phi_1\in\mappings[0]{\ts{X}}{\ts{X}}$\@.  Where convenient,
we will drop the reference to $\Phi$ and simply refer to $\phi$\@.  For
instance, we will write $\ChRec(\phi)\eqdef\ChRec(\Phi)$\@.

First, it is useful to give a refined equivalent characterisation of complete
Lyapunov functions in this case.
\begin{lemma}
Let\/ $(\ts{X},\d)$ be a metric space and let\/ $\phi=\Phi_1$\@.  A mapping\/
$L\in\mappings[0]{\ts{X}}{\interval[0,1]}$ is a complete Lyapunov function if
and only if
\begin{compactenum}[(i)]
\item $L\scirc\phi(x)\le L(x)$ for\/ $x\in\ts{X}$\@,
\item $L\scirc\phi(x)<L(x)$ for\/ $x\in\ts{X}\setminus\ChRec(\Phi)$\@,
\item $L\scirc\phi(x)=L(x)$ for\/ $x\in\ChRec(\Phi)$\@, and
\item points\/ $x,y\in\ChRec(\Phi)$ are chain equivalent if and only if\/
$L(x)=L(y)$\@.
\end{compactenum}
\begin{proof}
Since the fourth property in the definition of a complete Lyapunov function
and in the statement of the lemma are the same, we need only work with the
first three properties in each case.

First suppose that $L$ is a complete Lyapunov function.  Then, since $0<1$\@,
$L\scirc\phi(x)\le L(x)$ for every $x\in\ts{X}$\@.  Similarly,
$L\scirc\phi(x)<L(x)$ for $x\in\ts{X}\setminus\ChRec(\Phi)$\@.  If
$x\in\ChRec(\Phi)$\@, then $\phi(x)\in\ChRec(\Phi)$ by
Proposition~\ref{prop:ChRec-props}\eqref{pl:ChRec1}\@.  Therefore,
$L\scirc\phi(x)=L(x)$\@.  This gives the first three properties in the
statement of the lemma.

Next suppose that $L$ has the first three properties from the statement of
the lemma.  Since $L\scirc\Phi(t,x)=L\scirc\phi^t(x)$\@, an elementary
induction gives each of the first three properties in the definition of a
complete Lyapunov function.
\end{proof}
\end{lemma}

We shall work with $\varepsilon$-$1$-chains for
$\varepsilon\in\mappings[0]{\ts{X}}{\realp}$\@, which we simply call
$\varepsilon$-chains according to Definition~\ref{def:vareps-chain}\@.  It
turns out to be convenient to modify slightly our notion of chain.  As
pointed out by \citet{MH:98}\@, this minor modification is required by the
condition that complete Lyapunov functions \emph{decrease}\@.  If we were
happy with them being allowed to \emph{increase}\@, we could get by with our
existing notion of chains.
\begin{definition}
Let $(\ts{X},\d)$ be a metric space and let
$\phi\in\mappings[0]{\ts{X}}{\ts{X}}$\@.  A \defn{$\sigma$chain} for $\phi$
is a finite sequence
\begin{equation*}
(x_0,y_0),(x_1,y_1),\dots,(x_k,y_k)
\end{equation*}
of ordered pairs in $\nbhd{X}$ such that $x_{j+1}=\phi(y_j)$\@,
$j\in\{0,1,\dots,k-1\}$\@.  The nonnegative integer $k$ is the \defn{length}
of the $\sigma$chain and the $\sigma$chain is said to be \defn{from $x_0$ to
$y_k$}\@.  If $\varepsilon\in\mappings[0]{\nbhd{X}}{\realp}$\@, then a
$\sigma$chain
\begin{equation*}
(x_0,y_0),(x_1,y_1),\dots,(x_k,y_k)
\end{equation*}
is an \defn{$\varepsilon$-$\sigma$chain} for $\phi$ if
$\d(x_j,y_j)<\varepsilon(x_j)$\@, $j\in\{0,1,\dots,k\}$\@.\oprocend
\end{definition}

Let us clarify the relationship between $\varepsilon$-chains and
$\varepsilon$-$\sigma$chains.
\begin{lemma}
Let\/ $(\ts{X},\d)$ be a metric space and let\/
$\phi\in\mappings[0]{\ts{X}}{\ts{X}}$\@.  Then the following statements hold:
\begin{compactenum}[(i)]
\item \label{pl:eps-epsJ1} if\/
$\varepsilon\in\mappings[0]{\ts{X}}{\realp}$\@, then there exists\/
$\varepsilon'\in\mappings[0]{\ts{X}}{\realp}$ such that, if
\begin{equation*}
(x_0,y_0),(x_1,y_1),\dots,(x_k,y_k)
\end{equation*}
is an\/ $\varepsilon'$-$\sigma$chain for\/ $\phi$\@, then
\begin{equation*}
x_0,x_1,\dots,x_k
\end{equation*}
is an\/ $\varepsilon$-chain for $\phi$\@;
\item \label{pl:eps-epsJ2} if\/
$\varepsilon\in\mappings[0]{\ts{X}}{\realp}$ and if
\begin{equation*}
x_0,x_1,\dots,x_k
\end{equation*}
is a\/ $\varepsilon$-chain for\/ $\phi$\@, then
\begin{equation*}
(x_0,x_0),(\phi(x_0),x_1),\dots,(\phi(x_{k-1}),x_k)
\end{equation*}
is a\/ $\varepsilon$-$\sigma$chain for\/ $\phi$\@.
\end{compactenum}
\begin{proof}
\eqref{pl:eps-epsJ1} By Lemma~\ref{lem:general-epsilon-delta}\@, let
$\varepsilon'\in\mappings[0]{\ts{X}}{\realp}$ be such that
\begin{equation*}
\d(x,y)<\varepsilon'(x)\implies\d(\phi(x),\phi(y))<\varepsilon(x).
\end{equation*}
This $\delta$ can be verified to have the asserted property.

\eqref{pl:eps-epsJ2} This is clear by definition.
\end{proof}
\end{lemma}

We shall also work with a modified notion of trapping region.  Specifically,
we will work with the \emph{usual} notion of trapping region from the
literature in the discrete-time case, as discussed in
Remark~\ref{rem:special-trap}--\ref{enum:special-trap1}\@.
\begin{definition}\label{def:strong-trap}
Let $(\ts{X},\d)$ be a metric space and let
$\phi\in\mappings[0]{\ts{X}}{\ts{X}}$\@.  A \defn{strong trapping region} for
$\phi$ is a subset $\nbhd{T}$ for which $\closure(\phi(\nbhd{T}))\subset\interior(\nbhd{T})$\@.\oprocend
\end{definition}

Clearly a strong trapping region is a trapping region.  Thus a strong
trapping region has an associated attracting set and, in the case of flows,
an associated repelling set.

\subsection{The discrete-time case: weak Lyapunov functions for attracting
sets}

In this part of our development, we do much of the technical heavy lifting.
We work with a fixed strong trapping region and build a weak Lyapunov
function (\ie~a function that is nonincreasing along trajectories) with
particular properties relative to the associated attracting set.

Our construction of a weak Lyapunov function for an attracting set is done in
stages. First, given $\varepsilon\in\mappings[0]{\ts{X}}{\realp}$ and a
$\sigma$chain
\begin{equation*}
(x_0,y_0),(x_1,y_1),\dots,(x_k,y_k)
\end{equation*}
for $\phi$ which we denote by $\gamma$\@, we define
\begin{equation*}
D_\varepsilon(\gamma)=\sum_{j=0}^k\frac{\d(x_j,y_j)}{\varepsilon(x_j)}.
\end{equation*}
Now, if $C\subset\ts{X}$ is a nonempty closed set, if $x\in\ts{X}$\@, and if
$\phi\in\mappings[0]{\ts{X}}{\ts{X}}$\@, denote by $\sC_\phi(C;x)$ the set of
all $\sigma$chains for $\phi$ from a point in $C$ to $x$\@.  Then, for
$\varepsilon\in\mappings[0]{\ts{X}}{\realp}$\@, define
\begin{equation*}
E_{\phi,\varepsilon}(C;x)=\inf\setdef{D_\varepsilon(\gamma)}
{\gamma\in\sC_\phi(C;x)}.
\end{equation*}
It is evident that, since $(x,x)\in\sC_\phi(C;x)$ if $x\in C$\@, that
$E_{\phi,\varepsilon}(C;x)=0$ in this case.

The following properties of $E_{\phi,\varepsilon}(C;x)$ are important for us.
\begin{lemma}\label{lem:Ephieps}
Let\/ $(\ts{X},\d)$ be a metric space and let\/
$\phi\in\mappings[0]{\ts{X}}{\ts{X}}$\@.  Let\/ $C\subset\ts{X}$ be a
nonempty closed set, let\/ $x\in\ts{X}$\@, and let\/
$\varepsilon\in\mappings[0]{\ts{X}}{\realp}$\@.  Then the following
statements hold:
\begin{compactenum}[(i)]
\item \label{pl:Ephieps1} $E_{\phi,\varepsilon}(C;\phi(x))\le
E_{\phi,\varepsilon}(C;x)$\@;
\item \label{pl:Ephieps2} the mapping\/ $x\mapsto E_{\phi,\varepsilon}(C;x)$
is continuous.
\end{compactenum}
\begin{proof}
\eqref{pl:Ephieps1} Let
\begin{equation*}
(x_0,y_0),(x_1,y_1),\dots,(x_k,x)
\end{equation*}
be a $\sigma$chain from $x_0\in C$ to $x$\@, denoted by $\gamma$\@.  Define a
$\sigma$chain $\gamma'$ from $x_0\in C$ to $\phi(x)$ by
\begin{equation*}
(x_0,y_0),(x_1,y_1),\dots,(x_k,x),(\phi(x),\phi(x)).
\end{equation*}
Since $D_\varepsilon(\gamma)=D_\varepsilon(\gamma')$\@, we must have
$E_{\phi,\varepsilon}(C;\phi(x))\le E_{\phi,\varepsilon}(C;x)$\@.

\eqref{pl:Ephieps2} Suppose first that $x\in C$\@.  In this case,
$E_{\phi,\varepsilon}(C;x)=0$ (as we observed just prior to the statement of
the lemma).  For $y\in C$\@, we have the $\sigma$chain $\gamma$ given by
\begin{equation*}
(y,x)
\end{equation*}
from a point in $C$ to $x$\@.  Note that
\begin{equation*}
E_{\phi,\varepsilon}(C;y)=D_\varepsilon(\gamma)=
\frac{\d(x,y)}{\varepsilon(x)},
\end{equation*}
and so $\lim_{y\to
x}E_{\phi,\varepsilon}(C;y)=0=E_{\phi,\varepsilon}(C;x)$\@, giving continuity
at $x$\@.

Now suppose that $x\not\in C$\@.  Consider a $\sigma$chain $\gamma$
\begin{equation*}
(x_0,y_0),(x_1,y_1),\dots,(x_k,x)
\end{equation*}
from a point in $C$ to $x$\@.  Let $y\in\ts{X}$ and consider the
$\sigma$chain $\gamma'$
\begin{equation*}
(x_0,y_0),(x_1,y_1),\dots,(x_k,y)
\end{equation*}
from the same point in $C$ to $y$\@.  We then have
\begin{equation}\label{eq:Ephieps1}
\snorm{D_\varepsilon(\gamma)-D_\varepsilon(\gamma')}=
\asnorm{\frac{\d(x_k,x)}{\varepsilon(x_k)}-\frac{\d(x_k,y)}{\varepsilon(x_k)}}
\le\frac{\d(x,y)}{\varepsilon(x_k)}
\end{equation}
using a standard metric identity~\cite[Theorem~1.1.2]{MOS:07}\@.  We claim
now that this shows that $x\mapsto E_{\phi,\varepsilon}(C;x)$ is upper
semicontinuous.  To see this, let $\delta\in\realp$ and let $\gamma$ be a
$\sigma$chain from a point in $C$ to $x$ such that
\begin{equation*}
D_\varepsilon(\gamma)<E_{\phi,\varepsilon}(C;x)+\frac{\delta}{2}.
\end{equation*}
For $y\in\ts{X}$ sufficiently close to $x$ and with $\gamma'$ as defined
above, we can ensure that
\begin{equation*}
D_\varepsilon(\gamma')<D_\varepsilon(\gamma')+\frac{\delta}{2}.
\end{equation*}
With $\gamma'$ so chosen, we have
\begin{equation*}
E_{\phi,\varepsilon}(C;y)\le D_\varepsilon(\gamma')<
D_\varepsilon(\gamma)+\frac{\delta}{2}<
E_{\phi,\varepsilon}(C;x)+\delta,
\end{equation*}
from which we can conclude the asserted upper semicontinuity.

Let us suppose that we do not have lower semicontinuity at $x$\@.  Then there
exists a sequence $\ifam{y_j}_{j\in\integerp}$ in $\ts{X}$ converging to $x$
and $\beta\in\realp$ such that
\begin{equation*}
E_{\phi,\varepsilon}(C;y_j)<E_{\phi,\varepsilon}(C;x)-\beta,
\qquad j\in\integerp.
\end{equation*}
If this last condition holds, then there exists a sequence
$\gamma'_j\in\sC_\phi(C;y_j)$\@, $j\in\integerp$\@, such that
\begin{equation*}
D_\varepsilon(\gamma'_j)<E_{\phi,\varepsilon}(C;x)-\beta,\qquad j\in\integerp.
\end{equation*}
Denote by $\gamma_j$ the $\sigma$chain obtained as above, replacing $y_j$
with $x$\@.  As in~\eqref{eq:Ephieps1}\@, we have
\begin{equation*}
\snorm{D_\varepsilon(\gamma_j)-D_\varepsilon(\gamma'_j)}\le
\frac{\d(x,y_j)}{\varepsilon(z_j)},
\end{equation*}
where the last pair in $\gamma'_j$ is $(z_j,y_j)$\@.  If it holds that
$\liminf_{j\to\infty}\frac{\d(x,y_j)}{\varepsilon(z_j)}=0$\@, then we arrive
at the contradiction that $D_\varepsilon(\gamma_j)<E_{\phi,\varepsilon}(C;x)$
for some sufficiently large $j$\@, in contradiction to the definition of
$E_{\phi,\varepsilon}(C;x)$\@.  This contradiction would then show that our
assumption that we do not have lower semicontinuity at $x$ must be false.

Thus it remains to show that
$\liminf_{j\to\infty}\frac{\d(x,y_j)}{\varepsilon(z_j)}=0$ under the
assumption that
\begin{equation*}
E_{\phi,\varepsilon}(C;y_j)<E_{\phi,\varepsilon}(C;x)-\beta,
\qquad j\in\integerp.
\end{equation*}
To this end, note that, as we have seen, the assumption implies that
\begin{equation*}
E_{\phi,\varepsilon}(C;x)-\beta>D_\varepsilon(\gamma'_j)\ge
\frac{\d(z_j,y_j)}{\varepsilon(z_j)}.
\end{equation*}
Therefore,
\begin{equation*}
\varepsilon(z_j)\ge M\d(z_j,y_j),\qquad j\in\integerp,
\end{equation*}
for a constant $M\in\realp$ independent of $j$\@.  Therefore,
\begin{equation*}
\frac{\d(x,y_j)}{\varepsilon(z_j)}\le M^{-1}\frac{\d(x,y_j)}{\d(z_j,y_j)},
\qquad j\in\integerp.
\end{equation*}
We now consider two cases: (1)~$\lim_{j\to0}\d(z_j,y_j)=0$ and (2)~(1) does
not hold.  In the first case, $\lim_{j\to\infty}z_j=x$\@.  Therefore, for $j$
sufficiently large, $\varepsilon(z_j)\ge\frac{1}{2}\varepsilon(x)$\@, in
which case we directly have
\begin{equation*}
\lim_{j\to\infty}\frac{\d(x,y_j)}{\varepsilon(z_j)}=0,
\end{equation*}
giving the desired conclusion in this case.  In the second case, we have
\begin{equation*}
\liminf_{j\to\infty}\frac{\d(x,y_j)}{\varepsilon(z_j)}\le
\liminf_{j\to\infty}M^{-1}\frac{\d(x,y_j)}{\d(z_j,y_j)}=0,
\end{equation*}
again giving the desired result.
\end{proof}
\end{lemma}

The way to view the function $x\mapsto E_\varepsilon(C;x)$ is that it
measures the ``$\varepsilon$-effort'' expended along any
$\varepsilon$-$\sigma$chain going from a point in $C$ to the point $x$\@,
where ``$\varepsilon$-effort'' is characterised by the ratio of the jumps of
a $\sigma$chain compared to the maximum jump permitted by $\varepsilon$\@.
(Thus ``$1$'' represents neutral effort.)  The closed set $C$ in the
definition has no relationship to the dynamics.

In the next stage in our construction, we fix a strong trapping region
$\nbhd{T}$ and consider our previous ``minimum effort'' function applied to
the sets $\closure(\phi^k(\nbhd{T}))$\@, $k\in\integerp$\@.  We tailor
$\varepsilon$ according to Lemma~\ref{lem:trap-erf}\@, and furthermore to be
bounded above by $1$\@.  Throughout the next part of our construction, we
understand this $\varepsilon$ to have been chosen and fixed.  With such a
$\varepsilon$ at hand and for $k\in\integerp$\@, we define
\begin{equation*}
E_{\varepsilon,k}(x)=E_{\phi^k,\varepsilon}(\closure(\phi^k(\nbhd{T}));x),
\qquad x\in\ts{X}.
\end{equation*}
Let us list the pertinent properties of $E_{\varepsilon,k}$\@.
\begin{lemma}\label{lem:Eepsk}
Let\/ $(\ts{X},\d)$ be a metric space and let\/
$\phi\in\mappings[0]{\ts{X}}{\ts{X}}$\@.  For a strong trapping region\/
$\nbhd{T}$ and with\/ $\varepsilon$ as above, and for\/ $k\in\integerp$\@,
the following statements hold:
\begin{compactenum}[(i)]
\item \label{pl:Eepsk1} $E_{\varepsilon,k}$ is nonnegative and continuous;
\item \label{pl:Eepsk2}
$E_{\varepsilon,k}\scirc\phi^k(x)\le E_{\varepsilon,k}(x)$\@,\/
$x\in\ts{X}$\@;
\item \label{pl:Eepsk3} if
\begin{equation*}
(x_0,y_0),(x_1,y_1),\dots,(x_k,y_k)
\end{equation*}
is an\/ $\varepsilon$-$\sigma$chain for\/ $\phi^k$ and if\/
$x_0\in\closure(\phi^k(\nbhd{T}))$\@, then\/ $y_k\in\nbhd{T}$\@;
\item \label{pl:Eepsk4}
$E_{\varepsilon,k}^{-1}(0)=\closure(\phi^k(\nbhd{T}))$\@.\savenum
\end{compactenum}
Additionally, no longer fixing\/ $k$\@,
\begin{compactenum}[(i)]\resumenum
\item \label{pl:Eepsk5} if\/ $A$ is the attracting set associated with\/
$\nbhd{T}$ and if\/ $x\in A$\@, then\/ $E_{\varepsilon,k}(x)=0$ for every\/
$k\in\integerp$\@;
\item \label{pl:Eepsk6} if\/ $x\not\in\nbhd{T}$\@, then\/
$E_{\varepsilon,k}(x)\ge1$ for every\/ $k\in\integerp$\@.
\end{compactenum}
\begin{proof}
Parts~\eqref{pl:Eepsk1} and~\eqref{pl:Eepsk2} follow directly from the
definitions and Lemma~\ref{lem:Ephieps}\@.

\eqref{pl:Eepsk3} Since $x_0\in\closure(\phi^k(\nbhd{T}))\subset\nbhd{T}$\@,
then $y_0\in\nbhd{T}$ by the fact that $\varepsilon$ satisfies the condition
of Lemma~\ref{lem:trap-erf}\@.  Then
$x_1=\phi^k(y_0)\in\closure(\phi^k(\nbhd{T}))$ and so $y_1\in\nbhd{T}$ by the
same argument.  Inductively, $y_k\in\nbhd{T}$\@.

\eqref{pl:Eepsk4} As we have observed above when working with
$E_{\phi,\varepsilon}(C;x)$\@, we have
\begin{equation*}
\closure(\phi^k(\nbhd{T}))\subset E_{\varepsilon,k}^{-1}(0).
\end{equation*}
It suffices, then, to show that $E_{\varepsilon,k}(x)>0$ for
$x\not\in\closure(\phi^k(\nbhd{T}))$\@.  If there are no
$\varepsilon$-$\sigma$chains from a point in $\closure(\phi^k(\nbhd{T}))$ to
$x$\@, then $E_{\varepsilon,k}(x)\ge1$\@, as we shall see in our proof of
part~\eqref{pl:Eepsk6} below.  So suppose that
\begin{equation*}
(x_0,y_0),(x_1,y_1),\dots,(x_k,x)
\end{equation*}
is a $\varepsilon$-$\sigma$chain from a point
$x_0\in\closure(\phi^k(\nbhd{T}))$ to $x$\@.  By part~\eqref{pl:Eepsk3}\@,
$x_k\in\closure(\phi^k(\nbhd{T}))$\@, whereupon
\begin{equation*}
D_\varepsilon(\gamma)\ge\frac{\d(x_k,x)}{\varepsilon(x_k)}
\ge\d(x_k,x)\ge\dist_{\closure(\phi^k(\nbhd{T}))}(x)>0,
\end{equation*}
noting that we have assumed that $\varepsilon$ takes values in
$\interval(0,1]$\@.

\eqref{pl:Eepsk5} If $x\in A$\@, then $x\in\closure(\phi^k(\nbhd{T}))$ for
$k\in\integerp$\@.  As we argued in the previous part of the proof, this
implies that $E_{\varepsilon,k}(x)=0$ for all $k\in\integerp$\@.

\eqref{pl:Eepsk6} If $x\not\in\nbhd{T}$ and if
$\gamma\in\sC_\phi(\closure(\phi^k(\nbhd{T}));x)$\@, then $\gamma$ is not an
$\varepsilon$-$\sigma$chain by part~\eqref{pl:Eepsk3}\@.  Therefore, if
$\gamma$ is given by
\begin{equation*}
(x_0,y_0),(x_1,y_1),\dots,(x_k,y_k),
\end{equation*}
then, for some $j\in\{0,1,\dots,k\}$\@, $\d(x_j,y_j)\ge\varepsilon(x_j)$\@,
from which we have
\begin{equation*}
D_\varepsilon(\gamma)\ge\frac{\d(x_j,y_j)}{\varepsilon(x_j)}\ge1.
\end{equation*}
Therefore, $E_{\phi^k,\varepsilon}(\closure(\phi^k(\nbhd{T}));x)\ge1$\@.
\end{proof}
\end{lemma}

The next step in our construction is to average $E_{\varepsilon,k}$ over the
first $k$ iterates of $\phi$\@:
\begin{equation*}
\ol{E}_{\varepsilon,k}(x)=
\frac{1}{k}\sum_{j=0}^{k-1}E_{\varepsilon,k}\scirc\phi^j(x).
\end{equation*}
Let us record the salient properties of this function.
\begin{lemma}\label{lem:olEepsk}
Let\/ $(\ts{X},\d)$ be a metric space and let\/
$\phi\in\mappings[0]{\ts{X}}{\ts{X}}$ with strong trapping region\/
$\nbhd{T}$ and associated error function\/ $\varepsilon$ as in
Lemma~\ref{lem:trap-erf}\@.  Then the following statements hold:
\begin{compactenum}[(i)]
\item \label{pl:olEepsk1} $\ol{E}_{\varepsilon,k}$ is continuous and\/
$\realnn$-valued;
\item \label{pl:olEepsk2}
$\ol{E}_{\varepsilon,k}\scirc\phi(x)\le\ol{E}_{\varepsilon,k}(x)$\@.
\end{compactenum}
\begin{proof}
\eqref{pl:olEepsk1} This is clear.

\eqref{pl:olEepsk2} We have
\begin{equation*}\eqqed
\ol{E}_{\varepsilon,k}\scirc\phi(x)-\ol{E}_{\varepsilon,k}(x)=
\frac{1}{k}(E_{\varepsilon,k}\scirc\phi(x)-E_{\varepsilon,k}(x))\le0.
\end{equation*}
\end{proof}
\end{lemma}

Now we can give the final part of the constructions we shall make for a fixed
strong trapping region $\nbhd{T}$\@.  We define
\begin{equation*}
E_{\nbhd{T}}(x)=\sum_{k=1}^\infty\frac{\min\{\ol{E}_{\varepsilon,k}(x),1\}}{2^k},
\qquad x\in\ts{X},
\end{equation*}
and then
\begin{equation*}
L_{\nbhd{T}}(x)=\sum_{j=0}^\infty\frac{E_{\nbhd{T}}\scirc\phi^j(x)}{2^j},
\qquad x\in\ts{X}.
\end{equation*}
Let us record the properties of the second of these functions.  In the
statement and proof of the lemma (and of many constructions in the remainder
of this section), it is insightful to keep in mind the characterisation of
repelling sets from Lemma~\ref{lem:repelling-negorb}\@.
\begin{lemma}\label{lem:LT}
Let\/ $(\ts{X},\d)$ be a metric space and let\/
$\phi\in\mappings[0]{\ts{X}}{\ts{X}}$\@.  Let\/ $\nbhd{T}$ be an open strong
trapping region for\/ $\phi$ with\/ $A$ the associated attracting set.  Then
the following statements hold:
\begin{compactenum}[(i)]
\item \label{pl:LT1} $L_{\nbhd{T}}$ is continuous with values in\/
$\interval[0,1]$\@;
\item \label{pl:LT2} $L_{\nbhd{T}}^{-1}(0)=A$ and\/
$L_{\nbhd{T}}^{-1}(1)=\ts{X}\setminus\Orb^-(\nbhd{T})$\@;
\item \label{pl:LT3} $L_{\nbhd{T}}\scirc\phi(x)<L_{\nbhd{T}}(x)$ for\/
$x\in\Orb^-(\nbhd{T})\setminus A$\@.
\end{compactenum}
\begin{proof}
\eqref{pl:LT1} The series defining the functions $E_{\nbhd{T}}$ is
uniformly convergent by the Weierstrass $M$-test, by virtue of which the
limit function is continuous.  The function $E_{\nbhd{T}}$ take values in
$\interval[0,1]$ since $\sum_{j=1}^\infty\frac{1}{2^j}=1$\@.  For the same
reason, the function $L_{\nbhd{T}}$ take values in $\interval[0,1]$\@.

\eqref{pl:LT2} If $x\in A$\@, then $\phi^j(x)\in A$ for every
$j\in\integernn$ by
Proposition~\ref{prop:attract-repell-props}\eqref{pl:attract-repell-props2}\@.
Therefore, $E_{\varepsilon,k}\scirc\phi^j(x)=0$ for every $k\in\integerp$ and
$j\in\integernn$ by Lemma~\ref{lem:Eepsk}\eqref{pl:Eepsk5}\@.  Therefore,
$\ol{E}_{\varepsilon,k}\scirc\phi^j(x)=0$ for every $k\in\integerp$ and
$j\in\integernn$ and so $E_{\nbhd{T}}\scirc\phi^j(x)=0$ for every
$j\in\integernn$\@.  Therefore, $L_{\nbhd{T}}(x)=0$\@.  Conversely, if
$L_{\nbhd{T}}(x)=0$\@, then $E_{\nbhd{T}}\scirc\phi^j(x)=0$ for
$j\in\integernn$\@.  This, in turn, implies that
$\ol{E}_{\varepsilon,k}\scirc\phi^j(x)=0$ for every $k\in\integerp$ and
$j\in\integernn$\@.  Considering the case of $j=0$\@, this implies that
$E_{\varepsilon,k}(\phi^l(x))=0$ for every $k\in\integerp$ and
$l\in\{0,1,\dots,k-1\}$\@.  In particular, $E_{\varepsilon,k}(x)=0$ for every
$k\in\integerp$\@.  Thus $x\in\cap_{k\in\integerp}\closure(\phi^k(\nbhd{T}))$
by Lemma~\ref{lem:Eepsk}\eqref{pl:Eepsk4}\@.  Thus $x\in A$\@.

Now suppose that $x\in\ts{X}\setminus\Orb^-(\nbhd{T})$ and note that
$\phi^j(x)\not\in\nbhd{T}$ for every $j\in\integernn$\@.  Therefore,
$E_{\varepsilon,k}\scirc\phi^j(x)\ge1$ for every $k\in\integerp$ and every
$j\in\integernn$ by Lemma~\ref{lem:Eepsk}\eqref{pl:Eepsk6}\@.  Therefore,
$\ol{E}_{\varepsilon,k}(x)\ge1$\@, and so $E_{\nbhd{T}}\scirc\phi^j(x)=1$ for
every $k\in\integerp$ and $j\in\integernn$\@.  Thus $L_{\nbhd{T}}(x)=1$\@.

Conversely, suppose that $L_{\nbhd{T}}(x)=1$\@.  It follows that
$E_{\nbhd{T}}\scirc\phi^j(x)=1$ for $j\in\integernn$ and so
$\ol{E}_{\varepsilon,k}\scirc\phi^j(x)\ge1$ for every $k\in\integerp$ and
$j\in\integernn$\@.  We claim that this prohibits $x\in\Orb^-(\nbhd{T})$ for
some $y\in\nbhd{T}$\@.  Indeed, if $x\in\Orb^-(y)$ for some $y\in\nbhd{T}$\@,
then $\phi^j(x)\in\nbhd{T}$ for some $j\in\integernn$ and so
$\phi^{j+k}(x)\in\phi^{j+k}(\nbhd{T})\subset\closure(\phi^{j+k}(\nbhd{T}))$
for $k\in\integerp$\@.  By Lemma~\ref{lem:Eepsk}\eqref{pl:Eepsk4}\@, this
implies that $E_{\varepsilon,k}\scirc\phi^j(x)=0$\@.  This, in turn, implies
that $\ol{E}_{\varepsilon,k}(x)<1$\@.  This shows that, indeed, if
$L_{\nbhd{T}}(x)=1$\@, then $x\not\in\Orb(y)$ for every $y\in\nbhd{T}$\@.
Thus $x\in\ts{X}\setminus\Orb^-(\nbhd{T})$\@.

\eqref{pl:LT3} By Lemma~\ref{lem:olEepsk}\eqref{pl:olEepsk2}\@,
\begin{equation*}
\ol{E}_{\varepsilon,k}\scirc\phi^{j+1}(x)\le\ol{E}_{\varepsilon,k}\scirc\phi^j(x),
\qquad x\in\ts{X},\ k\in\integerp,\ j\in\integernn.
\end{equation*}
It follows that
\begin{equation*}
E_{\nbhd{T}}(\phi^{j+1}(x))\le E_{\nbhd{T}}(\phi^j(x)),\qquad
x\in\ts{X},\ j\in\integernn.
\end{equation*}
We will have $L_{\nbhd{T}}\scirc\phi(x)<L_{\nbhd{T}}(x)$ when there exists
$j\in\integernn$ such that
\begin{equation*}
E_{\nbhd{T}}(\phi^{j+1}(x))<E_{\nbhd{T}}(\phi^j(x)).
\end{equation*}
This will happen when there exists $k\in\integerp$ and $j\in\integernn$ such
that
\begin{compactenum}
\item \label{enum:LT1} $\ol{E}_{\varepsilon,k}\scirc\phi^j(x)<1$ and
\item \label{enum:LT2}
$\ol{E}_{\varepsilon,k}\scirc\phi^{j+1}(x)<\ol{E}_{\varepsilon,k}\scirc\phi^j(x)$\@.
\end{compactenum}
Therefore, to establish this part of the lemma, it suffices to show that, for
$x\in\Orb^-(\nbhd{T})$\@, there exists $k\in\integerp$ and $j\in\integernn$
such that conditions~\ref{enum:LT1} and~\ref{enum:LT2} are satisfied.

First suppose that $x\in\phi(\nbhd{T})\setminus A$\@.  Since $x\not\in A$\@,
there exists $k\in\integerp$ such that
$x\not\in\closure(\phi^{k+1}(\nbhd{T}))$ but $x\in\closure(\phi^j(\nbhd{T}))$
for $j\in\{1,\dots,k\}$\@.  Thus $\ol{E}_{\varepsilon,j}(x)=0$ for
$j\in\{1,\dots,k\}$ and $\ol{E}_{\varepsilon,k+1}(x)>0$ by
Lemma~\ref{lem:Eepsk}\eqref{pl:Eepsk4}\@.  Since
$\phi(x)\in\closure(\phi^j(\nbhd{T}))$ for $j\in\{1,\dots,k+1\}$\@, we can
again use Lemma~\ref{lem:Eepsk}\eqref{pl:Eepsk4} to see that
\begin{equation*}
\ol{E}_{\varepsilon,k+1}\scirc\phi(x)=0<\ol{E}_{\varepsilon,k+1}(x),
\end{equation*}
giving this part of the lemma in the case that
$x\ni\phi(\nbhd{T})\setminus A$\@.  Now suppose that
$x\in\Orb^-(\nbhd{T})\setminus\phi(\nbhd{T})$\@.  Then
$\phi^j(x)\in\phi(\nbhd{T})$ for some $j\in\integerp$\@, and then the above
argument gives
\begin{equation*}
\ol{E}_{\varepsilon,k+1}\scirc\phi^{j+1}(x)=0<
\ol{E}_{\varepsilon,k+1}\scirc\phi^j(x)
\end{equation*}
for some $k\in\integerp$\@.  This gives the result.
\end{proof}
\end{lemma}

\subsection{The discrete-time case: the complete Lyapunov function}

In the preceding part of our construction, we fixed a strong trapping region
$\nbhd{T}$ and constructed a weak Lyapunov function $L_{\nbhd{T}}$ with some
useful properties.  We now use this construction to build a complete Lyapunov
function for a continuous mapping.  To do so, we ``sum over attracting
sets,'' and so this requires being able to sum in a useful way.  Thus we
first establish some countability for strong trapping regions.  To do this,
we shall make a connection between strong trapping regions and points that
can be connected by chains, rather as we did in the proof of
Theorem~\ref{the:conley-decomp}\@.

This being said, we begin with some constructions with chains.
\begin{lemma}\label{lem:robust-eps}
Let\/ $(\ts{X},\d)$ and let\/ $\phi\in\mappings[0]{\ts{X}}{\ts{X}}$\@.  Let\/
$\varepsilon\in\mappings[0]{\ts{X}}{\realp}$ and\/ $x,y\in\ts{X}$ be such
that there is no\/ $3\varepsilon$-chain of length at least\/ $2$ for\/ $\phi$
from\/ $x$ to\/ $y$\@.  Then there are neighbourhoods\/ $\nbhd{U}$ of\/ $x$
and\/ $\nbhd{V}$ of\/ $y$\@, and\/ $\delta\in\mappings[0]{\ts{X}}{\realp}$
such that there is no\/ $\delta$ chain for\/ $\phi$ from a point in\/
$\nbhd{U}$ to a point in\/ $\nbhd{V}$\@.
\begin{proof}
We first claim that there exists $\delta_1\in\mappings[0]{\ts{X}}{\realp}$
such that, if
\begin{equation*}
x_0,x_1,\dots,x_k
\end{equation*}
is a $\varepsilon$-chain for $\phi$\@, and if $\d(z,x_0)<\delta_1(x_0)$\@,
then
\begin{equation*}
(z,x_1,\dots,x_n)
\end{equation*}
is a $3\varepsilon$-chain for $\phi$\@.  By Lemma~\ref{lem:Papprox}\@, let
$\delta\in\mappings[0]{\ts{X}}{\realp}$ such that
\begin{equation*}
\d(z_1,z_2)<\delta(z_1)\implies
\tfrac{1}{2}\varepsilon(z_1)<\varepsilon(z_2)<\tfrac{3}{2}\varepsilon(z_1).
\end{equation*}
Without loss of generality, suppose that $\delta\le\frac{1}{2}\varepsilon$\@.
By Lemma~\ref{lem:general-epsilon-delta}\@, let
$\delta_1\in\mappings[0]{\ts{X}}{\realp}$ be such that
\begin{equation*}
\d(z_1,z_2)<\delta_1(z_1)\implies\d(\phi(z_1),\phi(z_2))<\delta(\phi(z_1)).
\end{equation*}
Note that
\begin{equation*}
\d(z,x_0)<\delta_1(x_0)\implies\d(\phi(x_0),\phi(z))<\delta(\phi(x_0))
\implies\frac{1}{2}\varepsilon(\phi(x_0))<\varepsilon(\phi(z)).
\end{equation*}
Then, if $\d(z,x_0)<\delta_1(x_0)$\@, we have
\begin{align*}
\d(\phi(z),x_1)\le&\;\d(\phi(z),\phi(x_0))+\d(\phi(x_0),x_1)\\
<&\;\delta(\phi(x_0))+\varepsilon(\phi(x_0))\le
\frac{3}{2}\varepsilon(\phi(x_0))\le3\varepsilon(\phi(z)),
\end{align*}
as desired.

Next we claim that there exists $\delta_2\in\mappings[0]{\ts{X}}{\realp}$
such that, if
\begin{equation*}
x_0,x_1,\dots,x_{k-1},x_k
\end{equation*}
is a $\delta_2$-chain with $k\ge2$ for $\phi$ and if
$\d(x_k,z)<\varepsilon(x_k)$\@, then
\begin{equation*}
x_0,x_1,\dots,x_{k-1},z
\end{equation*}
is a $3\varepsilon$-chain for $\phi$\@.  Indeed, by
Lemma~\ref{lem:Papprox}\@, let $\delta_2\in\mappings[0]{\ts{X}}{\realp}$ be
such that
\begin{equation*}
\d(z_1,z_2)<\delta_2(z_1)\implies\frac{1}{2}\varepsilon(z_1)<
\varepsilon(z_2)<\frac{3}{2}\varepsilon(z_1).
\end{equation*}
Without loss of generality, suppose that $\delta_2<\varepsilon$\@.  Note that
\begin{equation*}
\d(\phi(x_{k-1}),x_k)<\varepsilon(\phi(x_{k-1}))\implies
\varepsilon(x_k)<\frac{3}{2}\varepsilon(\phi(x_{k-1})).
\end{equation*}
Therefore, if $\d(x_k,z)<\varepsilon(x_k)$\@, we have
\begin{align*}
\d(\phi(x_{k-1}),z)\le&\;\d(\phi(x_{k-1}),x_k)+\d(x_k,z)\\
<&\;\delta_2(\phi(x_{k-1}))+\varepsilon(x_k)<
\frac{5}{2}\varepsilon(\phi(x_{k-1}))<3\varepsilon(\phi(x_{k-1})),
\end{align*}
as desired.

Now take $\delta=\frac{1}{2}\min\{\delta_1,\delta_2\}$ with $\delta_1$ and
$\delta_2$ as in the preceding two paragraphs.  Let $\nbhd{U}$ be a
neighbourhood of $x$ such that
\begin{compactenum}
\item $\nbhd{U}\subset\oball{\delta(x)}{x}$ and
\item $2\delta(x')>\delta(x)$ for $x'\in\nbhd{U}$\@,
\end{compactenum}
and specify a neighbourhood $\nbhd{V}$ of $y$ similarly.  Now, if
\begin{equation*}
x',x_1,\dots,x_{k-1},y'
\end{equation*}
is a $\delta$-chain for $\phi$ from $x'\in\nbhd{U}$ to $y'\in\nbhd{V}$\@,
then
\begin{equation*}
\d(x,x')<\frac{1}{2}\delta_1(x)<\delta_1(x')
\end{equation*}
and
\begin{equation*}
\d(y,y')<\frac{1}{2}\delta_2(y)<\delta_2(y'),
\end{equation*}
our constructions ensure that
\begin{equation*}
x,x_1,\dots,x_{k-1},y
\end{equation*}
is a $3\varepsilon$-chain for $\phi$ from $x$ to $y$\@.  This proves the
lemma by contraposition.
\end{proof}
\end{lemma}

Using the lemma, we devise a countable subset of strong trapping regions with
desired properties.  As we shall see\textemdash{}and as we have seen already
in our proof of Theorem~\ref{the:conley-decomp}\textemdash{}there is a
connection between strong trapping regions and chains.  Thus the first step
in reducing to a countable number of strong trapping regions is to reduce to
a countable number of $\varepsilon$'s used to define these strong trapping
regions.
\begin{lemma}\label{lem:sPdef}
Let\/ $(\ts{X},\d)$ be a separable metric space and let\/
$\phi\in\mappings[0]{\ts{X}}{\ts{X}}$\@.  Then there exists a countable
subset\/ $\sP\subset\mappings[0]{\ts{X}}{\realp}$ with the following
property: if\/ $x,y\in\ts{X}$ and if, for some\/
$\varepsilon\in\mappings[0]{\ts{X}}{\realp}$\@, there is no\/
$\varepsilon$-chain of length at least\/ $2$ for $\phi$ from\/ $x$ to\/
$y$\@, then there exists\/ $\delta\in\sP$ such that there is no\/
$\delta$-chain of length at least\/ $2$ for $\phi$ from\/ $x$ to\/ $y$\@.
\begin{proof}
Let us denote
\begin{multline*}
\sNC=\{(x,y)\in\ts{X}\times\ts{X}|\enspace\textrm{there exists\/ $\varepsilon\in\mappings[0]{\ts{X}}{\realp}$ such
that}\\\textrm{there is no\/ $\varepsilon$-chain of length at least $2$ from
$x$ to $y$}\}.
\end{multline*}
For $(x,y)\in\sNC$\@, there exist neighbourhoods $\nbhd{U}_x$ of $x$ and
$\nbhd{V}_y$ of $y$\@, and $\delta_{x,y}\in\mappings[0]{\ts{X}}{\realp}$
such that there is no $\delta_{x,j}$-chain from a point in $\nbhd{U}$ to a
point in $\nbhd{V}$\@.  Now note that the collection
$\nbhd{U}_x\times\nbhd{V}_y$\@, $(x,y)\in\sNC$\@, of open sets covers
$\sNC$\@.  There is a countable collection of points
$\ifam{(x_j,y_j)}_{j\in\integerp}$ from $\sNC$ such that the sets
$\nbhd{N}_j\eqdef\nbhd{U}_{x_j}\times\nbhd{V}_{y_j}$\@, $j\in\integerp$\@,
covers $\sNC$~\cite[Theorem~16.9]{SW:70}\@.  Taking
$\sP=\setdef{\delta_j\eqdef\delta_{x_j,y_j}}{j\in\integerp}$ gives the
result.
\end{proof}
\end{lemma}

Let us extract from the proof the notation
\begin{multline*}
\sNC=\{(x,y)\in\ts{X}\times\ts{X}|\enspace\textrm{there exists\/ $\varepsilon\in\mappings[0]{\ts{X}}{\realp}$ such
that}\\\textrm{ there is no\/ $\varepsilon$-chain of length at least $2$ from
$x$ to $y$}\}.
\end{multline*}
Let $D\subset\ts{X}$ be a countable dense subset and let $\sP$ be as
prescribed by the lemma.  Now, for $x\in\ts{X}$ and
$\varepsilon\in\mappings[0]{\ts{X}}{\realp}$\@, denote
\begin{equation*}
\nbhd{T}(x,\varepsilon)=\setdef{y\in\ts{X}}
{\textrm{there exists an $\varepsilon$-chain of length at least $2$ from $x$
to $y$}}.
\end{equation*}
By Lemma~\ref{lem:chain->trap}\@, $\nbhd{T}(x,\varepsilon)$ is a nonempty,
open strong trapping region.  Denote
\begin{equation*}
\sT=\setdef{\nbhd{T}(x,\varepsilon)}{x\in D,\ \varepsilon\in\sP}.
\end{equation*}
As $\sT$ is countable, we enumerate its elements as
$\ifam{\nbhd{T}_j}_{j\in\integerp}$\@.  We let $A_j$ be the attracting set
for $\nbhd{T}_j$ and we let $L_j=L_{\nbhd{T}_j}$ be the Lyapunov function for
$\nbhd{T}_j$ as constructed in the preceding section.  Then define
\begin{equation}\label{eq:clf-map-def}
L(x)=\sum_{j=1}^\infty\frac{2L_j(x)}{3^j}.
\end{equation}
Since $L_j$ takes values in $\interval[0,1]$\@, this series converges
uniformly by the Weierstrass $M$-test, and so converges uniformly to a
continuous function.  Certainly $L\scirc\phi(x)\le L(x)$ for
$x\in\ts{X}$\@.  The following lemma lists some other pertinent properties
of $L$\@.
\begin{lemma}
Let\/ $(\ts{X},\d)$ be a separable metric space and let\/
$\phi\in\mappings[0]{\ts{X}}{\ts{X}}$\@.  With the notation introduced above,
the following statements hold for\/ $(x,y)\in\sNC$\@:
\begin{compactenum}[(i)]
\item \label{pl:clf1} there exist\/ $z\in D$ and\/ $\delta_j\in\sP$ such
that\/ $y\not\in\nbhd{T}(z,\delta_j)$ and\/
$x\in\Orb^-(\nbhd{T}(z,\delta_j))$\@.\savenum
\end{compactenum}
With\/ $z\in D$ and\/ $j\in\integerp$ as in the preceding statement, let\/
$k\in\integerp$ be such that\/ $\nbhd{T}(z,\delta_j)=\nbhd{T}_k$\@.  Then the
following statements hold;
\begin{compactenum}[(i)]\resumenum
\item \label{pl:clf2} if\/ $y\in\ChRec(\phi)$\@, then\/ $L_k(x)<L_k(y)$\@;
\item \label{pl:clf3} if\/ $y=x$\@, then\/ $L_k\scirc\phi(x)<L_k(x)$\@;
\item \label{pl:clf4} if\/ $x\not\in\ChRec(\phi)$\@, then\/ $L\scirc\phi(x)<L(x)$\@;
\item \label{pl:clf5} if\/ $x,y\in\ChRec(\phi)$ but\/ $x$ and\/ $y$ are not
chain equivalent, then\/ $L(x)\not=L(y)$\@;
\item \label{pl:clf6} if\/ $C,C'\subset\ChRec(\phi)$ are distinct chain
components and if, for each\/
$\varepsilon\in\mappings[0]{\ts{X}}{\realp}$\@, there is a\/
$\varepsilon$-chain for\/ $\phi$ from a point in\/ $C$ to a point in\/
$C'$\@, then\/ $L(C)>L(C')$\@.
\end{compactenum}
\begin{proof}
\eqref{pl:clf1} Adopting the notation from the proof of
Lemma~\ref{lem:sPdef}\@, let $\nbhd{N}_j=\nbhd{U}_{x_j}\times\nbhd{V}_{y_j}$
be such that $(x,y)\in\nbhd{N}_j$\@.  Note that $D\cap\nbhd{U}_{x_j}$ is
dense in $\nbhd{U}_{x_j}$\@.  For $z\in D\cap\nbhd{U}_{x_j}$\@, there is no
$\delta_j$-chain from $z$ to $y$\@, as we can conclude from
Lemma~\ref{lem:sPdef}\@.  Thus $y\not\in\nbhd{T}(z,\delta_j)$\@.  Now, if $z$
is sufficiently close to $x$\@, then $\phi(z)$ will be close enough to
$\phi(x)$ that
\begin{equation*}
z,\phi(x),\phi^2(x)
\end{equation*}
is a $\delta_j$-chain from $z$ to $\phi^2(x)$\@.  Thus
$\phi^2(x)\in\nbhd{T}(z,\delta_j)$ and so $x\in\cup_{j\in\integerp}\phi^{-j}(\nbhd{T}(z,\varepsilon))$\@.

\eqref{pl:clf2} If $y\in\ChRec(\phi)$\@, then
$y\in A_k\cup(\ts{X}\setminus\Orb^-(\nbhd{T}_k))$ by
Theorem~\ref{the:conley-decomp}\@.  Since $y\not\in\nbhd{T}_k\supset A_k$\@,
we must have $y\in\ts{X}\setminus\Orb^-(\nbhd{T}_k)$\@.  Since
$x\in\Orb^-(\nbhd{T}_k)$\@, $L_k(x)<1$\@,~\cf~the proof of
Lemma~\ref{lem:LT}\eqref{pl:LT2}\@, which gives this part of the result.

\eqref{pl:clf3} With the stated hypotheses, we have
$x\in\Orb(\nbhd{T}_k)\setminus\nbhd{T}_k$\@, whereupon
$L_k\scirc\phi(x)<L_k(x)$ by Lemma~\ref{lem:LT}\eqref{pl:LT3}\@.

\eqref{pl:clf4} Note that $x\in\ChRec(\phi)$ if and only if $(x,x)\in\sNC$\@.
This being the case, this part of the result follows from the previous one.

\eqref{pl:clf5} Note that $x$ and $y$ are not chain equivalent if and only if
$(x,y)\in\sNC$ or $(y,x)\in\sNC$\@.  Let us consider the case
$(x,y)\in\sNC$\@.  By part~\eqref{pl:clf2}\@, $L_k(x)<L_y(y)$\@.  By
Lemma~\ref{lem:LT}\eqref{pl:LT2}\@, $L_k(x)=0$ and $L_k(y)=1$\@.  Thus $L(x)$
and $L(y)$ do not agree since they necessarily have different ternary (that
is, base $3$) expansions.

\eqref{pl:clf6} Note that
\begin{equation*}
C,C'\subset\ChRec(\phi)=\bigcap_{\nbhd{T}\in\sT}
\setdef{A\cup(\ts{X}\setminus\Orb^-(\nbhd{T}))}{\textrm{$A$ is the attracting
set for\/ $\nbhd{T}$}}
\end{equation*}
by Theorem~\ref{the:conley-decomp}\@.  Suppose that $C\subset A_k$ for some
$k\in\integerp$\@.  We claim that, with the hypotheses of this part of the
lemma, $C'\subset A_k$\@.  Indeed, if $C'\not\subset A_k$\@, then
$C_k\subset\ts{X}\setminus\Orb^-(\nbhd{T}_k)$\@.  However, by
Lemma~\ref{lem:trap-erf}\@, there exists
$\varepsilon\in\mappings[0]{\ts{X}}{\realp}$ such that every
$\varepsilon$-chain starting in $\nbhd{T}_k$ ends in $\nbhd{T}_k$\@, in
contradiction with the current hypotheses.  Thus, indeed, $C'\subset A_k$\@.
This shows that, if $L_k(C)=0$\@, then $L_k(C')=0$\@.  Since, for each
$k\in\integerp$\@, $L_k(C),L_k(C')\in\{0,1\}$ by
Lemma~\ref{lem:LT}\eqref{pl:LT2}\@, we have $L_k(C)\ge L_k(C')$ for each
$k\in\integerp$\@.  This part of the result now follows from the previous
part.
\end{proof}
\end{lemma}

This gives the following theorem which is the second part of the Fundamental
Theorem of Dynamical Systems for discrete-time flows and semiflows.
\begin{theorem}\label{the:clf-maps}
Let\/ $(\ts{X},\d)$ be a separable metric space and let\/
$\phi\in\mappings[0]{\ts{X}}{\ts{X}}$\@.  Then there exists a complete
Lyapunov function for\/ $\phi$\@.
\end{theorem}

\subsection{The continuous-time case: the complete Lyapunov function}

Now we use the preceding results concerning mappings to obtain the existence
of complete Lyapunov functions for flows.  We use the idea of
\citet{MP:11}\@.  However, the proof of \citeauthor{MP:11} contains a number
of errors, including making use of connectedness of chain components (as far
as we know, this has only been proved in the compact case) and using
Lemma~\ref{lem:comp-inv} for semiflows (the lemma is not true for semiflows).

The theorem we prove is the following.
\begin{theorem}\label{the:clf-flows}
Let\/ $(\ts{X},\d)$ be a metric space and let\/ $\Phi$ be a continuous-time
topological flow or semiflow on\/ $\ts{X}$\@.  If\/
$\ell\in\mappings[0]{\ts{X}}{\interval[0,1]}$ is a complete Lyapunov function
for\/ $\Phidisc{1}$\@, then the function\/ $\map{L}{\ts{X}}{\real}$ given by
\begin{equation*}
L(x)=\int_0^1\ell\scirc\Phi(s,x)\,\d{s}
\end{equation*}
is a complete Lyapunov function for\/ $\Phi$\@.  In particular, if\/ $\ts{X}$
is separable, then there exists a complete Lyapunov function for\/ $\Phi$\@.
\begin{proof}
We note that the complete Lyapunov function $\ell$ that we constructed in the
proof of Theorem~\ref{the:clf-maps} has a property of which we shall make
use.  Namely, in the definition~\eqref{eq:clf-map-def} of $\ell$\@, we see
that points in $\ChRec(\Phidisc{1})$\@, being points where the functions
$L_j$ take value $0$\@, are points whose ternary (\ie~base~$3$) expansion
contains only $1$s and $2$s.  This means that $\ell(\ChRec(\Phidisc{1}))$ is
a subset of the classical middle-thirds Cantor set, which is closed and
nowhere dense.

First note that $L$ is continuous by standard results concerning swapping
limits and integrals~\cite[\eg][Theorem~7.16]{WR:76}.  Since $\ell(x)\in\interval[0,1]$ for
every $x\in\ts{X}$\@, we also have $L(x)\in\interval[0,1]$ for every
$x\in\ts{X}$\@.

The following calculation will be useful in the remainder of the proof.  Let
$t\in\real$ (if $\Phi$ is a flow) or $t\in\realnn$ (if $\Phi$ is a semiflow)
and suppose that $t\in\interval[j,{j+1})$ for some $j\in\integer$\@.  Let
$x\in\ts{X}$\@.  Then we have
\begin{align}\notag
L\scirc\Phi(t,x)=&\;
\int_0^1\ell\scirc\Phi(s,\Phi(t-j+j,x))\,\d{s}\\\notag
=&\;\int_0^1\ell\scirc\Phi(s+t-j,\Phi(j,x))\,\d{s}\\\notag
=&\;\int_{t-j}^{1+t-j}\ell\scirc\Phi(s,\Phi(j,x))\,\d{s}\\\notag
=&\;\int_{t-j}^1\ell\scirc\Phi(s,\Phi(j,x))\,\d{s}+
\int_1^{1+t-j}\ell\scirc\Phi(s,\Phi(j,x))\,\d{s}\\\label{eq:dt->ctclt1}
=&\;\int_{t-j}^1\ell\scirc\Phi(s,\Phi(j,x))\,\d{s}
+\int_0^{t-j}\ell\scirc\Phi(1+s,\Phi(j,x))\,\d{s}.
\end{align}

Let $x\in\ts{X}$ and let $t\in\interval[0,1)$\@.  Then
\begin{equation*}
\ell\scirc\Phi(1+s,x)\le\ell\scirc\Phi(s,x),\qquad s\in\realnn,
\end{equation*}
since $\ell$ is a complete Lyapunov function for $\Phidisc{1}$\@.  Then,
using~\eqref{eq:dt->ctclt1} for $j=0$\@, we have
\begin{align*}
L\scirc\Phi(t,x)=&\;\int_t^1\ell\scirc\Phi(s,x)\,\d{s}
+\int_0^t\ell\scirc\Phi(1+s,x)\,\d{s}\\
\le&\;\int_t^1\ell\scirc\Phi(s,x)\,\d{s}
+\int_0^t\ell\scirc\Phi(s,x)\,\d{s}\\
=&\;\int_0^1\Phi(s,x)\,\d{s}=L(x).
\end{align*}
In particular, $L\scirc\Phi(1,x)\le L(x)$\@.  Now, if $t\in\realp$ satisfies
$t\in\interval[j,{j+1})$ for some $j\in\integernn$\@, then
\begin{equation*}
L\scirc\Phi(t,x)=L\scirc\Phi(t-j,\Phi(j,x))\le
L\scirc\Phi(j,x),
\end{equation*}
and so, by an elementary induction, $L\scirc\Phi(t,x)\le L(x)$\@.

Let $t\in\real$ (if $\Phi$ is a flow) or $t\in\realnn$ (if $\Phi$ is a
semiflow) satisfy $t\in\interval[j,{j+1})$ for some $j\in\integer$\@, let
$x\in\ChRec(\Phi)=\ChRec(\Phidisc{1})$\@, this last equality by
Corollary~\ref{cor:chan-recurrent}\@.  By
Proposition~\ref{prop:ChRec-props}\eqref{pl:ChRec1} and
Corollary~\ref{cor:chan-recurrent} we have
\begin{multline*}
\Phi(s,\Phi(k,x))\in\ChRec(\Phi)=\ChRec(\Phidisc{1}),\\
k\in\integer\ \textrm{(for flows) or}\ k\in\integernn\ \textrm{(for semiflows)},\ s\in\interval[0,1).
\end{multline*}
Therefore, by Theorem~\ref{the:clf-maps} we have
\begin{equation*}
\ell\scirc\Phi(1+s,\Phi(j,x))=
\ell\scirc\Phi(s,\Phi(j,x)),\qquad s\in\interval[0,1).
\end{equation*}
Thus we have, by~\eqref{eq:dt->ctclt1}\@,
\begin{align*}
L\scirc\Phi(t,x)=&\;
\int_{t-j}^1\ell\scirc\Phi(s,\Phi(j,x)))\,\d{s}+
\int_0^{t-j}\ell\scirc\Phi(1+s,\Phi(j,x))\,\d{s}\\
=&\;\int_0^1\ell\scirc\Phi(s,\Phi(j,x))\,\d{s}
=\int_0^1\ell\scirc\Phi(s,x)\,\d{s}=L(x).
\end{align*}
Thus $L$ is constant along forward trajectories through points in
$\ChRec(\Phi)$\@.

Let $x\in\ts{X}\setminus\ChRec(\Phi)$\@.  Since $\ts{X}\setminus\ChRec(\Phi)$
is open by Proposition~\ref{prop:ChRec-props}\eqref{pl:ChRec2}\@, there
exists $\tau\in\realp$ such that
\begin{equation*}
\Phi(t,x)\in\ts{X}\setminus\ChRec(\Phi)=
\ts{X}\setminus\ChRec(\Phidisc{1}),\qquad t\in\interval[{t_1},{t_1+\tau}].
\end{equation*}
Therefore,
\begin{equation*}
\ell\scirc\Phi(1+t,x)<\ell(\Phi(t,x)),\qquad t\in\interval[0,\tau]
\end{equation*}
since $\ell$ is a complete Lyapunov function for $\Phidisc{1}$\@.  Therefore,
for $t\in\interval[0,\tau]$\@, we can use~\eqref{eq:dt->ctclt1} to get
\begin{align*}
L\scirc\Phi(t,x)=&\;\int_t^1\ell\scirc\Phi(s,x)\,\d{s}+
\int_0^t\ell\scirc\Phi(1+s,x)\,\d{s}\\
<&\;\int_t^1\ell\scirc\Phi(s,x)\,\d{s}+
\int_0^t\ell\scirc\Phi(s,x)\,\d{s}=L(x).
\end{align*}
If $t>\tau$\@, we have
\begin{equation*}
L\scirc\Phi(t,x)=L\scirc\Phi(t-\tau,\Phi(\tau,x))\le
L\scirc\Phi(\tau,x)<L(x)
\end{equation*}
since $\Phi(t-\tau,x)\in\ts{X}\setminus\ChRec(\Phi)$\@.

Now let $\alpha\in\ell(\ChRec(\Phidisc{1}))$\@.  Since $\ell$ is a complete
Lyapunov function for $\Phidisc{1}$\@, $\ell^{-1}(\alpha)$ is a chain
component for $\Phidisc{1}$\@, and so also a chain component of $\Phi$ by
Theorem~\ref{the:chain-equiv}\@.  By
Proposition~\ref{prop:chain-comp}\eqref{pl:chain-comp2}\@,
$\ell^{-1}(\alpha)$ is invariant under $\Phi$\@.  Thus
$\Phi(t,x)\in\ell^{-1}(\alpha)$ for all $t\in\real$ (if $\Phi$ is a flow) or
for all $t\in\realnn$ (if $\Phi$ is a semiflow) if $x\in\ell^{-1}(\alpha)$\@.
Therefore, for $x\in\ell^{-1}(\alpha)$\@, we have
\begin{equation*}
L(x)=\int_0^1\ell\scirc\Phi(t,x)\,\d{t}=\alpha,
\end{equation*}
and so
\begin{equation}\label{eq:dt->ctclt2}
\ell^{-1}(\alpha)\subset L^{-1}(\alpha),
\qquad\alpha\in\ell(\ChRec(\Phidisc{1})).
\end{equation}
To show that this inclusion is equality, we claim that it is sufficient to
show that $L^{-1}(\alpha)\subset\ChRec(\Phi)$\@.  To prove this claim, we
proceed by contradiction.  Thus suppose that
\begin{compactenum}
\item \label{enum:ctclf1} there exists $\alpha\in\ell(\ChRec(\Phidisc{1}))$
such that $\ell^{-1}(\alpha)\subsetneq L^{-1}(\alpha)$ and
\item \label{enum:ctclf2} $L^{-1}(\alpha)\subset\ChRec(\Phi)$ for each
$\alpha\in\ell(\ChRec(\Phidisc{1}))$\@.
\end{compactenum}
Since chain equivalence defines an equivalence relation on
$\ChRec(\Phi)=\ChRec(\Phidisc{1})$\@, and since the equivalence classes are
determined by the unique value that $\ell$ takes on each equivalence class,
the assumptions~\ref{enum:ctclf1} and~\ref{enum:ctclf2} imply that there
exists $\alpha'\in\ell(\ChRec(\Phidisc{1}))$ such that $\alpha\not=\alpha'$
and such that $L^{-1}(\alpha)\cap\ell^{-1}(\alpha')\not=\emptyset$\@.  Thus
there exists $x\in\ChRec(\Phidisc{1})$ such that $L(x)=\alpha$ and
$\ell(x)=\alpha'$\@.  However, since
$\ell^{-1}(\alpha')\subset L^{-1}(\alpha')$\@, this implies that
$L(x)=\alpha'\not=\alpha$\@.  This contradiction shows that, if we can show
that $L^{-1}(\alpha)\subset\ChRec(\Phi)$ for every
$\alpha\in\ell(\ChRec(\Phidisc{1}))$\@, it will follow that
$\ell^{-1}(\alpha)=L^{-1}(\alpha)$ for every
$\alpha\in\ell(\ChRec(\Phidisc{1}))$\@.

To show that $L^{-1}(\alpha)\subset\ChRec(\Phi)$ for every
$\alpha\in\ell(\ChRec(\Phidisc{1})$\@, let
$x\in\ts{X}\setminus\ChRec(\Phi)$\@.  We claim that
$L(x)\not\in\ell(\ChRec(\Phidisc{1}))$\@.  Let
\begin{equation*}
\tau_x=\sup\setdef{t\in\real}{\Phi(t,x)\in\ts{X}\setminus\ChRec(\Phi)},
\end{equation*}
noting that $\tau_x\in\realp$ since $\ts{X}\setminus\ChRec(\Phi)$ is open by
Proposition~\ref{prop:ChRec-props}\eqref{pl:ChRec2}\@.  If $\tau_x<\infty$\@,
we claim that $\Phi(\tau_x,x)\in\ChRec(\Phi)$\@.  Indeed, if
$\Phi(\tau_x.x)\not\in\ChRec(\Phi)$\@, then there must be $\tau'>\tau_x$ with
$\Phi(t,x)\not\in\ChRec(\Phi)$ for $t\in\interval[{\tau_x},{\tau'}]$ by
openness of $\ts{X}\setminus\ChRec(\Phi)$\@, contradicting the definition of
$\tau_x$\@.  Let
$\alpha=\ell\scirc\Phi(\tau_x,x)\in\ell(\ChRec(\Phidisc{1}))$\@.  Thus,
using~\eqref{eq:dt->ctclt2}\@,
\begin{equation*}
\Phi(\tau_x,x)\in\ell^{-1}(\alpha)\subset L^{-1}(\alpha),
\end{equation*}
and so $L\scirc\Phi(\tau_x,x)=\alpha$\@.  By the Mean Value Theorem for
integrals~\cite[Theorem~5.10]{WR:76}\@, there exists $s_0\in\interval(0,1)$
such that
\begin{equation}\label{eq:dt->ctclt3}
L(x)=\int_0^1\ell\scirc\Phi(s,x)\,\d{s}=\ell\scirc\Phi(s_0,x)(1-0)=
\ell\scirc\Phi(s_0,x).
\end{equation}
If $\tau_x\ge1$\@,~\eqref{eq:dt->ctclt3} immediately gives
$L(x)\not\in\ell(\ChRec(\Phidisc{1}))$ since the preceding formula holds with
$s_0<1<\tau_x$ and using the definition of $\tau_x$\@.  On the other hand, if
$\tau_x\in\interval(0,1)$\@, we proceed as follows.  Since
$x\in\ts{X}\setminus\ChRec(\Phi)$\@, we have
\begin{equation*}
L(x)>L\scirc\Phi(\tau_x,x)
\end{equation*}
using already proved properties of $L$\@.  By~\eqref{eq:dt->ctclt3} we thus
have, for some $s_0\in\interval(0,1)$\@,
\begin{equation*}
L(x)=\ell\scirc\Phi(s_0,x)>\alpha=\ell\scirc\Phi(\tau_x,x).
\end{equation*}
We claim that $s_0\in\interval(0,{\tau_x})$\@.  To see this, first note that
\begin{equation}\label{eq:dt->ctclt4}
L\scirc\Phi(t,x)\le L\scirc\Phi(\tau_x,x)=\alpha,
\qquad t\in\interval[{\tau_x},1),
\end{equation}
using the fact, already proved, that $L$ is nonincreasing along trajectories.
Suppose now that there exists $t_0\in\interval[{\tau_x},1)$ such that
$\ell\scirc\Phi(t_0,x)>\alpha$\@.  By
Proposition~\ref{prop:ChRec-props}\eqref{pl:ChRec1} and the fact that
$\Phi(\tau_x,x)\in\ChRec(\Phi)$\@, we have
\begin{equation*}
\Phi(t_0,x)\in\ChRec(\Phi)=\ChRec(\Phidisc{1})\implies
\alpha'\eqdef\ell\scirc\Phi(t_0,x)\in\ell(\ChRec(\Phidisc{1})).
\end{equation*}
Therefore, by~\eqref{eq:dt->ctclt2}\@, we have
\begin{equation*}
\ell\scirc\Phi(t_0,x)=\alpha'\implies
\Phi(x,t_0)\in\ell^{-1}(\alpha')\subset L^{-1}(\alpha')\implies
L\scirc\Phi(t_0,x)=\alpha'>\alpha,
\end{equation*}
in contradiction with~\eqref{eq:dt->ctclt4}\@.  Therefore, we have shown
that, for any $x\in\ts{X}\setminus\ChRec(\Phi)$ and irregardless of the value
of $\tau_x\in\realp$\@, there exists $s_0\in\interval(0,{\tau_x})$ such that
\begin{equation*}
L(x)=\ell\scirc\Phi(s_0,x)\not\in\ChRec(\Phidisc{1}).
\end{equation*}

Now we have the following logical implications:
\begin{align*}
&\bigl(x\in\ts{X}\setminus\ChRec(\Phi)\implies
L(x)\not\in\ell(\ChRec(\Phidisc{1}))\bigr)\\
\iff&\bigl(L(x)\in\ell(\ChRec(\Phidisc{1}))\implies
x\in\ChRec(\Phi)\bigr)\\
\iff&\bigl(L^{-1}(\alpha)\subset\ChRec(\Phi),\ \alpha\in\ell(\ChRec(\Phidisc{1}))\bigr)\\
\implies&\ell^{-1}(\alpha)=L^{-1}(\alpha),\
\alpha\in\ell(\ChRec(\Phidisc{1})),
\end{align*}
where the last implication was proved above.  This shows that the partition
of $\ChRec(\Phi)=\ChRec(\Phidisc{1})$ by level sets $\ell^{-1}(\alpha)$\@,
$\alpha\in\ell(\ChRec(\Phidisc{1}))$\@, is the same as the partition of
$\ChRec(\Phi)$ by the level sets $L^{-1}(\alpha)$\@,
$\alpha\in\ell(\ChRec(\Phidisc{1}))$\@.  This shows that $L(C)=L(C')$ for
chain components $C$ and $C'$ if and only if $C=C'$\@.
\end{proof}
\end{theorem}

\printbibliography[heading=bibintoc]

\end{document}